\documentclass[10pt,a4paper]{article}

\usepackage[T1]{fontenc}
\usepackage[english]{babel}
\usepackage[a4paper,margin=2cm,headheight=14pt,headsep=14pt,footskip=25pt]{geometry}
\usepackage{times}
\usepackage{titlesec}
\usepackage{fancyhdr}
\usepackage{tikz}
\usetikzlibrary{arrows.meta,calc,positioning,decorations.pathreplacing}
\usepackage[font=small,labelfont=bf,labelsep=period]{caption}
\usepackage{microtype}
\usepackage{amsmath,amssymb,amsthm,mathtools,calrsfs}
\usepackage{csquotes}
\usepackage{enumitem}
\usepackage{xcolor}
\usepackage{cite}
\usepackage{hyperref}
\usepackage{xurl}
\ifdefined\newcounteralias
\else
  \usepackage{aliascnt}
\fi
\usepackage[nameinlink,capitalize,noabbrev]{cleveref}

\hypersetup{
  colorlinks=true,
  linkcolor=magenta,
  citecolor=green,
  urlcolor=blue,
  pdftitle={Singular Tracking of Fractional Brownian Motion under Proportional Costs},
  pdfauthor={Alexandre Autran},
  pdfkeywords={fractional Brownian motion, singular control, causal tracking, total variation, past information, filtered tangent, Volterra homogenization, Markovian lift}
}

\makeatletter
\renewcommand\maketitle{\begingroup
  \parindent=0pt
  \vspace*{-20pt}{\fontsize{18}{20}\selectfont\scshape\@title\par}\vspace{20pt}{\normalfont\normalsize\@author\par}\vspace{24pt}\endgroup}
\makeatother
\renewenvironment{abstract}
  {\begin{center}\large\scshape Abstract\end{center}\begin{list}{}{\leftmargin=0.35in\rightmargin=0.15in}\item\relax\normalfont\normalsize}
  {\end{list}\vspace{7pt}}
\titleformat{\section}{\large\normalfont\scshape}{\thesection}{1em}{}
\titleformat{\subsection}{\normalsize\normalfont\scshape}{\thesubsection}{1em}{}
\titleformat{\paragraph}[runin]{\normalsize\bfseries}{}{0pt}{}
\titlespacing*{\section}{0pt}{14pt plus 3pt minus 2pt}{8pt plus 1pt}
\titlespacing*{\subsection}{0pt}{10pt plus 2pt minus 1pt}{6pt plus 1pt}
\titlespacing*{\paragraph}{0pt}{6pt plus 1pt}{0.55em}
\setlist{nosep,leftmargin=*}
\allowdisplaybreaks[2]
\newcommand{\appdef}[5]{#1\colon\left|\begin{array}{rcl}
  #2&\longrightarrow&#3\\[2pt]
  #4&\longmapsto&\displaystyle #5
  \end{array}\right.}

\newtheorem{theorem}{Theorem}[section]
\ifdefined\newcounteralias
  \newcommand{\newsharedtheorem}[2]{\newtheorem{#1}[theorem]{#2}}
\else
  \newcommand{\newsharedtheorem}[2]{\newaliascnt{#1}{theorem}\newtheorem{#1}[#1]{#2}\aliascntresetthe{#1}}
\fi
\newsharedtheorem{proposition}{Proposition}
\newsharedtheorem{lemma}{Lemma}
\newsharedtheorem{corollary}{Corollary}
\theoremstyle{definition}
\newsharedtheorem{definition}{Definition}
\newsharedtheorem{assumption}{Assumption}
\theoremstyle{remark}
\newsharedtheorem{remark}{Remark}

\crefname{theorem}{theorem}{theorems}
\Crefname{theorem}{Theorem}{Theorems}
\crefname{proposition}{proposition}{propositions}
\Crefname{proposition}{Proposition}{Propositions}
\crefname{lemma}{lemma}{lemmas}
\Crefname{lemma}{Lemma}{Lemmas}
\crefname{corollary}{corollary}{corollaries}
\Crefname{corollary}{Corollary}{Corollaries}
\crefname{definition}{definition}{definitions}
\Crefname{definition}{Definition}{Definitions}
\crefname{assumption}{assumption}{assumptions}
\Crefname{assumption}{Assumption}{Assumptions}
\crefname{remark}{remark}{remarks}
\Crefname{remark}{Remark}{Remarks}
\crefname{section}{section}{sections}
\Crefname{section}{Section}{Sections}
\crefname{equation}{equation}{equations}
\Crefname{equation}{Equation}{Equations}
\crefname{figure}{figure}{figures}
\Crefname{figure}{Figure}{Figures}
\crefname{table}{table}{tables}
\Crefname{table}{Table}{Tables}

\DeclareMathOperator{\TV}{TV}
\DeclareMathOperator{\Var}{Var}
\DeclareMathOperator{\Cov}{Cov}
\DeclareMathOperator{\sign}{sign}
\DeclareMathOperator{\Span}{span}
\DeclareMathOperator{\diag}{diag}
\newcommand{\CTV}{\TV^{\mathrm{aug}}}
\newcommand{\E}{\mathbf E}
\newcommand{\Pp}{\mathbf P}
\newcommand{\mR}{\mathbf R}

\newcommand{\cA}{\mathcal A}

\newcommand{\cE}{\mathcal E}
\newcommand{\cW}{\mathcal W}
\newcommand{\dd}{\,\mathrm d}

\newcommand{\ind}{\mathbf 1}
\newcommand{\past}{\mathrm{past}}
\newcommand{\loc}{\mathrm{loc}}
\newcommand{\Borel}{\mathcal B}
\newcommand{\Nn}{\mathbf N}
\newcommand{\Qq}{\mathbf Q}
\newcommand{\indep}{\mathrel{\perp\!\!\!\perp}}

\title{Singular Tracking of Fractional Brownian\\ Motion under Proportional Costs}
\author{\textbf{Alexandre Autran}\\
Department of Mathematics
\\
École normale supérieure de Rennes, France\\
\href{mailto:alexandre.autran@ens-rennes.fr}{\texttt{alexandre.autran@ens-rennes.fr}}
}
\date{}

\begin{document}
\maketitle
\begin{abstract}
We study the tracking of fractional Brownian motion by adapted controls
of finite variation, with quadratic loss and proportional adjustment
costs. Let $B^H$ be standard fractional Brownian motion with Hurst
parameter $H\in(0,1)$. We write $V_\lambda(T;a,b,\sigma)$ for the
optimal expected cost of tracking $\sigma B^H$ over a horizon $T>0$,
where $\sigma>0$ is the target amplitude, $a>0$ weights the quadratic
loss, and $\lambda b$ is the cost per unit of adjustment, with
$\lambda,b>0$. We prove
\[
 V_\lambda(T;a,b,\sigma)\underset{\lambda\to0}{\sim}\kappa_HT
 a^{(1-H)/(1+H)}(\lambda b)^{2H/(1+H)}\sigma^{2/(1+H)}.
\]
Here $\kappa_H>0$ is
the long-run optimal cost per unit time in the normalized problem
$a=b=\sigma=\lambda=1$. Observing the entire past does not change $\kappa_H$. For the normalized
problem with full past information, a horizon $L>0$ and a terminal
adjustment to the target, the difference between its value and
$\kappa_H L$ is bounded uniformly in $L$. This estimate also yields
continuity of $\kappa_H$ in $H$. In the Brownian case,
$\kappa_{1/2}=\frac12(\frac32)^{2/3}$. We extend the asymptotic formula to nonstationary targets under a
coupling condition on local paths and their observation filtrations.
We verify this condition for Gaussian Volterra targets under Brownian
observation. All asymptotically optimal controls have the same
leading-order split between tracking and adjustment costs. We also
give primal and dual bounds on $\kappa_H$. For $H<1/2$, we construct
approximations based on exponential sums with positive weights and
quantify their error as the horizon grows.
\end{abstract}


\section{Introduction}
\label{sec:introduction}

Fix a time horizon $T>0$ and an observation filtration
$\mathbf F=(\mathcal F_t)_{0\le t\le T}$, where $\mathcal F_t$
represents the information available at time $t$. A controller tracks
a continuous adapted target $\Theta=(\Theta_t)_{0\le t\le T}$ by
choosing a control $\phi=(\phi_t)_{0\le t\le T}$. \emph{Adaptedness} means
that the target and control at time $t$ are measurable with respect
to $\mathcal F_t$. The positive deterministic weights $a_t$ and $b_t$ specify the relative
importance of tracking error and adjustments. The parameter $\lambda>0$
sets the overall adjustment-cost level. Writing $\E$ for expectation,
we minimize
\begin{equation}
 J_\lambda(\phi;\Theta)
 :=\E\left[\frac12\int_0^T a_t(\phi_t-\Theta_t)^2\dd t
       +\lambda\int_{[0,T]}b_t\dd|\phi|_t\right].
 \label{eq:intro-objective}
\end{equation}
We minimize over adapted, c\`adl\`ag controls of finite variation and
write
\[
 V_\lambda(\Theta):=\inf_\phi J_\lambda(\phi;\Theta).
\]
Here \emph{c\`adl\`ag} means right-continuous with left limits, and \emph{finite
variation} means that the total amount of adjustment is finite.
The notation $\dd|\phi|_t$ denotes the total-variation measure of
the control increments. It includes the initial adjustment from the
prescribed pre-initial position $\phi_{0-}=\Theta_0$, with cost
$\lambda b_0|\phi_0-\Theta_0|$.
\Cref{sec:model-rate} gives the precise admissible class, integrability
requirements and coefficient assumptions. Suppose that target increments over a small time interval of length
$h>0$ have typical size $h^H$, where $H\in(0,1)$ describes their
scaling. An error tolerance $\delta>0$ then suggests an adjustment
interval $h\asymp\delta^{1/H}$. We use $\asymp$ for comparison up to
positive multiplicative constants independent of the small parameter.
The tracking cost per unit time is then of order $\delta^2$, whereas
the variation-cost rate satisfies
\[
 \lambda\frac\delta h\asymp\lambda\delta^{1-1/H}.
\]
Balancing these two rates suggests a tracking-error scale
$\delta_\lambda$ and a value of order
\begin{equation}
 \delta_\lambda\asymp\lambda^{H/(1+H)},\qquad
 V_\lambda(\Theta)\asymp\lambda^{2H/(1+H)}.
 \label{eq:heuristic}
\end{equation}
An increment--energy inequality proves the lower bound even for
\emph{noncausal controls}, which have access to future information. Foresight
therefore leaves the exponent unchanged. For the fractional Brownian target, let $B^H$ be standard fractional
Brownian motion with Hurst parameter $H$, normalized by $B_0^H=0$ and
\[
 \E\left[(B_{t+h}^H-B_t^H)^2\right]=h^{2H},
 \qquad t\ge0,\quad h>0.
\]
The constant $\sigma>0$ scales the target amplitude. For
$\Theta=\sigma B^H$ and constant weights $a_t=a>0$ and $b_t=b>0$,
we write
\[
 V_\lambda(T;a,b,\sigma):=V_\lambda(\sigma B^H),
\]
using the usual natural filtration of $B^H$, so that controls use
only observations of the target from time zero up to the current time.
The natural time scale is
\begin{equation}
 \ell_\lambda:=\left(\frac{\lambda b}{a\sigma}\right)^{1/(1+H)}.
 \label{eq:natural-scales-intro}
\end{equation}
Rescaling time, the target, the control and the cost transforms the
small-cost problem into a normalized tracking problem, called the
\emph{cell problem}, with $a=b=\sigma=\lambda=1$. Its horizon is
$L:=T/\ell_\lambda$, which tends to infinity as $\lambda\to0$.
The limiting optimal cost per unit time is denoted by $\kappa_H$. A \emph{terminal reset} is a final adjustment of the control to the target.
Adding its cost makes the cell value subadditive, and the extra cost
is sublinear in $L$. Consequently,
\begin{equation}
 V_\lambda(T;a,b,\sigma)
 \underset{\lambda\to0}{\sim} \kappa_H T a^{(1-H)/(1+H)}
 (\lambda b)^{2H/(1+H)}\sigma^{2/(1+H)},\qquad \kappa_H\in(0,\infty).
 \label{eq:intro-main-equivalent}
\end{equation}

\emph{Full past information} means that the controller also observes the
target at all negative times, using a two-sided fractional Brownian
motion. This additional information can change finite-horizon values
but not $\kappa_H$. For any normalized horizon $L>0$, let
$u_H^{\past}(L)$ denote the optimal expected cost of the normalized
cell problem with full past information and an additional terminal
reset cost. Averaging over shifted partitions gives
\[
 0\le u_H^{\past}(L)-\kappa_HL\le2\sqrt{2\kappa_H},
\]
which also yields continuity in $H$ when fractional Brownian motions
with different Hurst parameters are constructed from the same
Brownian motion. For nonstationary targets, let
$\rho:[0,T]\to(0,\infty)$ be the continuous \emph{local amplitude}:
near time $t$, target increments over a small interval of length $h$
have characteristic scale $\rho(t)h^H$. Thus $\rho(t)$ plays the
local role of $\sigma$ in the homogeneous model.
The filtered tangent condition in \cref{ass:strong-tangent}
couples the normalized local paths with standard fractional Brownian
motion and matches their observation filtrations. Under this condition,
for continuous positive weights $a$ and $b$,
\[
 \lambda^{-2H/(1+H)}V_\lambda(\Theta)
 \xrightarrow[\lambda\to0]{}\kappa_H\int_0^T
 a_t^{(1-H)/(1+H)}b_t^{2H/(1+H)}\rho(t)^{2/(1+H)}\dd t.
\]
We verify the condition for Gaussian Volterra
targets when their driving Brownian motion is observed.
Perturbing the coefficients then identifies the leading tracking and
variation costs for every \emph{asymptotically optimal family of controls},
meaning one whose cost divided by $V_\lambda(\Theta)$ tends to one.
\paragraph{Contributions.}
The asymptotic formula for fractional Brownian motion in \cref{thm:cell}
extends to local Volterra targets through
\cref{thm:exact-local,prop:volterra-strong-tangent}. The proof uses the information and reset estimates in
\cref{thm:past-robustness,thm:past-finite-horizon}. Prediction states
account for the past information carried by fractional targets. To estimate the cell constant, we derive primal and
dual bounds and approximations based on exponential sums with positive weights
in \cref{sec:bounds,sec:factor-convergence}.
\paragraph{Related work.}
Our Brownian benchmark builds on singular-control results
\cite{karatzas1983,harrison1983,taksar1985} and small-cost tracking of
continuous It\^o targets \cite{cai2017lp,cai2017feedback}. Fractional linear--quadratic control with
absolutely continuous controls is treated in
\cite{hu2005,kleptsyna2003,kleptsyna2005,duncan2012,duncan2013}.
Schelling \cite{schelling2019} links rough-target tracking to an ergodic
fractional cell under quadratic trading costs. In our setting, the penalty
is linear in variation and permits singular controls. One-sided fractional singular control is
studied in \cite{ghosh2010,leeweerasinghe2011}. Truncated variation and hard-tube approximation
\cite{bednorz2015,lochowski2013,lochowski2014,ghomrasni2015}
impose pathwise constraints. Related
level-crossing asymptotics for fractional Brownian motion are studied in
\cite{das2025}.
Prediction states and Markovian fractional or Volterra representations
appear in \cite{carmona1998,carmona2000,hairer2005,viens2019,cuchiero2020,kirchner2025}, while
positive exponential approximations are studied in
\cite{abijaber2019,harms2021,bayer2023,alfonsi2024}.
We approximate control values under a fixed observation filtration, with
stationary-past and growing-horizon error bounds.
The hedging application in \cref{sec:hedging} gives a penalized mean-square
identity. For self-financing models with proportional transaction costs
in fractional finance, see \cite{guasoni2006}.
\paragraph{Organization.}
We first establish the rate and asymptotic formula in
\cref{sec:model-rate,sec:cell}, before turning to Brownian reflection
and past information in \cref{sec:brownian-benchmark,sec:memory-state}.
These results lead to the local limit theorem and its Volterra
application in \cref{sec:local-homogenization}.
We then develop the bounds and factor approximations in
\cref{sec:bounds,sec:factor-convergence}.
The appendices give the filtration and path-space arguments, the
prediction-state construction, and the technical proofs used in the
main text. They also contain the multi-window bounds and the hedging
identity in \cref{app:multiwindow-bounds,sec:hedging}.
\section{Model and rate bounds}
\label{sec:model-rate}
\subsection*{Notation and admissibility}
Fix a complete filtered probability space
$(\Omega,\mathcal F,\mathbf F,\Pp)$, where
$\mathbf F=(\mathcal F_t)_{t\ge0}$, and a horizon $T>0$.
All filtrations are complete and right-continuous. Expectations, variances,
and covariances are denoted by $\E$, $\Var$, and $\Cov$.
A bar denotes completion under the indicated probability measure.
For a standard Borel space $E$, write $\mathcal B(E)$ for its Borel
sigma-field and $\mathcal P(E)$ for the set of probability measures on $E$.
We use $L^0(\Omega,\mathcal H,\Pp;E)$ for the $\mathcal H$-measurable
$E$-valued random variables, identified almost surely. For a separable Banach space
$(E,\|\cdot\|_E)$ and $1\le p<\infty$, the subspace
$L^p(\Omega,\mathcal H,\Pp;E)$ of $L^0(\Omega,\mathcal H,\Pp;E)$
consists of the variables with
$\E[\|X\|_E^p]<\infty$ and carries the norm
$\|X\|_{L^p}:=(\E[\|X\|_E^p])^{1/p}$.
Likewise, $L^\infty$ uses the essential-supremum norm, while deterministic $L^p(I;E)$
uses Lebesgue measure on $I$. Omitted sigma-fields and measures are those
of the ambient space. An omitted scalar codomain is $\mR$.
$\E_\mu$ denotes expectation under a specified law $\mu$.
On $\mR^d$, $x\cdot y=x^\top y$ and $\|x\|=(x^\top x)^{1/2}$ are
the Euclidean inner product and norm. Inner products and norms with subscript $2$ refer
to the stated $L^2$ space. We write $\mathcal L(X)$ for the law of $X$,
$\overset d=$ for equality in law, and
$U\indep\mathcal H\mid Y$ for conditional independence given $Y$.
Conditional essential infima are taken in the almost-sure order. $\mathcal C(I;E)$ and $\mathcal D(I;E)$ denote continuous and c\`adl\`ag
paths, with the uniform and Skorokhod $J_1$ topologies on compact $I$.
$\mathcal C_{\loc}$ uses uniform convergence on compact intervals.
For H\"older estimates, use
\[
 [f]_\alpha:=\sup_{s\ne t}\frac{|f(t)-f(s)|}{|t-s|^\alpha},
 \qquad 0<\alpha\le1.
\]
For $0<\alpha<1$, $\mathcal C^\alpha$ is the space of continuous
functions with finite $[f]_\alpha$. When derivatives are specified,
$\mathcal C^k$ denotes $k$ continuous derivatives.
Set $\mathcal C_T:=\mathcal C([0,T];\mR)$ and
$\mathcal D_T:=\mathcal D([0,T];\mR)$.
For $x\in\mathcal D_T$ and $0\le r<s\le T$, define
\[
 \TV_{(r,s]}(x):=
 \sup_{r=t_0<\cdots<t_n=s}\sum_{j=1}^n|x_{t_j}-x_{t_{j-1}}|,
\]
and let $\mathrm{BV}([0,T];\mR)$ be the set of paths in $\mathcal D_T$
with $\TV_{(0,T]}(x)<\infty$. Given a left value $x_{r-}$, include
the initial adjustment by setting
\[
 \CTV_{[r,s]}(x):=|x_r-x_{r-}|+\TV_{(r,s]}(x).
\]
We write $\Delta x_t:=x_t-x_{t-}$ for the jump at time $t$.
The signed Stieltjes measure $\dd x$ includes the prescribed initial
atom. Its total-variation measure $|\dd x|$ is also written $\dd|x|_t$. The continuous part is $x^c$.
$\mR$, $\Qq$, and $\Nn=\{1,2,\ldots\}$ denote the real, rational, and
positive integer numbers, and we set $\Qq_+:=\Qq\cap[0,\infty)$.
For real $x,y$, $x_+:=\max(x,0)$, $x\wedge y:=\min(x,y)$,
$x\vee y:=\max(x,y)$, and $\sign(0):=0$.
We use $\ind_A$ for the indicator of $A$ and $\lfloor x\rfloor$ and
$\lceil x\rceil$ for the floor and ceiling. In linear algebra,
$\Span$ denotes the linear span, $\diag$ a diagonal matrix, and $\top$ transpose.
$N(0,1)$ denotes a standard normal variable, and
$m_1:=\E\left[|N(0,1)|\right]=\sqrt{2/\pi}$.
Constants $C$ may change between lines. Subscripts indicate parameter
dependence. The symbols $O$, $o$, $\sim$, and $\asymp$ mean a bounded
ratio, a vanishing ratio, a ratio tending to one, and two-sided comparison
in the stated limit. Weak convergence of finite measures is tested
against bounded continuous functions.
\begin{definition}[augmented total variation and admissible controls]
\label{def:admissible-controls}
Let $\Theta\in L^0(\Omega,\mathcal F,\Pp;\mathcal C_T)$ be
$\mathbf F$-adapted. A \emph{control} is a jointly
$\mathcal F\otimes\Borel([0,T])$-measurable map
\[
 \appdef{\phi}{\Omega\times[0,T]}{\mR}
 { (\omega,t)}{\phi_t(\omega)}
\]
such that
\[
 \phi(\omega,\cdot)\in\mathrm{BV}([0,T];\mR)
 \quad\Pp\text{-a.s.},
 \qquad
 \phi_t\in L^0(\Omega,\mathcal F_t,\Pp),
 \quad 0\le t\le T.
\]
With the pre-initial value $\phi_{0-}:=\Theta_0$, the \emph{augmented variation} is
\begin{equation}
 \CTV_{[0,T]}(\phi)
 :=
 |\phi_0-\phi_{0-}|+\TV_{(0,T]}(\phi).
 \label{eq:variation-convention}
\end{equation}
The measure $|\dd\phi|$ on $[0,T]$ is determined by
\[
 |\dd\phi|(\{0\})=|\phi_0-\phi_{0-}|,
 \qquad
 |\dd\phi|((s,t])=\TV_{(s,t]}(\phi),
 \quad 0\le s<t\le T.
\]
Let $\cA_T(\Theta)$ be the class of these controls satisfying
\[
 \E\left[\int_0^T(\phi_t-\Theta_t)^2\dd t
       +\CTV_{[0,T]}(\phi)\right]<\infty.
\]
When $\Theta$ is fixed, write $\cA_T:=\cA_T(\Theta)$.
\end{definition}
\begin{definition}[tracking value]
\label{def:tracking-value}
Let $a,b:[0,T]\to(0,\infty)$ be Borel functions with
$a_-\le a_t\le a_+$ for $\dd t$-a.e. $t$ and
$b_-\le b_t\le b_+$ for every $t\in[0,T]$, where
$a_-,b_->0$ and $a_+,b_+<\infty$. For $\lambda>0$, the \emph{tracking value} is
\begin{equation}
 V_\lambda(\Theta)
 :=
 \inf_{\phi\in\cA_T(\Theta)}
 \E\left[
 \frac12\int_0^T a_t(\phi_t-\Theta_t)^2\dd t
 +\lambda\int_{[0,T]}b_t\dd|\phi|_t
 \right].
 \label{eq:model}
\end{equation}
For $x_0\in\mR$, $V_\lambda^{x_0}(\Theta)$ denotes the same infimum with
$\phi_{0-}=x_0$.
\end{definition}
\begin{proposition}[effect of the initial position]
\label{prop:initial-inventory}
Assume $\cA_T(\Theta)\ne\varnothing$ and $\E[|\Theta_0|]<\infty$. Then
\begin{equation}
 \left|V_\lambda^{x_0}(\Theta)-V_\lambda(\Theta)\right|
 \le
 \lambda b_+\E\left[|x_0-\Theta_0|\right].
 \label{eq:initial-inventory-bound}
\end{equation}
Therefore, every asymptotic equivalent of order $\lambda^\alpha$ with
$\alpha\in(0,1)$ is unchanged by a fixed initial position.
\end{proposition}
\begin{proof}
Write $J_\lambda^{x_0}$ for the objective with pre-initial value $x_0$.
For the same control path, only the initial variation atom changes:
\begin{align*}
 |J_\lambda^{x_0}(\phi;\Theta)-J_\lambda(\phi;\Theta)|
 &=\lambda b_0\left|\E\bigl[|\phi_0-x_0|-|\phi_0-\Theta_0|\bigr]\right|\\
 &\le
 \lambda b_+\E[|x_0-\Theta_0|].
\end{align*}
By the triangle inequality, the two admissible classes contain the same
paths. Taking infima in both directions proves
\eqref{eq:initial-inventory-bound}. The normalized error is bounded by
$b_+\E\left[|x_0-\Theta_0|\right]\lambda^{1-\alpha}\xrightarrow[\lambda\to0]{}0$.
\end{proof}
\begin{proposition}[finite-variation perturbations]
\label{prop:BV-perturbation}
Assume $\cA_T(\Theta)\ne\varnothing$. Let $A$ be continuous, adapted, and of finite variation, with
$\E[\TV_{(0,T]}(A)]<\infty$. Using each target's pre-initial value,
\begin{equation}
 |V_\lambda(\Theta+A)-V_\lambda(\Theta)|
 \le\lambda\E\left[\int_{(0,T]}b_t\,|\dd A|_t\right].
 \label{eq:BV-perturbation}
\end{equation}
Thus an integrable finite-variation perturbation does not change a leading
term of order $\lambda^\alpha$ with $0<\alpha<1$.
\end{proposition}
\begin{proof}
Adding $A$ defines a bijection from $\cA_T(\Theta)$ to
$\cA_T(\Theta+A)$. It preserves adaptedness and satisfies
\[
 (\phi_t+A_t)-(\Theta_t+A_t)=\phi_t-\Theta_t,
 \qquad
 |(\phi_0+A_0)-(\Theta_0+A_0)|=|\phi_0-\Theta_0|.
\]
For the signed Stieltjes measures on $(0,T]$,
\begin{align*}
 |\dd(\phi+A)|&\le|\dd\phi|+|\dd A|.
\end{align*}
Then,
\begin{align*}
 \int_{(0,T]} b_t|\dd(\phi+A)|_t
 &\le\int_{(0,T]} b_t|\dd\phi|_t+\int_{(0,T]} b_t|\dd A|_t.
\end{align*}
Hence
\[
 V_\lambda(\Theta+A)\le V_\lambda(\Theta)
 +\lambda\E\left[\int_{(0,T]}b_t|\dd A|_t\right].
\]
Replace $(\Theta,A)$ by $(\Theta+A,-A)$ for the reverse bound. Division
by $\lambda^\alpha$ leaves an $O(\lambda^{1-\alpha})$ error.
\end{proof}
\begin{assumption}[two-sided increment scale]
\label{ass:increments}
There are $H\in(0,1)$, $h_0>0$, and $c_1,C_2>0$ such that, for $0<h<\min(h_0,T)$,
\begin{equation*}
 \int_0^{T-h}\E\left[|\Theta_{t+h}-\Theta_t|\right]\dd t
 \ge c_1(T-h)h^H,
\end{equation*}
and
\begin{equation}
 \sup_{0\le t\le T-h}\E\left[|\Theta_{t+h}-\Theta_t|^2\right]
 \le C_2h^{2H}.
 \label{eq:increment-assumption}
\end{equation}
\end{assumption}
\begin{definition}[noncausal comparison class]
\label{def:noncausal-class}
Let $\cA_T^{\rm all}(\Theta)$ be the set of
$\mathcal D_T$-valued random variables $\phi$ satisfying the pathwise,
pre-initial-value, and integrability conditions of
$\cA_T(\Theta)$, but without the conditions
$\phi_t\in L^0(\Omega,\mathcal F_t,\Pp)$.
When $\Theta$ is fixed, write $\cA_T^{\rm all}:=\cA_T^{\rm all}(\Theta)$.
\end{definition}
\begin{lemma}[increment versus energy]
\label{lem:increment}
For $\phi\in\cA_T^{\rm all}$, let
$M=\E[\int_0^T(\phi_t-\Theta_t)^2\dd t]$ and
$R=\E[\CTV_{[0,T]}(\phi)]$. For $0<h<T$,
\begin{equation}
 \int_0^{T-h}\E[|\Theta_{t+h}-\Theta_t|]\dd t
 \le hR+2\sqrt{TM}.
 \label{eq:increment-energy}
\end{equation}
\end{lemma}
\begin{proof}
Put $e:=\phi-\Theta\in L^2(\Omega\times[0,T],\Pp\otimes\dd t)$.
Pathwise,
\[
 |\Theta_{t+h}-\Theta_t|
 \le\int_{(t,t+h]}|\dd\phi|_s+|e_t|+|e_{t+h}|.
\]
Tonelli's theorem and the window length give
\begin{align*}
 \int_0^{T-h}\int_{(t,t+h]}|\dd\phi|_s\dd t
 &=
 \int_{(0,T]}\left(\int_0^{T-h}
       \ind_{\{t<s\le t+h\}}\dd t\right)|\dd\phi|_s\\
 &=\int_{(0,T]}|[s-h,s)\cap[0,T-h]|\,|\dd\phi|_s\\
 &\le h\CTV_{[0,T]}(\phi),
\end{align*}
where $|I|$ denotes the Lebesgue length of an interval $I$.
For either $I=[0,T-h]$ or $I=[h,T]$, Cauchy--Schwarz yields
\begin{align*}
 \int_I\E[|e_t|]\dd t
 &\le
 \sqrt{|I|}\left(\int_I\E[|e_t|]^2\dd t\right)^{1/2}\\
 &\le
 \sqrt{|I|}\left(\int_I\E[e_t^2]\dd t\right)^{1/2}
 \le\sqrt{TM}.
\end{align*}
Taking expectations and adding the two error bounds proves the claim.
\end{proof}
\begin{theorem}[rate bounds]
\label{thm:general-rate}
Suppose \cref{ass:increments} holds. With the cutoff
\begin{equation}
 \lambda_0:=\left(1\wedge\frac{h_0}{2}\wedge\frac T2\right)^{1+H},
 \label{eq:rate-lambda0}
\end{equation}
define the lower and upper constants by
\begin{equation}
 C_- :=
 T\min\left\{
 \frac{a_-c_1^2}{128},
 \frac{b_-c_1}{4}
 \right\},
 \label{eq:rate-Cminus}
\end{equation}
and
\begin{equation}
 C_+ :=
 \frac{a_+C_2T}{2}
 +b_+\sqrt{C_2}\,(T+1).
 \label{eq:rate-Cplus}
\end{equation}
For $0<\lambda\le\lambda_0$, the value satisfies
\begin{equation}
 C_-\lambda^{2H/(1+H)}
 \le V_\lambda(\Theta)
 \le C_+\lambda^{2H/(1+H)}.
 \label{eq:general-rate}
\end{equation}
The lower bound also holds when the infimum is taken over
$\cA_T^{\rm all}$.
\end{theorem}
\begin{proof}
Put $h:=\lambda^{1/(1+H)}\le1\wedge h_0/2\wedge T/2$ and use
$M,R$ from \cref{lem:increment}. For any admissible path, causal or not,
\[
 J_\lambda(\phi;\Theta)\ge\frac{a_-}{2}M+\lambda b_-R.
\]
Moreover,
\begin{align*}
 hR+2\sqrt{TM}
 &\overset{\eqref{eq:increment-energy}}\ge
 \int_0^{T-h}\E[|\Theta_{t+h}-\Theta_t|]\dd t\\
 &\overset{\cref{ass:increments}}\ge
 c_1(T-h)h^H
 \ge\frac{c_1T}{2}h^H.
\end{align*}
If $M\ge c_1^2Th^{2H}/64$, then
\[
 J_\lambda(\phi;\Theta)
 \ge\frac{a_-M}{2}
 \ge\frac{a_-c_1^2T}{128}h^{2H}.
\]
Otherwise, the increment bound gives
\[
 R\ge\frac{c_1T}{2}h^{H-1}-\frac{2\sqrt{TM}}h
 \ge\frac{c_1T}{4}h^{H-1},
\]
and therefore
\[
 J_\lambda(\phi;\Theta)
 \ge\lambda b_-R
 \ge\frac{b_-c_1T}{4}\lambda h^{H-1}
 =\frac{b_-c_1T}{4}h^{2H}.
\]
Consequently,
\begin{equation}
 J_\lambda(\phi;\Theta)
 \ge T\min\left\{\frac{a_-c_1^2}{128},\frac{b_-c_1}{4}\right\}
 h^{2H}=C_-\lambda^{2H/(1+H)}.
 \label{eq:generic-lower}
\end{equation}
For the upper bound, take a uniform mesh with step $h$ and a possibly
shorter last interval. Set $N:=\lceil T/h\rceil$ and
$t_k:=\min\{kh,T\}$ for $k=0,\ldots,N$.
On $[t_k,t_{k+1})$, set $\phi_t:=\Theta_{t_k}$,
and choose $\phi_T:=\phi_{T-}$ to avoid a terminal jump. The increment
bounds and Cauchy--Schwarz give
\begin{align*}
 \E\left[\int_0^T(\phi_t-\Theta_t)^2\dd t\right]
 \overset{\eqref{eq:increment-assumption}}\le
 C_2\sum_{k=0}^{N-1}\int_{t_k}^{t_{k+1}}(t-t_k)^{2H}\dd t\le C_2Th^{2H},
\end{align*}
and
\begin{align*}
 \E\left[\CTV_{[0,T]}(\phi)\right]
 \le\sum_{k=1}^{N-1}\E\left[|\Theta_{t_k}-\Theta_{t_{k-1}}|\right]\le
 \sqrt{C_2}(T/h+1)h^H.
\end{align*}
Therefore
\begin{align*}
 J_\lambda(\phi;\Theta)
 &\le\frac{a_+C_2T}{2}h^{2H}
 +\lambda b_+\sqrt{C_2}(Th^{H-1}+h^H)\\
 &=
 h^{2H}\left[\frac{a_+C_2T}{2}+b_+\sqrt{C_2}(T+h)\right]\\
 &\le C_+\lambda^{2H/(1+H)}.
\end{align*}
\end{proof}
\section{The fractional Brownian cell and the asymptotic constant}
\label{sec:cell}
Fix $H\in(0,1)$. Standard fBm is the continuous centered Gaussian process
$B^H$ with $B_0^H=0$ and covariance
\begin{equation}
 \appdef{R_H}{\mR^2}{\mR}{(s,t)}
 {\tfrac12\bigl(|s|^{2H}+|t|^{2H}-|t-s|^{2H}\bigr)}.
 \label{eq:fractional-covariance-kernel}
\end{equation}
In particular,
\begin{equation}
 \E\left[(B_t^H-B_s^H)^2\right]=|t-s|^{2H},\qquad
 \E\left[B_s^HB_t^H\right]=R_H(s,t),\quad s,t\ge0.
 \label{eq:fbm-covariance}
\end{equation}
For constant $a,b,\sigma>0$, we retain the notation
$V_\lambda(T;a,b,\sigma)$ introduced in
\cref{sec:introduction}, using the usual natural filtration of $B^H$.
On $\mathcal C_L=\mathcal C([0,L];\mR)$, let
$\mu_H^L$ be the fBm law and let $\mathbf F^{H,L}$ be the usual augmented
coordinate filtration. We use $B^H$ for the coordinate process and
$\mathcal D_L=\mathcal D([0,L];\mR)$ for control paths.
\begin{definition}[one-sided fractional cells]
\label{def:one-sided-cells}
Let $\cA_L^H$ be the set of
$\psi\in L^0(\mathcal C_L,\Borel(\mathcal C_L),\mu_H^L;\mathcal D_L)$
that are $\mathbf F^{H,L}$-adapted, have pre-initial value $\psi_{0-}=0$
and finite variation $\CTV_{[0,L]}(\psi)<\infty$ $\mu_H^L$-a.s.,
and satisfy the integrability condition
\[
 \E\left[\int_0^L(\psi_s-B_s^H)^2\dd s+\CTV_{[0,L]}(\psi)\right]<\infty.
\]
Their tracking and reset values are, respectively,
\begin{equation}
 v_H(L)
 :=
 \inf_{\psi\in\cA_L^H}
 \E\left[\frac12\int_0^L(\psi_s-B_s^H)^2\dd s
 +\CTV_{[0,L]}(\psi)\right],
 \label{eq:v-cell}
\end{equation}
and
\begin{equation}
 u_H(L)
 :=
 \inf_{\psi\in\cA_L^H}
 \E\left[\frac12\int_0^L(\psi_s-B_s^H)^2\dd s
 +\CTV_{[0,L]}(\psi)+|\psi_L-B_L^H|\right].
 \label{eq:u-cell}
\end{equation}
\end{definition}
\begin{proposition}[self-similar scaling identity]
\label{prop:scaling}
With $\ell_\lambda$ from \eqref{eq:natural-scales-intro},
\begin{equation}
 V_\lambda(T;a,b,\sigma)
 =a\sigma^2\ell_\lambda^{2H+1}
 v_H\left(\frac{T}{\ell_\lambda}\right).
 \label{eq:exact-scaling}
\end{equation}
Equivalently, for $L:=T/\ell_\lambda$,
\begin{equation}
 \frac{V_\lambda(T;a,b,\sigma)}
 {T a^{(1-H)/(1+H)}(\lambda b)^{2H/(1+H)}\sigma^{2/(1+H)}}
 =\frac{v_H(L)}L.
 \label{eq:normalized-scaling}
\end{equation}
\end{proposition}
\begin{proof}
For $\ell>0$, the continuous bijection
\[
 \appdef{\mathsf S_\ell}{\mathcal C_T}{\mathcal C_{T/\ell}}
 {x}{\bigl(s\mapsto\ell^{-H}x(\ell s)\bigr)}
\]
transports fBm law and its natural filtration. Transport causal rules by
\cref{thm:canonical-transport}. For $0\le s\le T/\ell$, rescale
the target and control by
\[
 \widetilde B_s^H:=\ell^{-H}B_{\ell s}^H,
 \qquad \psi_s:=\frac{\phi_{\ell s}}{\sigma\ell^H}.
\]
The inverse is $\phi_t=\sigma\ell^H\psi_{t/\ell}$. Both maps preserve
adaptedness, finite variation, and the zero pre-initial value. With
$t=\ell s$,
\begin{equation*}
 \frac a2\E\left[\int_0^T(\phi_t-\sigma B_t^H)^2\dd t\right]
 =\frac{a\sigma^2\ell^{2H+1}}2
   \E\left[\int_0^{T/\ell}(\psi_s-\widetilde B_s^H)^2\dd s\right],
\end{equation*}
and
\begin{equation*}
 \lambda b\E\left[\CTV_{[0,T]}(\phi)\right]
 =\lambda b\sigma\ell^H\E\left[\CTV_{[0,T/\ell]}(\psi)\right].
\end{equation*}
Equality of the two coefficients is equivalent to
\[
 a\sigma^2\ell^{2H+1}=\lambda b\sigma\ell^H
 \quad\iff\quad\ell^{H+1}=\frac{\lambda b}{a\sigma}
 \quad\iff\quad\ell=\ell_\lambda.
\]
Taking infima in both directions gives \eqref{eq:exact-scaling}. Finally,
\begin{align*}
 a\sigma^2\ell_\lambda^{2H}
 =a\sigma^2\left(\frac{\lambda b}{a\sigma}\right)^{2H/(1+H)}
 =a^{(1-H)/(1+H)}(\lambda b)^{2H/(1+H)}\sigma^{2/(1+H)}.
\end{align*}
Divide \eqref{eq:exact-scaling} by $Ta\sigma^2\ell_\lambda^{2H}$ to
obtain \eqref{eq:normalized-scaling}.
\end{proof}
\begin{proposition}[reset concatenation]
\label{prop:reset-concatenation}
For every $S,T>0$,
\begin{equation}
 u_H(S+T)\le u_H(S)+u_H(T).
 \label{eq:u-subadditive}
\end{equation}
\end{proposition}
\begin{proof}
Choose Borel reset rules $\psi^S,\psi^T$ with expected costs at most
$u_H(S)+\eta$ and $u_H(T)+\eta$, respectively. Apply the second rule
to the increment path $Y_t:=B_{S+t}^H-B_S^H$, $0\le t\le T$, and concatenate by
\[
 \psi_t:=
 \begin{cases}
 \psi_t^S,&0\le t<S,\\
 B_S^H+\psi_{t-S}^T(Y),&S\le t\le S+T.
 \end{cases}
\]
Stationary increments give $\mathcal L(Y)=\mu_H^T$, and the second rule is
adapted to the original filtration. At the interface,
\[
 |B_S^H+\psi_0^T(Y)-\psi_{S-}^S|
 \le|\Delta\psi_S^S|+|\psi_S^S-B_S^H|+|\psi_0^T(Y)|.
\]
The first rule's variation, its reset cost, and the second rule's initial
adjustment cost cover these terms, respectively. Running costs add, and
$|\psi_{S+T}-B_{S+T}^H|=|\psi_T^T(Y)-Y_T|$.
Thus \cref{lem:canonical-concatenation} gives
\[
 u_H(S+T)\le u_H(S)+u_H(T)+2\eta.
\]
Letting $\eta\to0$ proves the claim.
\Cref{fig:reset-concatenation} summarizes the construction.
\end{proof}
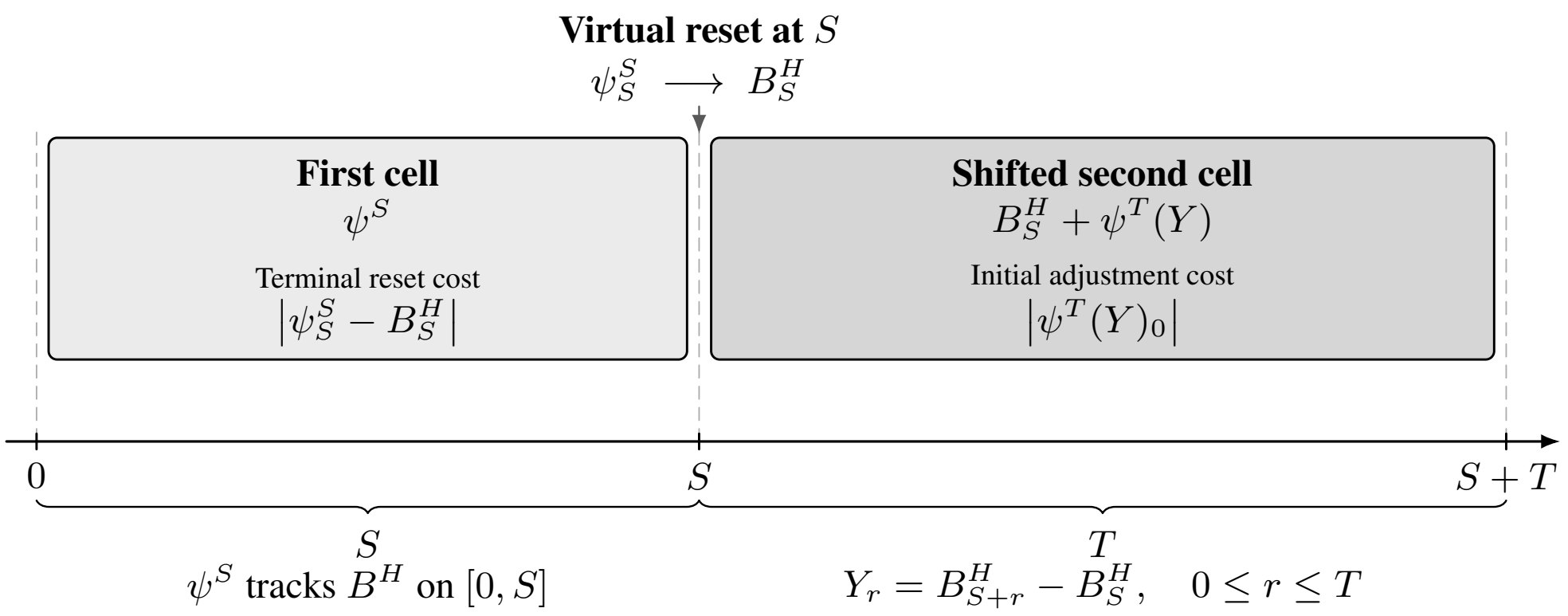
\begin{figure}[tbp]
\centering
\begin{tikzpicture}[
  x=1cm,y=1cm,
  font=\small,
  line cap=round,line join=round,
  >={Latex[length=1.6mm]},
  every node/.style={text=black,inner sep=2pt}
]
  \draw[black!32,densely dashed,line width=.35pt] (0,0) -- (0,2.58);
  \draw[black!32,densely dashed,line width=.35pt] (5.5,0) -- (5.5,2.58);
  \draw[black!32,densely dashed,line width=.35pt] (12.2,0) -- (12.2,2.58);

  \draw[fill=black!8,draw=black,line width=.55pt,rounded corners=2pt]
    (.10,.68) rectangle (5.40,2.52);
  \draw[fill=black!16,draw=black,line width=.55pt,rounded corners=2pt]
    (5.60,.68) rectangle (12.10,2.52);

  \node[text=black,font=\small\bfseries] at (2.75,2.23) {First cell};
  \node at (2.75,1.84) {$\psi^S$};
  \node[text=black,font=\scriptsize] at (2.75,1.36) {Terminal reset cost};
  \node at (2.75,.99) {$\left|\psi^S_S-B_S^H\right|$};

  \node[text=black,font=\small\bfseries] at (8.85,2.23) {Shifted second cell};
  \node at (8.85,1.84) {$B_S^H+\psi^T(Y)$};
  \node[text=black,font=\scriptsize] at (8.85,1.36) {Initial adjustment cost};
  \node at (8.85,.99) {$\left|\psi^T(Y)_0\right|$};

  \node[font=\small\bfseries] at (5.5,3.42) {Virtual reset at $S$};
  \node at (5.5,3.00) {$\psi^S_S\;\longrightarrow\;B_S^H$};
  \draw[->,black!65,line width=.45pt] (5.5,2.78) -- (5.5,2.55);

  \draw[->,black,line width=.6pt] (-.25,0) -- (12.65,0);
  \foreach \p in {0,5.5,12.2}{
    \draw[black,line width=.55pt] (\p,-.065) -- (\p,.065);
  }
  \node[below=3pt] at (0,0) {$0$};
  \node[below=3pt] at (5.5,0) {$S$};
  \node[below=3pt] at (12.2,0) {$S+T$};

  \draw[decorate,decoration={brace,mirror,amplitude=3pt},
        black,line width=.45pt]
    (0,-.49) -- (5.5,-.49)
    node[midway,below=5pt,text=black] {$S$};
  \draw[decorate,decoration={brace,mirror,amplitude=3pt},
        black,line width=.45pt]
    (5.5,-.49) -- (12.2,-.49)
    node[midway,below=5pt,text=black] {$T$};

  \node at (2.75,-1.20) {$\psi^S\text{ tracks }B^H\text{ on }[0,S]$};
  \node at (8.85,-1.20) {$Y_r=B_{S+r}^H-B_S^H,\quad 0\le r\le T$};
\end{tikzpicture}
\space
\caption{Reset concatenation of two causal cell rules. The second rule
tracks the increment path $Y_r=B_{S+r}^H-B_S^H$ and is translated by
$B_S^H$. The virtual reset and the second rule's initial adjustment
bound the interface variation.}
\label{fig:reset-concatenation}
\end{figure}
\begin{theorem}[cell limit and small-cost asymptotic]
\label{thm:cell}
For every $H\in(0,1)$, there is a constant $\kappa_H\in(0,\infty)$ such that
\begin{equation}
 \lim_{L\to\infty}\frac{v_H(L)}L
 =\lim_{L\to\infty}\frac{u_H(L)}L
 =\kappa_H,
 \qquad
 \kappa_H=\inf_{L>0}\frac{u_H(L)}L.
 \label{eq:kappa-cell}
\end{equation}
Therefore,
\begin{equation}
 V_\lambda(T;a,b,\sigma)
 \underset{\lambda\to0}{\sim}
 \kappa_H T a^{(1-H)/(1+H)}(\lambda b)^{2H/(1+H)}
 \sigma^{2/(1+H)}.
 \label{eq:sharp-fbm-equivalent}
\end{equation}
\end{theorem}
\begin{proof}
\emph{Reset-cell limit.}
The zero rule gives
\[
 u_H(L)\le\frac12\int_0^Ls^{2H}\dd s+\E\left[|B_L^H|\right]
 =\frac{L^{2H+1}}{2(2H+1)}+m_1L^H.
\]
Hence $u_H$ is locally bounded. Subadditivity and \cref{app:fekete} yield
\begin{equation}
 \lim_{L\to\infty}\frac{u_H(L)}L
 =\inf_{L>0}\frac{u_H(L)}L=:\kappa_H.
 \label{eq:fekete-u}
\end{equation}
A rule that samples the target at integer times and holds its value
between samples, followed by a reset, satisfies
\begin{align*}
 \frac12\E\left[\int_0^L(\psi_s-B_s^H)^2\dd s\right]&\le L/2,
\end{align*}
and
\begin{align*}
 \E\left[\CTV_{[0,L]}(\psi)+|\psi_L-B_L^H|\right]
 &\le m_1(L+1),
\end{align*}
so $v_H(L)\le u_H(L)=O(1+L)$.\\
\emph{Removing the reset.}
For $L\ge1$, take $\psi$ with non-reset cost at most $v_H(L)+\eta$.
Since
\[
 \int_{L-1}^L\E\left[|\psi_r-B_r^H|^2\right]\dd r\le2(v_H(L)+\eta),
\]
choose a deterministic $r\in[L-1,L)$ with
$\E\left[|\psi_r-B_r^H|^2\right]\le2(v_H(L)+\eta)$.
Keep $\psi$ before $r$, jump to $B_r^H$, and hold until $L$. The
interface inequality
\[
 |B_r^H-\psi_{r-}|
 \le|\Delta\psi_r|+|\psi_r-B_r^H|
\]
shows that the extra variation is at most $|\psi_r-B_r^H|$.
\begin{samepage}
The other added terms satisfy
\begin{align*}
 \frac12\E\left[\int_r^L(B_t^H-B_r^H)^2\dd t\right]
 \overset{\eqref{eq:fbm-covariance}}=
 \frac{(L-r)^{2H+1}}{2(2H+1)}\le\frac1{2(2H+1)},
\end{align*}
and
\begin{align*}
 \E[|B_L^H-B_r^H|]&=m_1(L-r)^H\le m_1.
\end{align*}
\end{samepage}
Cauchy--Schwarz and $\eta\to0$ therefore give
\begin{equation}
 0\le u_H(L)-v_H(L)
 \le\sqrt{2v_H(L)}+\frac1{2(2H+1)}+m_1,
 \qquad L\ge1.
 \label{eq:reset-sublinear}
\end{equation}
The gap is $O(\sqrt L)$, so $v_H(L)/L\xrightarrow[L\to\infty]{}\kappa_H$.\\
\emph{Positivity and scaling.}
For standard fBm on $[0,1]$, the increment assumptions hold with
$c_1=m_1$ and $C_2=1$. Taking $\lambda=L^{-(1+H)}$ gives
$\ell_\lambda=L^{-1}$ and, for large $L$,
\[
 L^{-(2H+1)}v_H(L)=V_\lambda(1;1,1,1)
 \ge
 C\lambda^{2H/(1+H)}=CL^{-2H}.
\]
Thus $v_H(L)/L\ge C>0$. Now use
\eqref{eq:normalized-scaling}, since $T/\ell_\lambda\xrightarrow[\lambda\to0]{}\infty$.
\end{proof}
\begin{corollary}[sublinear reset gap]
\label{cor:reset-gap}
There is $C_H<\infty$ such that
\begin{equation}
 0\le u_H(L)-v_H(L)\le C_H(1+\sqrt L),\qquad L\ge1.
 \label{eq:reset-gap}
\end{equation}
The reset-gap estimate follows from \eqref{eq:reset-sublinear} and the
upper bound obtained with a time step of one.
\end{corollary}
\section{The Brownian case}
\label{sec:brownian-benchmark}
For $H=1/2$, the tracking error is a sufficient state for the controlled
Brownian dynamics.
\begin{proposition}[Brownian cell constant]
\label{prop:brownian}
Reflection of $e=\psi-B$ on $[-d_*,d_*]$ is optimal in the long-run cell,
where
\begin{equation}
 d_*=(3/2)^{1/3},\qquad
 \kappa_{1/2}=\frac{d_*^2}{2}=\frac12(3/2)^{2/3}.
 \label{eq:brownian-kappa}
\end{equation}
\end{proposition}
For reflection on $[-d,d]$, the long-run tracking and variation cost
rates are $d^2/6$ and $1/(2d)$, respectively. The total cost rate is
\[
 \appdef{c}{(0,\infty)}{(0,\infty)}{d}{\frac{d^2}{6}+\frac1{2d}}.
\]
Its unique minimum occurs at $d^3=3/2$.
\begin{proof}
See \cref{app:brownian-verification}.
\end{proof}
\section{Past information and prediction states}
\label{sec:memory-state}
An interior cell has access to observations made before its local starting
time. We show that this information does not change the long-run value.
\subsection{The cell with past information}
Let $\mathbf B^H=(\mathbf B_t^H)_{t\in\mR}$ be a two-sided standard fBm
with $\mathbf B_0^H=0$ and covariance $R_H$.
Its canonical space is
$\mathcal C_0(\mR)=\{x\in\mathcal C(\mR;\mR):x(0)=0\}$,
with law $\mu_H^{\past}$ and the local-uniform topology.
For control times $t\ge0$, use the full-past filtration
$\mathcal F_t^{\mathbf B}:=
\overline{\sigma(\mathbf B_s^H:s\le t)}^{\mu_H^{\past}}$,
which is right-continuous by \cref{lem:mvn-filtration}.
\begin{definition}[cells with past information]
\label{def:past-cells}
Let $\cA_L^{H,\past}$ be the admissible class on $[0,L]$ with pre-initial
value zero and filtration $\mathbf F^{\mathbf B}$. Its tracking and reset
values are, respectively,
\begin{equation}
 v_H^{\past}(L)
 :=\inf_{\psi\in\cA_L^{H,\past}}
 \E\left[\frac12\int_0^L(\psi_t-\mathbf B_t^H)^2\dd t
                    +\CTV_{[0,L]}(\psi)\right],
 \label{eq:v-past}
\end{equation}
and
\begin{equation}
 u_H^{\past}(L)
 :=\inf_{\psi\in\cA_L^{H,\past}}
 \E\left[\frac12\int_0^L(\psi_t-\mathbf B_t^H)^2\dd t
        +\CTV_{[0,L]}(\psi)+|\psi_L-\mathbf B_L^H|\right].
 \label{eq:u-past}
\end{equation}
For $M>0$, let $\mathbf F^{(M)}$ be the usual augmentation of
$\sigma(\mathbf B_s^H:-M\le s\le t)$, $0\le t\le L$.
Restricting controls to this filtration defines
$u_{H,M}^{\past}(L)$ and $v_{H,M}^{\past}(L)$.
Write $\cA_{L,M}^{H,\past}$ for the corresponding admissible class.
\end{definition}
\begin{lemma}[finite-memory approximation]
\label{lem:finite-memory-density}
For every fixed $L>0$,
\begin{equation}
 u_{H,M}^{\past}(L)\underset{M\to\infty}{\downarrow} u_H^{\past}(L).
 \label{eq:finite-memory-density}
\end{equation}
The same result holds for $v_H^{\past}(L)$.
\end{lemma}
\begin{proof}
See \cref{app:finite-memory}.
\end{proof}
\begin{theorem}[effect of past information]
\label{thm:past-robustness}
For every $H\in(0,1)$,
\begin{equation}
 \lim_{L\to\infty}\frac{v_H^{\past}(L)}L
 =\lim_{L\to\infty}\frac{u_H^{\past}(L)}L
 =\inf_{L>0}\frac{u_H^{\past}(L)}L=\kappa_H.
 \label{eq:past-kappa}
\end{equation}
\end{theorem}
\begin{proof}
The process
$(\mathbf B_{S+t}^H-\mathbf B_S^H)_{t\in\mR}$ has the two-sided
fBm law and is observed up to index $t$ at time $S+t$.
Reset concatenation gives
\[
 u_H^{\past}(S+T)\le u_H^{\past}(S)+u_H^{\past}(T).
\]
The zero rule and the rule with time step one give local boundedness and an $O(1+L)$
upper bound. The reset-removal proof of \cref{thm:cell} uses only
increment moments and yields
\[
 0\le u_H^{\past}(L)-v_H^{\past}(L)
 \le\sqrt{2v_H^{\past}(L)}+\frac1{2(2H+1)}+m_1.
\]
Thus both normalized values have the common limit
\[
 \inf_{L>0}\frac{u_H^{\past}(L)}L\le\kappa_H,
\]
where the inequality follows by ignoring the past.
For the converse, fix $L,M>0$ and a finite-memory reset rule of cost
at most $u_{H,M}^{\past}(L)+\eta$. On a one-sided path, first run a
finite-cost reset rule on $[0,M]$, with cost $C_M$. On each block
starting at $s_k:=M+kL$, apply the selected rule to the segment
\[
 (B_{s_k+t}^H-B_{s_k}^H)_{-M\le t\le L},
\]
and add $B_{s_k}^H$ to its output. This history is observed because
$s_k-M\ge0$, and stationary increments give the required segment law.
The reset costs cover the variation at every interface, hence
\[
 u_H(M+nL)\le C_M+n\{u_{H,M}^{\past}(L)+\eta\}.
\]
After division by $M+nL$ and $n\to\infty$,
\[
 \kappa_H\le\frac{u_{H,M}^{\past}(L)+\eta}{L}
 \xrightarrow[\eta\to0,\ M\to\infty]{}
 \frac{u_H^{\past}(L)}L
\]
by \cref{lem:finite-memory-density}. Taking the infimum over $L$ proves
equality.
\end{proof}
\begin{theorem}[a bounded full-past reset correction]
\label{thm:past-finite-horizon}
For every $H\in(0,1)$ and $L>0$,
\begin{equation}
 0\le u_H^{\past}(L)-\kappa_H L\le2\sqrt{2\kappa_H}.
 \label{eq:past-finite-horizon}
\end{equation}
In particular,
\[
 0\le\frac{u_H^{\past}(L)}L-\kappa_H\le\frac{C_0}{L},
 \qquad\text{with}\quad C_0:=2\sqrt{2(1/2+\sqrt{2/\pi})},
\]
uniformly in $H\in(0,1)$.
\end{theorem}
\begin{proof}
The lower bound is \eqref{eq:past-kappa}. Fix $L>0$. For $R>2L$, take
a past-adapted rule with non-reset cost $J_R$ satisfying
$J_R/R\xrightarrow[R\to\infty]{}\kappa_H$. For the error
$e:=\psi-\mathbf B^H$, the energy
$M_R:=\E[\int_0^Re_t^2\dd t]$ therefore satisfies $M_R\le2J_R$.
For a deterministic $s\in[0,L)$, set
$n_s:=\lfloor(R-s)/L\rfloor$ and use grid points $s+jL$,
$0\le j\le n_s$.
Each complete block, translated and given a virtual initial adjustment
and terminal reset, is a full-past cell. Disjointness of the original
costs and \cref{lem:virtual-reset-decomposition} imply, for almost every $s$,
\[
 n_su_H^{\past}(L)\le J_R+2\sum_{j=0}^{n_s}\E\left[|e_{s+jL}|\right],
 \qquad n_s\ge R/L-2.
\]
Averaging over the grid offset and applying Tonelli's theorem gives
\begin{align*}
 \int_0^L\sum_{j=0}^{\lfloor(R-s)/L\rfloor}\E[|e_{s+jL}|]\dd s&=
 \sum_{j=0}^{\infty}\int_0^L\ind_{\{s+jL\le R\}}\E\left[|e_{s+jL}|\right]\dd s\\
 &=\sum_{j=0}^{\infty}\int_{[jL,(j+1)L)\cap[0,R]}\E[|e_t|]\dd t\\
 &=\int_0^R\E[|e_t|]\dd t.
\end{align*}
Thus
\begin{align*}
 \left(\frac RL-2\right)u_H^{\past}(L)
 \le J_R+\frac2L\int_0^R\E[|e_t|]\dd t
 \le
 J_R+\frac2L\sqrt{RM_R}.
\end{align*}
Multiply by $L/R$ to get
\begin{align*}
 \left(1-\frac{2L}{R}\right)u_H^{\past}(L)
 \le L\frac{J_R}{R}+2\sqrt{\frac{M_R}{R}}
 \le
 L\frac{J_R}{R}+2\sqrt{\frac{2J_R}{R}}.
\end{align*}
Letting $R\to\infty$ proves the upper bound. Finally,
\[
 \kappa_H\le u_H(1)\le\frac1{2(2H+1)}+m_1\le\frac12+m_1,
\]
which gives the uniform constant $C_0$.
\end{proof}
\begin{remark}[finite-horizon lower and upper bounds]
Solving the quadratic inequality from
\cref{thm:past-finite-horizon} in $\sqrt{\kappa_H}\ge0$ gives
\[
 u_H^{\past}(L)\le L\kappa_H+2\sqrt{2\kappa_H}
 \qquad\iff\qquad
 \sqrt{\kappa_H}\ge
 \frac{\sqrt{2+Lu_H^{\past}(L)}-\sqrt2}{L}.
\]
Together with $\kappa_H\le u_H^{\past}(L)/L$, this yields
\begin{equation}
 \left(\frac{\sqrt{2+Lu_H^{\past}(L)}-\sqrt2}{L}\right)^2
 \le\kappa_H\le\frac{u_H^{\past}(L)}L.
 \label{eq:past-computable-bracket}
\end{equation}
A certified approximation of a finite-horizon value therefore gives
lower and upper bounds for the long-run constant.
\end{remark}
\paragraph{Prediction states.}
For the one-sided fractional Brownian motion with its natural filtration,
the observed past determines the prediction curve
\[
 G_t(r):=\E[B_{t+r}^H\mid\mathcal F_t^{B^H}],\qquad r\ge0.
\]
The conditional law of the future depends on the past only through
$(t,G_t)$. Its Markov property and its role in continuation values are
established in \cref{app:prediction-states}. A one-step tracking problem
there identifies information relevant to a decision beyond the current
scalar error.
\section{Local homogenization}
\label{sec:local-homogenization}
\subsection{A filtered fractional tangent}
In the rest of this section, the original probability space is standard up
to completion and each $\mathcal F_t$ is the completion of a countably
generated sigma-field.  These conditions hold on the canonical Gaussian path
spaces used below.
\begin{assumption}[filtered $L^2$ fractional tangent]
\label{ass:strong-tangent}
There are $H\in(0,1)$ and
$\rho\in\mathcal C([0,T];(0,\infty))$ such that the following three
conditions hold.
\smallskip
\noindent\textnormal{(FT1) Filtered coupling.}
For every compact $J\subset(0,T)$ and every
$0<c_-<c_+<\infty$, there are $h_*=h_*(J,c_-,c_+)>0$ and a common extension
$(\overline\Omega,\overline{\mathcal F},\overline\Pp)$ of the original
probability space.  If $\mathcal H\subset\overline{\mathcal F}$, write
$\overline{\mathcal H}$ for its $\overline\Pp$-completion.  On this
extension, for every $t\in J$, every
$0<h<h_*$, and every $c\in[c_-,c_+]$, there is an auxiliary random variable $U^{t,h,c}$ taking values in a standard Borel space, independent of the
pullback of $\mathcal F_T$, and a two-sided fBm $\mathbf B^{t,h,c}$ with
$\mathcal L(\mathbf B^{t,h,c})=\mu_H^{\past}$.  For every
$u\in[0,(T-t)/(ch)]$,
\begin{equation}
 \overline{\mathcal F_{t+chu}\vee\sigma(U^{t,h,c})}
 =
 \overline{\sigma(\mathbf B_v^{t,h,c}:v\le u)}.
 \label{eq:filtered-tangent-filtration}
\end{equation}
\smallskip
\noindent\textnormal{(FT2) Uniform $L^2$ tangent approximation.}
For every $R<\infty$,
\begin{equation}
 \sup_{\substack{t\in J,\ c\in[c_-,c_+]\\0\le u\le R}}
 \left\|
 \frac{\Theta_{t+chu}-\Theta_t}{\rho(t)(ch)^H}
 -\mathbf B_u^{t,h,c}
 \right\|_{L^2(\overline\Pp)}
 \xrightarrow[h\to0]{}0.
 \label{eq:strong-tangent-L2}
\end{equation}
\smallskip
\noindent\textnormal{(FT3) Global increment bound.}
For some $C,h_0>0$,
\begin{equation}
 \sup_{0<h<h_0}\sup_{0\le t\le T-h}
 h^{-2H}\E\left[|\Theta_{t+h}-\Theta_t|^2\right]
 \le C.
 \label{eq:strong-tangent-global-bound}
\end{equation}
\end{assumption}
Assume from now on that
$a,b,\rho\in\mathcal C([0,T];(0,\infty))$.  Define
$\alpha_H:=2H/(1+H)\in(0,1)$ and
\begin{equation}
 \appdef{g_H}{[0,T]}{(0,\infty)}{t}
 {a_t^{(1-H)/(1+H)}b_t^{2H/(1+H)}\rho(t)^{2/(1+H)}}.
 \label{eq:local-density}
\end{equation}
\begin{theorem}[local homogenization]
\label{thm:exact-local}
Suppose \cref{ass:strong-tangent} holds and $a,b,\rho$ are continuous and bounded away from zero.  Then
\begin{equation}
 \lambda^{-\alpha_H}V_\lambda(\Theta)
 \xrightarrow[\lambda\to0]{}
 \kappa_H\int_0^T g_H(t)\dd t.
 \label{eq:exact-local-homogenization}
\end{equation}
\end{theorem}
The proof is given below, after the reset-value estimates.
\subsection{Stability of reset values}
Fix $L>0$, write $\mathbf G=(\mathcal G_t)_{0\le t\le L}$, and let
$(\Omega,\mathcal F,\mathbf G,\Pp)$ be a complete filtered probability
space satisfying the usual conditions. Let
$X\in L^0(\Omega;\mathcal C([0,L];\mR))$ be $\mathbf G$-adapted, with
$X_0=0$ and $\E[\int_0^LX_t^2\dd t+|X_L|]<\infty$.
Let $\cA_L(\mathbf G)$ consist of the $\mathbf G$-adapted processes
$\psi\in L^0(\Omega;\mathcal D_L)$ such that $\psi_{0-}=0$,
$\TV_{(0,L]}(\psi)<\infty$ almost surely, and
\[
 \E\left[\int_0^L\psi_t^2\dd t+\CTV_{[0,L]}(\psi)+|\psi_L|\right]<\infty.
\]
For these controls, define the \emph{pathwise reset cost} by
\[
 \appdef{\mathfrak C_L^X}{\cA_L(\mathbf G)}{L^1(\Omega;[0,\infty))}{\psi}
 {\displaystyle\frac12\int_0^L(\psi_t-X_t)^2\dd t
 +\CTV_{[0,L]}(\psi)+|\psi_L-X_L|}.
\]
Its expectation defines
\[
 \appdef{J_L^X}{\cA_L(\mathbf G)}{[0,\infty)}{\psi}
 {\E[\mathfrak C_L^X(\psi)]},
\]
so the reset value is
\begin{equation}
 \mathcal U_L(X;\mathbf G)
 :=\inf_{\psi\in\cA_L(\mathbf G)}J_L^X(\psi).
 \label{eq:general-reset-value}
\end{equation}
When the filtration is clear, write $\mathcal U_L(X):=\mathcal U_L(X;\mathbf G)$.
\begin{lemma}[stability under a common filtration]
\label{lem:value-stability}
Let $X,Y$ start at zero and be adapted to the same filtration $\mathbf G$,
with finite reset values. Put
$A_L:=\E[\int_0^L|X_t-Y_t|^2\dd t]$ and
$B_L:=\E[|X_L-Y_L|]$. Then, for every $\delta>0$,
\begin{equation}
 \mathcal U_L(X;\mathbf G)
 \le(1+\delta)\mathcal U_L(Y;\mathbf G)
       +\tfrac12(1+\delta^{-1})A_L+B_L.
 \label{eq:value-stability}
\end{equation}
More precisely, writing $U_X=\mathcal U_L(X;\mathbf G)$ and
$U_Y=\mathcal U_L(Y;\mathbf G)$,
\begin{equation}
 |U_X-U_Y|\le B_L+\frac{A_L}{2}
                  +\sqrt{2A_L\min(U_X,U_Y)}.
 \label{eq:value-stability-quantitative}
\end{equation}
Therefore, if both values are $O(L)$ and
\begin{equation}
 A_L/L\xrightarrow[L\to\infty]{}0,\qquad B_L/L\xrightarrow[L\to\infty]{}0,
 \label{eq:average-coupling}
\end{equation}
then $|U_X-U_Y|=o(L)$.
\end{lemma}
\begin{proof}
See \cref{app:reset-stability}.
\end{proof}
\begin{lemma}[removing conditionally independent auxiliary variables]
\label{lem:independent-randomization}
Let $E$ and $S$ be standard Borel spaces.  Assume that the probability space
is standard up to completion and that each $\mathcal G_t$ is the completion
of a countably generated sigma-field.  Let
$Y\in L^0(\Omega,\mathcal G_0,\Pp;E)$ and
$U\in L^0(\Omega,\mathcal F,\Pp;S)$ satisfy
\begin{equation}
 U\indep\mathcal G_L\mid Y.
 \label{eq:conditional-seed-independence}
\end{equation}
Write $\mathcal G_t^U:=\overline{\mathcal G_t\vee\sigma(U)}^{\Pp}$
and let $\mathbf G^U:=(\mathcal G_t^U)_{0\le t\le L}$ be the enlarged
filtration. Then
\begin{equation}
 \mathcal U_L(X;\mathbf G^U)=\mathcal U_L(X;\mathbf G).
 \label{eq:conditional-seed-value}
\end{equation}
More precisely,
\begin{equation}
 \operatorname*{ess\,inf}_{\psi\in\cA_L(\mathbf G^U)}
 \E[\mathfrak C_L^X(\psi)\mid Y]
 =
 \operatorname*{ess\,inf}_{\psi\in\cA_L(\mathbf G)}
 \E[\mathfrak C_L^X(\psi)\mid Y]
 \qquad\Pp\text{-a.s.}
 \label{eq:conditional-seed-value-given-state}
\end{equation}
The conditional values at the initial time also agree:
\begin{equation}
 \operatorname*{ess\,inf}_{\psi\in\cA_L(\mathbf G^U)}
 \E[\mathfrak C_L^X(\psi)\mid\mathcal G_0^U]
 =\operatorname*{ess\,inf}_{\psi\in\cA_L(\mathbf G)}
 \E[\mathfrak C_L^X(\psi)\mid\mathcal G_0].
 \label{eq:conditional-seed-initial-value}
\end{equation}
The same statements hold for
the deterministic weighted tracking costs used in this paper and for a
pre-initial position measurable with respect to $\mathcal G_0$,
if some admissible control has finite cost.
\end{lemma}
\begin{proof}
See \cref{app:conditional-randomization}.
\end{proof}
\subsection{Uniform growing-cell approximation}
\begin{lemma}[uniform growing-cell approximation]
\label{lem:growing-cell-transfer}
Under \cref{ass:strong-tangent}, fix a compact
$J\subset(0,T)$ and $0<c_-<c_+<\infty$.  For admissible
$(t,h,c)$, define the local target
\[
 \appdef{X^{t,h,c}}{[0,(T-t)/(ch)]}{L^2(\Omega)}{u}
 {\dfrac{\Theta_{t+chu}-\Theta_t}{\rho(t)(ch)^H}},
\]
and define the local filtration
$\overline{\mathbf G}^{t,h,c}
=(\overline{\mathcal G}_u^{t,h,c})_{0\le u\le(T-t)/(ch)}$ by
\[
 \overline{\mathcal G}_u^{t,h,c}
 :=\overline{\mathcal F_{t+chu}\vee\sigma(U^{t,h,c})}.
\]
By \cref{lem:independent-randomization}, the value of the original target
is the same with or without this enlargement.  There is
$L_h\underset{h\to 0}{\uparrow}\infty$ such that
$hL_h\xrightarrow[h\to 0]{}0$ and, uniformly in $t\in J$ and
$c\in[c_-,c_+]$,
\begin{equation}
 \frac1{L_h}\mathcal U_{L_h}\left(
 \left(\frac{\Theta_{t+chu}-\Theta_t}{\rho(t)(ch)^H}\right)_{0\le u\le L_h};
 (\mathcal F_{t+chu})_{u\le L_h}
 \right)
 \xrightarrow[h\to 0]{}\kappa_H.
 \label{eq:growing-cell-transfer}
\end{equation}
\end{lemma}
\begin{proof}
See \cref{app:growing-cells}.
\end{proof}
\subsection{Proof of the local homogenization theorem}
The block reset decomposition, separator points, and short-interval
control used below are established in \cref{app:block-estimates}.
\begin{proof}[Proof of \cref{thm:exact-local}]
We first choose the local scales and establish a uniform block comparison.
For $f\in\mathcal C([0,T];\mR)$, write its \emph{modulus of continuity} as
\[
 \appdef{\omega_f}{[0,\infty)}{[0,\infty)}{r}
 {\displaystyle\sup_{\substack{s,t\in[0,T]\\|t-s|\le r}}|f(t)-f(s)|}.
\]
Put $h_\lambda:=\lambda^{1/(1+H)}$, so
$\lambda=h_\lambda^{1+H}$ and $\lambda^{\alpha_H}=h_\lambda^{2H}$, and define
\begin{equation}
 \appdef{\ell_\lambda}{[0,T]}{(0,\infty)}{t}
 {\left(\dfrac{\lambda b_t}{a_t\rho(t)}\right)^{1/(1+H)}}.
 \label{eq:local-microscopic-scale}
\end{equation}
\Cref{fig:local-scaling} illustrates this time and amplitude normalization.
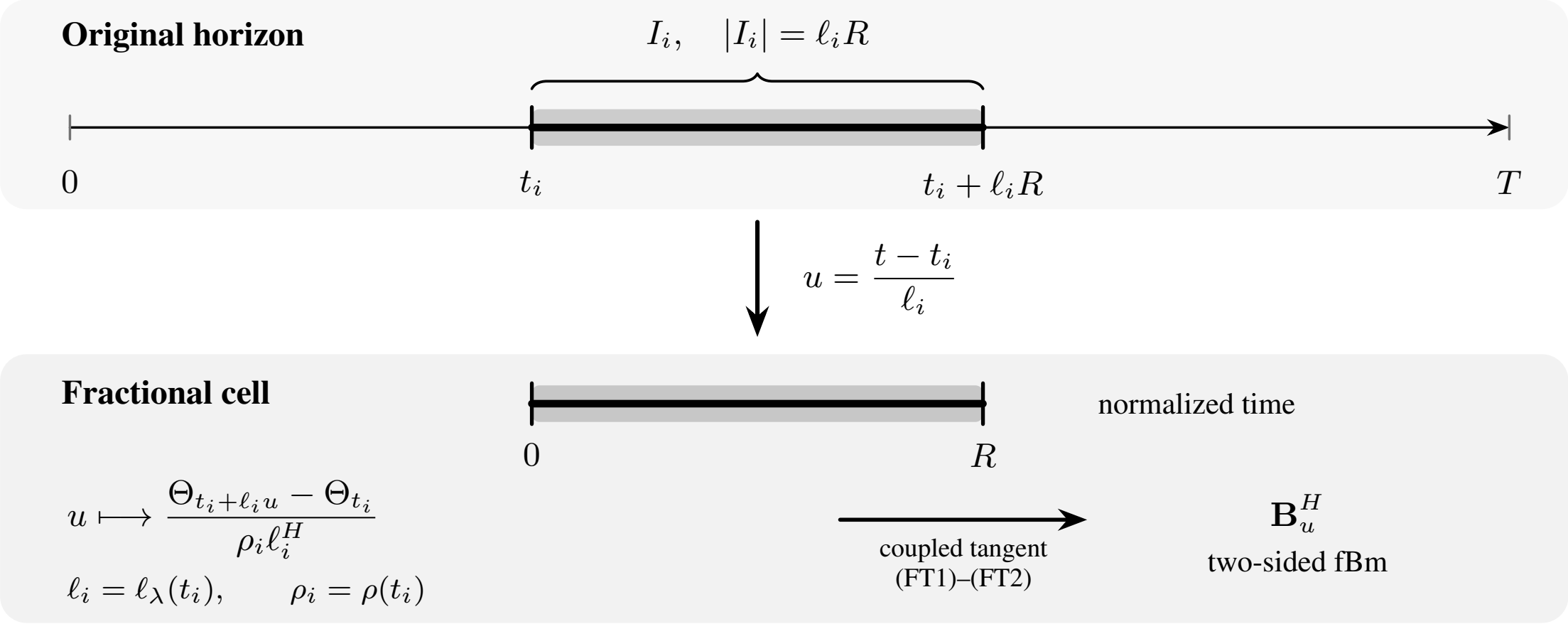
\begin{figure}[tbp]
\centering
\begin{tikzpicture}[
  x=1cm,y=1cm,
  font=\small,
  every node/.style={text=black},
  >=Stealth,
  line cap=round,
  line join=round
]
  \fill[rounded corners=7pt,black!3]
    (0,0) rectangle (14.6,1.95);
  \fill[rounded corners=7pt,black!5]
    (0,-3.85) rectangle (14.6,-1.35);

  \node[anchor=west,text=black,font=\bfseries\small]
    at (0.45,1.60) {Original horizon};
  \draw[black,line width=0.65pt,->] (0.65,0.76) -- (14.05,0.76);
  \fill[rounded corners=2pt,black!20]
    (4.95,0.59) rectangle (9.15,0.93);
  \draw[black,line width=2.0pt] (4.95,0.76) -- (9.15,0.76);
  \foreach \x in {0.65,14.05}
    \draw[black!55,line width=0.65pt] (\x,0.65) -- (\x,0.87);
  \foreach \x in {4.95,9.15}
    \draw[black,line width=1.0pt] (\x,0.57) -- (\x,0.95);
  \node[anchor=north] at (0.65,0.48) {$0$};
  \node[anchor=north] at (4.95,0.48) {$t_i$};
  \node[anchor=north] at (9.15,0.48) {$t_i+\ell_iR$};
  \node[anchor=north] at (14.05,0.48) {$T$};
  \draw[black,decorate,
        decoration={brace,amplitude=4pt},line width=0.7pt]
    (4.95,1.12) -- (9.15,1.12)
    node[midway,above=6pt,text=black] {$I_i,\quad |I_i|=\ell_iR$};

  \draw[black,line width=1.15pt,->]
    (7.05,-0.12) -- (7.05,-1.19);
  \node[anchor=west,text=black] at (7.35,-0.65)
    {$\displaystyle u=\frac{t-t_i}{\ell_i}$};

  \node[anchor=west,text=black,font=\bfseries\small]
    at (0.45,-1.70) {Fractional cell};
  \fill[rounded corners=2pt,black!22]
    (4.95,-1.98) rectangle (9.15,-1.64);
  \draw[black,line width=2.0pt] (4.95,-1.81) -- (9.15,-1.81);
  \foreach \x in {4.95,9.15}
    \draw[black,line width=1.0pt] (\x,-2.00) -- (\x,-1.62);
  \node[anchor=north] at (4.95,-2.06) {$0$};
  \node[anchor=north] at (9.15,-2.06) {$R$};
  \node[anchor=west,text=black,font=\footnotesize]
    at (10.10,-1.81) {normalized time};

  \node[anchor=west] at (0.50,-2.89)
    {$\displaystyle
       u\longmapsto
       \frac{\Theta_{t_i+\ell_i u}-\Theta_{t_i}}
            {\rho_i\ell_i^H}$};
  \draw[black,line width=1.05pt,->]
    (7.82,-2.89) -- (10.12,-2.89);
  \node[text=black,font=\scriptsize,align=center]
    at (8.97,-3.31) {coupled tangent\\(FT1)--(FT2)};
  \node[text=black] at (12.08,-2.83)
    {$\displaystyle\mathbf B_u^H$};
  \node[text=black,font=\footnotesize]
    at (12.08,-3.28) {two-sided fBm};
  \node[anchor=west,text=black,font=\footnotesize]
    at (0.50,-3.55)
    {$\ell_i=\ell_\lambda(t_i),\qquad\rho_i=\rho(t_i)$};
\end{tikzpicture}
\space
\caption{Local rescaling of a tracking block. The physical interval
$I_i=[t_i,t_i+\ell_iR]$ becomes $[0,R]$, and target increments are
normalized by $\rho_i\ell_i^H$. The lower arrow represents the filtered
$L^2$ coupling with the two-sided fBm tangent.}
\label{fig:local-scaling}
\end{figure}
Positivity and continuity give constants $0<c_-<c_+<\infty$ with
\[
 \frac{\ell_\lambda(t)}{h_\lambda}
 =\left(\frac{b_t}{a_t\rho(t)}\right)^{1/(1+H)}\in[c_-,c_+].
\]
Fix $\eta\in(0,T/2)$ and apply \cref{lem:growing-cell-transfer} on
$J_\eta:=[\eta/2,T-\eta/2]$. Write $L_\lambda:=L_{h_\lambda}$ and
$\Delta_\lambda:=h_\lambda L_\lambda$. Then
\[
 L_\lambda\xrightarrow[\lambda\to0]{}\infty,\qquad\Delta_\lambda\xrightarrow[\lambda\to0]{}0.
\]
The extension to comparable horizons is proved in \cref{app:growing-cells}.
Uniformly for starts in $J_\eta$ and horizons
$R=qL_\lambda+O(1)$, with $q$ bounded above and away from zero,
\[
 \mathcal U_R\left(
 \left(\frac{\Theta_{s+\ell_\lambda(s)u}-\Theta_s}
 {\rho(s)\ell_\lambda(s)^H}\right)_{0\le u\le R};
 (\mathcal F_{s+\ell_\lambda(s)u})_{0\le u\le R}\right)
 =(\kappa_H+o(1))R.
\]
To bound the effect of holding the coefficients constant, put $C_*:=2(1\vee c_+)$
and set
\[
 \omega_\lambda:=\max_{f\in\{a,b,\rho\}}
 \frac{\omega_f(C_*\Delta_\lambda)}{\min_{[0,T]}f}\xrightarrow[\lambda\to0]{}0.
\]
Holding the coefficients constant on a block of length at most
$C_*\Delta_\lambda$ changes every nonnegative cost term by a relative
factor between $1-\omega_\lambda$ and $1+\omega_\lambda$. 
To obtain the upper bound, start at $s_0:=\eta$ and take consecutive
blocks $I_i:=[s_i,s_{i+1}]\subset[\eta,T-\eta]$ of length
$\ell_iL_\lambda$, where $\ell_i:=\ell_\lambda(s_i)$.
Let $M_\lambda$ be the number of these blocks, indexed by
$i=0,\ldots,M_\lambda-1$.
For each block, write $a_i:=a_{s_i}$, $b_i:=b_{s_i}$ and
$\rho_i:=\rho(s_i)$. There are $O(\Delta_\lambda^{-1})$ blocks, and the
interior remainder has length $O(\Delta_\lambda)$.
With $t=s_i+\ell_i u$ and
$\phi_t=\Theta_{s_i}+\rho_i\ell_i^H\psi_u$, the reset cost
$J_{i,\lambda}^{\mathrm{fix}}$ with coefficients held constant satisfies
\begin{equation}
\begin{aligned}
 \inf_\psi J_{i,\lambda}^{\mathrm{fix}}(\psi)
 =a_i\rho_i^2\ell_i^{2H+1}\mathcal U_{L_\lambda}\left(
 \left(\frac{\Theta_{s_i+\ell_i u}-\Theta_{s_i}}
 {\rho_i\ell_i^H}\right)_{0\le u\le L_\lambda}\right).
\end{aligned}
 \label{eq:block-scaling}
\end{equation}
Indeed, $\ell_i^{H+1}=\lambda b_i/(a_i\rho_i)$, so the running and
variation/reset coefficients agree:
\[
 a_i\rho_i^2\ell_i^{2H+1}=\lambda b_i\rho_i\ell_i^H.
\]
Moreover,
\begin{equation}
 a_i\rho_i^2\ell_i^{2H}
 =\lambda^{\alpha_H}
 a_i^{1-\alpha_H}b_i^{\alpha_H}\rho_i^{2-\alpha_H}
 =\lambda^{\alpha_H}g_H(s_i).
 \label{eq:block-density-identity}
\end{equation}
Since $|I_i|=\ell_iL_\lambda$, the uniform cell comparison gives
\[
 \inf_\psi J_{i,\lambda}^{\mathrm{fix}}(\psi)
 =(\kappa_H+o(1))\lambda^{\alpha_H}g_H(s_i)|I_i|.
\]
Choose original-filtration controls within $\delta_\lambda L_\lambda$
of the normalized infima, where $\delta_\lambda\underset{\lambda\to0}{\to}0$.
\Cref{lem:independent-randomization} removes the auxiliary enlargement.
Their total approximation error is
\begin{align*}
 \sum_{i=0}^{M_\lambda-1} a_i\rho_i^2\ell_i^{2H+1}\delta_\lambda L_\lambda
 =\delta_\lambda\lambda^{\alpha_H}\sum_{i=0}^{M_\lambda-1} g_H(s_i)|I_i|
=O(\delta_\lambda\lambda^{\alpha_H})=o(\lambda^{\alpha_H}).
\end{align*}
The cost with the original coefficients is at most $1+\omega_\lambda$ times the cost with coefficients held constant.
Terminal resets and initial atoms cover the cost of all jumps between blocks
(\cref{lem:canonical-concatenation}). Therefore the interior cost is at most
\[
 (\kappa_H+o(1))\lambda^{\alpha_H}\sum_{i=0}^{M_\lambda-1} g_H(s_i)|I_i|.
\]
The Riemann sums satisfy
\[
 \sum_{i=0}^{M_\lambda-1} g_H(s_i)|I_i|\xrightarrow[\lambda\to0]{}
 \int_\eta^{T-\eta}g_H(t)\dd t.
\]
Use \cref{lem:short-interval-rule} on $[0,\eta]$, $[T-\eta,T]$, and
the interior remainder. Their combined cost is at most
$C\lambda^{\alpha_H}(2\eta+O(\Delta_\lambda)+3h_\lambda)$, whence
\[
 \limsup_{\lambda\to0}\lambda^{-\alpha_H}V_\lambda(\Theta)
 \le\kappa_H\int_\eta^{T-\eta}g_H(t)\dd t+C\eta.
\]
For the lower bound, we first control the extra cost of separating blocks.
Take $\phi^\lambda$ with cost
$V_\lambda(\Theta)+o(\lambda^{\alpha_H})$. The upper bound and
$a\ge a_->0$ give
\begin{equation}
\begin{aligned}
 \E\left[\int_0^T|\phi_t^\lambda-\Theta_t|^2\dd t\right]
 \le\frac2{a_-}J_\lambda(\phi^\lambda;\Theta)
 \le C\lambda^{\alpha_H}=Ch_\lambda^{2H}.
\end{aligned}
 \label{eq:near-optimal-energy}
\end{equation}
For sufficiently small $\lambda$, set
$K_\lambda:=\lfloor(T-2\eta-h_\lambda)/\Delta_\lambda\rfloor$
and use grid points $z_j:=\eta+j\Delta_\lambda$,
$j=0,\ldots,K_\lambda$, so that
$J_j:=[z_j,z_j+h_\lambda]\subset[\eta,T-\eta]$.
For small $\lambda$, these intervals are disjoint. Their number is
$O(\Delta_\lambda^{-1})$.
By \cref{lem:separator-points}, choose deterministic $r_j\in J_j$ such that
\begin{align}
 \sum_{j=0}^{K_\lambda}\E\left[|\phi_{r_j}^\lambda-\Theta_{r_j}|\right]
 &\overset{\eqref{eq:separator-sum}}\le
 \left(\frac{K_\lambda+1}{h_\lambda}
 \E\left[\int_0^T|\phi_t^\lambda-\Theta_t|^2\dd t\right]\right)^{1/2}\notag\\
 &\overset{\eqref{eq:near-optimal-energy}}\le
 C\left(\frac{h_\lambda^{2H}}
 {h_\lambda^2L_\lambda}\right)^{1/2}
 =C\frac{h_\lambda^{H-1}}{\sqrt{L_\lambda}}.
 \label{eq:boundary-error-sum}
\end{align}
Let $I_i:=[r_{i-1},r_i]$, $i=1,\ldots,K_\lambda$. Add a virtual initial adjustment from
$\Theta_{r_{i-1}}$ and a terminal reset to $\Theta_{r_i}$.
By \eqref{eq:weighted-virtual-reset-decomposition}, their total excess
cost is at most
\begin{equation}
\begin{aligned}
 2\lambda b_+\sum_{j=0}^{K_\lambda}\E\left[|\phi_{r_j}^\lambda-\Theta_{r_j}|\right]
 \overset{\eqref{eq:boundary-error-sum}}=
 O\left(\frac{h_\lambda^{1+H}h_\lambda^{H-1}}
 {\sqrt{L_\lambda}}\right)
 =O\left(\frac{h_\lambda^{2H}}{\sqrt{L_\lambda}}\right)
 =o(\lambda^{\alpha_H}).
\end{aligned}
 \label{eq:virtual-reset-negligible}
\end{equation}
It remains to compare the values of the resulting blocks.
Since $r_j\in[z_j,z_j+h_\lambda]$ and $z_i-z_{i-1}=\Delta_\lambda$,
\[
 |I_i|=\Delta_\lambda+O(h_\lambda)
 =h_\lambda L_\lambda+O(h_\lambda).
\]
For each block, set $t_i:=r_{i-1}$ and $\ell_i^-:=\ell_\lambda(t_i)$.
Its normalized horizon is
\[
 R_{i,\lambda}:=\frac{|I_i|}{\ell_i^-}
 =\frac{h_\lambda}{\ell_i^-}L_\lambda+O(1),
\]
with $h_\lambda/\ell_i^-\in[1/c_+,1/c_-]$. Let $J_i^\lambda$ be the expected weighted reset cost when restricting
$\phi^\lambda$ to $I_i$. Holding the coefficients constant and applying the
cell estimate for comparable horizons give, uniformly in $i$,
\begin{align*}
 J_i^\lambda
 &\ge
 (1-\omega_\lambda)a_{t_i}\rho(t_i)^2(\ell_i^-)^{2H+1}
 \times\mathcal U_{R_{i,\lambda}}\left(
 \left(\frac{\Theta_{t_i+\ell_i^-u}-\Theta_{t_i}}
 {\rho(t_i)(\ell_i^-)^H}\right)_{0\le u\le R_{i,\lambda}}\right)\\
 &\ge
 (\kappa_H-o(1))a_{t_i}\rho(t_i)^2(\ell_i^-)^{2H+1}R_{i,\lambda}\\
 &=
 (\kappa_H-o(1))a_{t_i}\rho(t_i)^2(\ell_i^-)^{2H}|I_i|\\
 &\overset{\eqref{eq:block-density-identity}}=
 (\kappa_H-o(1))\lambda^{\alpha_H}g_H(t_i)|I_i|.
\end{align*}
All discarded original costs are nonnegative. Hence
\[
 J_\lambda(\phi^\lambda;\Theta)
 \ge\sum_{i=1}^{K_\lambda}J_i^\lambda
 -2\lambda b_+\sum_{j=0}^{K_\lambda}\E\left[|\phi_{r_j}^\lambda-\Theta_{r_j}|\right].
\]
The selected blocks leave $O(\Delta_\lambda+h_\lambda)$ of
$[\eta,T-\eta]$ uncovered, and their mesh tends to zero. Therefore,
\[
 \sum_{i=1}^{K_\lambda} g_H(t_i)|I_i|\xrightarrow[\lambda\to0]{}
 \int_\eta^{T-\eta}g_H(t)\dd t.
\]
Together with \eqref{eq:virtual-reset-negligible}, this yields
\[
 \liminf_{\lambda\to0}\lambda^{-\alpha_H}V_\lambda(\Theta)
 \ge\kappa_H\int_\eta^{T-\eta}g_H(t)\dd t.
\]
Let $\eta\to0$ in both bounds and use continuity of $g_H$.
\end{proof}
\subsection{The cost split of near-optimal controls}
Set $p_H=(1-H)/(1+H)$, so that $p_H+\alpha_H=1$.
\begin{theorem}[local split between tracking and variation costs]
\label{thm:cost-allocation}
Under the assumptions of \cref{thm:exact-local}, let $\phi^\lambda$ satisfy
$J_\lambda(\phi^\lambda;\Theta)
 =V_\lambda(\Theta)+o(\lambda^{\alpha_H})$.
Define the finite tracking measure on $[0,T]$ by
\begin{equation}
 \mu_\lambda(\dd t):=\lambda^{-\alpha_H}\frac{a_t}{2}
       \E[(\phi_t^\lambda-\Theta_t)^2]\dd t,
 \label{eq:tracking-measure}
\end{equation}
and define the variation measure, for $A\in\mathcal B([0,T])$, by
\begin{equation}
 \nu_\lambda(A):=\lambda^{1-\alpha_H}
       \E\left[\int_A b_t\,|\dd\phi^\lambda|_t\right].
 \label{eq:variation-measure}
\end{equation}
Then, weakly as finite measures,
\begin{equation}
\begin{aligned}
 \mu_\lambda\underset{\lambda\to0}{\Longrightarrow}
 p_H\kappa_Hg_H(t)\dd t,\qquad\text{and}\qquad
 \nu_\lambda\underset{\lambda\to0}{\Longrightarrow}
 \alpha_H\kappa_Hg_H(t)\dd t.
\end{aligned}
 \label{eq:local-cost-allocation}
\end{equation}
In particular, the leading fractions of the total cost are
$p_H$ for tracking and $\alpha_H$ for variation.
\end{theorem}
\begin{proof}
Fix $f\in\mathcal C([0,T];\mR)$ and
$|\varepsilon|\|f\|_\infty<1$. Let $V_\lambda(\varepsilon)$ be the
value with $a$ replaced by $a(1+\varepsilon f)$. Testing with
$\phi^\lambda$ gives
\[
 V_\lambda(\varepsilon)-V_\lambda(0)
 \le o(\lambda^{\alpha_H})
 +\varepsilon\lambda^{\alpha_H}\int_{[0,T]}f\dd\mu_\lambda.
\]
The perturbed coefficients satisfy the same assumptions, so
\[
 \lambda^{-\alpha_H}V_\lambda(\varepsilon)
 \xrightarrow[\lambda\to0]{}K(\varepsilon),
\]
where
\[
 \appdef{K}{\{\varepsilon\in\mR:|\varepsilon|\|f\|_\infty<1\}}{\mR}
 {\varepsilon}
 {\displaystyle\kappa_H\int_0^Tg_H(t)(1+\varepsilon f(t))^{p_H}\dd t}.
\]
For fixed $\varepsilon$, division by
$\varepsilon\lambda^{\alpha_H}$ gives the following lower and upper bounds.
For $\varepsilon>0$,
\begin{equation*}
 \frac{K(\varepsilon)-K(0)}{\varepsilon}
 \le\liminf_{\lambda\to0}\int_{[0,T]}f\dd\mu_\lambda,
\end{equation*}
and, for $\varepsilon<0$,
\begin{equation*}
 \limsup_{\lambda\to0}\int_{[0,T]}f\dd\mu_\lambda
 \le\frac{K(\varepsilon)-K(0)}{\varepsilon}.
\end{equation*}
For $\varepsilon$ near zero, differentiation under the integral gives:
\[
 K'(\varepsilon)
 =p_H\kappa_H\int_0^Tf(t)g_H(t)(1+\varepsilon f(t))^{p_H-1}\dd t.
\]
Then,
\[
 K'(0)=p_H\kappa_H\int_0^Tf(t)g_H(t)\dd t.
\]
Letting $\varepsilon\to0$ from either side proves the tracking limit.
For the perturbation $b(1+\varepsilon f)$, the test term is
$\varepsilon\lambda^{\alpha_H}\int_{[0,T]}f\dd\nu_\lambda$ and the limiting
value is
\[
 \kappa_H\int_0^Tg_H(t)(1+\varepsilon f(t))^{\alpha_H}\dd t.
\]
The same argument, with exponent $\alpha_H$ in place of $p_H$, proves
\[
 \int_{[0,T]}f\dd\nu_\lambda\xrightarrow[\lambda\to0]{}
 \alpha_H\kappa_H\int_0^Tf(t)g_H(t)\dd t.
\]
Taking $f\equiv1$ gives the total cost fractions.
\end{proof}
\begin{remark}
For Brownian targets the split is $1/3$ for tracking and $2/3$ for
variation. Formula \eqref{eq:local-cost-allocation} also rules out a
leading-order cost atom at either endpoint for near-optimal controls.
\end{remark}
\subsection{A Gaussian Volterra class}
Fix $H\in(0,1)$ and let $W$ be a standard Brownian motion
on $[0,T]$, with augmented filtration $\mathcal F_t:=\mathcal F_t^W$.
For $F,q:[0,T]\to\mR$ and $K:(0,T]\to\mR$, define
\begin{equation}
 \appdef{\Theta}{[0,T]}{L^2(\Omega)}{t}
 {\displaystyle F(t)+\int_0^tK(t-s)q(s)\dd W_s}.
 \label{eq:volterra-target}
\end{equation}
Suppose that $q\in\mathcal C([0,T];\mR)$ and
\begin{equation}
 0<q_-\le|q(t)|\le q_+\qquad(0\le t\le T).
 \label{eq:q-nonzero}
\end{equation}
The drift is assumed to satisfy $F\in\mathcal C^{\alpha_F}([0,T];\mR)$
for some $\alpha_F\in(H,1]$. For the kernel, assume
$K\in\mathcal C^1((0,T];\mR)\cap L^2((0,T))$ and
\begin{equation}
 K(u)=c_Ku^{H-1/2}(1+o(1))\quad\text{as }u\to0,
 \qquad c_K\ne0.
 \label{eq:kernel-asymptotics}
\end{equation}
Finally, require the bounds
\begin{equation}
 |K(u)|\le Cu^{H-1/2}\quad\text{and}\quad
 |K'(u)|\le Cu^{H-3/2}\qquad(0<u\le T).
 \label{eq:kernel-bounds}
\end{equation}
Set
\[
 \Gamma_H:=\frac1{2H}
 +\int_0^\infty\left((1+z)^{H-1/2}-z^{H-1/2}\right)^2\dd z.
\]
With this normalization, define the tangent amplitude by
\begin{equation}
 \appdef{\rho}{[0,T]}{(0,\infty)}{t}{|c_Kq(t)|\sqrt{\Gamma_H}}.
 \label{eq:Gamma-rho}
\end{equation}
The integral defining $\Gamma_H$ is finite because its integrand is
$O(1+z^{2H-1})$ at zero and $O(z^{2H-3})$ at infinity.  When $H=1/2$, the integrand is set equal to zero at every $z>0$. Thus $\Gamma_{1/2}=1$.
Put $\beta:=H-1/2$. The \emph{Mandelbrot--Van Ness
kernel} \cite{mandelbrot1968} is
\begin{equation}
 \appdef{\mathfrak g_u^H}{\mR}{\mR}{z}
 {(u-z)_+^\beta-(-z)_+^\beta},\qquad u\in\mR.
 \label{eq:mvn-kernel}
\end{equation}
Here $x_+^\beta=x^\beta$ for $x>0$ and is zero for $x\le0$,
including $x_+^0=\ind_{\{x>0\}}$. Thus, for $H=1/2$,
\[
 \mathfrak g_u^{1/2}=
 \begin{cases}\ind_{(0,u)},&u>0,\\0,&u=0,\\-\ind_{(u,0)},&u<0,\end{cases}
 \quad\text{as elements of }L^2(\mR).
\]
\begin{proposition}[filtered tangent for Gaussian Volterra targets]
\label{prop:volterra-strong-tangent}
For every $H\in(0,1)$, the process in
\eqref{eq:volterra-target}--\eqref{eq:kernel-bounds} has a continuous
adapted modification satisfying \cref{ass:strong-tangent} with $\rho$ in
\eqref{eq:Gamma-rho}.  Therefore,
\begin{equation}
\begin{aligned}
\lambda^{-2H/(1+H)}V_\lambda(\Theta)\xrightarrow[\lambda\to0]{}
 \kappa_H\int_0^T
 a_t^{(1-H)/(1+H)}b_t^{2H/(1+H)}
 |c_Kq(t)|^{2/(1+H)}\Gamma_H^{1/(1+H)}\dd t.
\end{aligned}
 \label{eq:volterra-exact-equivalent}
\end{equation}
\end{proposition}
The coupling extends the driving Brownian motion by an independent
past and applies the Mandelbrot--Van Ness kernel. The kernel estimates
then give the uniform $L^2$ approximation and the global increment bound.
\begin{proof}
See \cref{app:volterra-verification}.
\end{proof}
\begin{corollary}[H\"older drifts do not affect the leading term]
\label{cor:regular-drift}
In the class of \cref{prop:volterra-strong-tangent}, the leading
asymptotic formula \eqref{eq:volterra-exact-equivalent} does not depend on the
drift $F$, if its H\"older exponent is larger than $H$.
\end{corollary}
\begin{proof}
For every fixed $R<\infty$,
\[
 \sup_{0\le u\le R}h^{-H}|F(t+hu)-F(t)|
 \le[F]_{\alpha_F}R^{\alpha_F}h^{\alpha_F-H}\xrightarrow[h\to0]{}0.
\]
Thus $F$ vanishes from the filtered tangent and does not enter $\rho$ or
the limiting density.
\end{proof}
\begin{remark}
At $t=0$ the tangent is generally Riemann--Liouville rather than
stationary fBm. The endpoint is handled by the short-interval estimate
\eqref{eq:short-interval-cost}: an interval of length $\eta$ contributes at
most $C\eta$ to the normalized upper bound, uniformly as
$\lambda\to0$.
\end{remark}
\section{Bounds and properties of the fractional cell constant}
\label{sec:bounds}
\subsection{A causal dual lower bound}
\begin{proposition}[causal dual lower bound]
\label{prop:causal-dual}
Fix $L>0$ and a complete filtered probability space
$(\Omega,\mathcal F,\mathbf F,\Pp)$ with
$\mathbf F=(\mathcal F_t)_{0\le t\le L}$.  Let
$X$ be a continuous, $\mathbf F$-adapted real process on $[0,L]$,
satisfying $X_0=0$ and the integrability assumptions required for
$\mathcal U_L(X;\mathbf F)$. Let $W^{1,2}(0,L)$ denote the space of
functions $f\in L^2(0,L)$ with weak derivative $f'\in L^2(0,L)$,
and let $p\in L^2(\Omega;W^{1,2}(0,L))$ have an almost surely absolutely
continuous version with $|p_t|\le1$. The field $p$ may be nonadapted,
and $\dot p\in L^2(\Omega\times[0,L],\Pp\otimes\dd t)$.  Write
\[
 L^2_{\mathrm{prog}}(\Omega\times[0,L])
 :=L^2(\Omega\times[0,L],\mathrm{Prog}(\mathbf F),
 \Pp\otimes\dd t),
\]
where $\mathrm{Prog}(\mathbf F)$ is the progressive sigma-field of
$\mathbf F$.  Let $m$ be the orthogonal projection of $\dot p$ onto
this closed subspace.  Equivalently, after choosing a jointly measurable
version,
\begin{equation}
 m_t=\E[\dot p_t\mid\mathcal F_t]
 \quad\text{for }\Pp\otimes\dd t\text{-a.e. }(\omega,t).
 \label{eq:dual-projection}
\end{equation}
Then $m\in L^2_{\mathrm{prog}}(\Omega\times[0,L])$ and
\begin{equation}
 \mathcal U_L(X;\mathbf F)
 \ge
 \E\left[p_LX_L-\int_0^L X_tm_t\dd t
 -\frac12\int_0^L m_t^2\dd t\right].
 \label{eq:causal-dual}
\end{equation}
\end{proposition}
\begin{proof}
For $\psi\in\cA_L(\mathbf F)$, the polar inequality and $|p|\le1$ give
\[
 \CTV_{[0,L]}(\psi)\ge\int_{[0,L]}p_t\dd\psi_t,
 \qquad |\psi_L-X_L|\ge-p_L(\psi_L-X_L).
\]
Integration by parts, including the initial atom $\psi_0-\psi_{0-}$,
uses $\psi_{0-}=0$ and continuity of $p$:
\begin{align*}
 \int_{[0,L]}p_t\dd\psi_t
 &=p_0\psi_0+\int_{(0,L]}p_t\dd\psi_t\\
 &=
 p_0\psi_0+p_L\psi_L-p_0\psi_0-\int_0^L\psi_t\dot p_t\dd t\\
 &=p_L\psi_L-\int_0^L\psi_t\dot p_t\dd t.
\end{align*}
The adapted c\`adl\`ag process $\psi$ is progressive and belongs to
$L^2(\Omega\times[0,L])$. Orthogonality of the projection therefore gives
\[
 \E\left[\int_0^L\psi_t(\dot p_t-m_t)\dd t\right]=0.
\]
The terminal terms cancel except for $p_LX_L$, and completion of the
square yields
\begin{align*}
 J_L^X(\psi)
 &\ge
 \E\left[p_LX_L+\int_0^L
       \left\{\frac12(\psi_t-X_t)^2-\psi_tm_t\right\}\dd t\right]\\
 &=
 \E\left[p_LX_L-\int_0^LX_tm_t\dd t-\frac12\int_0^Lm_t^2\dd t\right]+\frac12\E\left[\int_0^L(\psi_t-X_t-m_t)^2\dd t\right]\\
 &\ge\E\left[p_LX_L-\int_0^LX_tm_t\dd t-\frac12\int_0^Lm_t^2\dd t\right].
\end{align*}
Taking infima proves \eqref{eq:causal-dual}.
\end{proof}
\begin{remark}
If $p$ is adapted and $\dot p$ is progressive, then $m=\dot p$.  If $p$ is
not adapted, Jensen gives $\E[m_t^2]\le\E[\dot p_t^2]$ and
$\E[X_tm_t]=\E[X_t\dot p_t]$.
\end{remark}
\begin{corollary}[ergodic causal dual bound]
\label{cor:ergodic-dual}
On the two-sided fBm space of \cref{def:past-cells}, let
$\mathfrak P_L^H$ be the set of all
$p\in L^2(\Omega;W^{1,2}(0,L))$ with an almost surely absolutely continuous
version satisfying $|p_t|\le1$ for every $t\in[0,L]$. For each such $p$, put
$m_t^p:=\E[\dot p_t\mid\mathcal F_t^{\mathbf B}]$ and define the functional
\[
 \appdef{\mathfrak D_{H,L}}{\mathfrak P_L^H}{\mR}{p}
 {\E\left[p_L\mathbf B_L^H-\int_0^L
       \left(\mathbf B_t^Hm_t^p+\tfrac12(m_t^p)^2\right)\dd t\right]}.
\]
Then, for every $L>0$,
\begin{equation}
 u_H^{\past}(L)\ge\sup_{p\in\mathfrak P_L^H}\mathfrak D_{H,L}(p),
 \label{eq:finite-horizon-dual-cell}
\end{equation}
and
\begin{equation}
 \kappa_H\ge
 \limsup_{L\to\infty}\frac1L
 \sup_{p\in\mathfrak P_L^H}\mathfrak D_{H,L}(p).
 \label{eq:ergodic-dual-cell}
\end{equation}
\end{corollary}
\begin{proof}
Applying \cref{prop:causal-dual} with
$(X,\mathbf F)=(\mathbf B^H,\mathbf F^{\mathbf B})$ gives
$\mathfrak D_{H,L}(p)\le u_H^{\past}(L)$ for every
$p\in\mathfrak P_L^H$. Taking the supremum proves the finite-horizon
claim and, after division by $L$, yields
\[
 \frac1L\sup_{p\in\mathfrak P_L^H}\mathfrak D_{H,L}(p)
 \le\frac{u_H^{\past}(L)}L.
\]
The upper limit and \cref{thm:past-robustness} then give the ergodic bound.
\end{proof}
\subsection{Causal linear-filter bounds}
\begin{proposition}[single-window bounds]
\label{prop:analytic-bounds}
For every $H\in(0,1)$,
\begin{equation}
 C_{\rm low}(H)\le\kappa_H\le C_{\rm ma}(H),
 \label{eq:analytic-bracket}
\end{equation}
where
\begin{equation}
 \appdef{C_{\rm low}}{(0,1)}{(0,\infty)}{H}
 {\frac{m_1(1+H)}2
  \left(\frac{m_1(1-H)}4\right)^{(1-H)/(1+H)}},
 \label{eq:C-low}
\end{equation}
and
\begin{equation}
 \appdef{C_{\rm ma}}{(0,1)}{(0,\infty)}{H}
 {\frac1{4(1-H)}
  \left(\frac{2m_1(1-H)(1+H)}H\right)^{2H/(1+H)}}.
 \label{eq:C-ma}
\end{equation}
\end{proposition}
\begin{proof}
See \cref{app:single-window-bounds}.
\end{proof}
\begin{theorem}[continuity of the fractional cell constant]
\label{thm:kappa-continuity}
The map $H\mapsto\kappa_H$ is continuous on $(0,1)$.
For each fixed $L>0$, $H\mapsto u_H^{\past}(L)$ is also continuous.
\end{theorem}
\begin{proof}
Fix $H_0\in(0,1)$ and a compact interval
$I=[H_-,H_+]\subset(0,1)$ containing $H_0$ in its interior. Couple
all $H\in I$ on one two-sided Brownian space by
\[
 B_t^H:=\Gamma_H^{-1/2}\int_\mR\mathfrak g_t^H(s)\dd W_s.
\]
For fixed $t$, the kernels converge pointwise almost everywhere as
$H\to H_0$. Near $s=0,t$, their squares have an integrable envelope
of the form
\[
 C_t\bigl(1+|s|^{2H_--1}+|s-t|^{2H_--1}\bigr),
 \qquad 2H_--1>-1.
\]
For sufficiently negative $s$, they are bounded by
$C_t|s|^{2H_+-3}$, with $2H_+-3<-1$.
Dominated convergence gives
\[
 \|\mathfrak g_t^H-\mathfrak g_t^{H_0}\|_2\xrightarrow[H\to H_0]{}0,\qquad
 \Gamma_H=\|\mathfrak g_1^H\|_2^2\xrightarrow[H\to H_0]{}\Gamma_{H_0}>0.
\]
The Wiener isometry and the triangle inequality therefore give
\begin{align*}
 \|B_t^H-B_t^{H_0}\|_{L^2}
 &=\|\Gamma_H^{-1/2}\mathfrak g_t^H
       -\Gamma_{H_0}^{-1/2}\mathfrak g_t^{H_0}\|_2\\
 &\le
 \Gamma_H^{-1/2}\|\mathfrak g_t^H-\mathfrak g_t^{H_0}\|_2
 +|\Gamma_H^{-1/2}-\Gamma_{H_0}^{-1/2}|
              \|\mathfrak g_t^{H_0}\|_2
 \xrightarrow[H\to H_0]{}0.
\end{align*}
On $[0,L]$,
\[
 \E\left[|B_t^H-B_t^{H_0}|^2\right]\le2t^{2H}+2t^{2H_0}
 \le4\bigl(t^{2H_-}\ind_{\{t\le1\}}+t^{2H_+}\ind_{\{t>1\}}\bigr).
\]
Another application of dominated convergence yields
\begin{equation*}
 \|B^H-B^{H_0}\|_{L^2(\Omega\times[0,L])}
 \xrightarrow[H\to H_0]{}0,
\end{equation*}
and
\begin{equation*}
 \E\left[|B_L^H-B_L^{H_0}|\right]
 \xrightarrow[H\to H_0]{}0.
\end{equation*}
By \cref{lem:mvn-filtration}, all coupled targets have the same
completed full-past filtration. The common-filtration stability
estimate in \cref{lem:value-stability} therefore gives
\begin{align*}
 |u_H^{\past}(L)-u_{H_0}^{\past}(L)|
 &\le
   \E\left[|B_L^H-B_L^{H_0}|\right]
   +\tfrac12\|B^H-B^{H_0}\|_{L^2(\Omega\times[0,L])}^2\\
 &\hspace{1em}+\sqrt{2u_{H_0}^{\past}(L)}
       \|B^H-B^{H_0}\|_{L^2(\Omega\times[0,L])}\\
 &\xrightarrow[H\to H_0]{}0.
\end{align*}
Finally, \cref{thm:past-finite-horizon} gives
\[
 |\kappa_H-\kappa_{H_0}|
 \le\frac{2C_0}{L}+
 \frac{|u_H^{\past}(L)-u_{H_0}^{\past}(L)|}{L}.
\]
First let $H\to H_0$, then $L\to\infty$.
\end{proof}
\begin{corollary}[endpoint behavior of $\kappa_H$]
\label{cor:kappa-endpoints}
The fractional cell constant satisfies
\begin{equation}
 \lim_{H\uparrow1}\kappa_H=\sqrt{\frac2\pi}.
 \label{eq:kappa-H-one}
\end{equation}
At the other endpoint,
\begin{equation}
 \frac1{4\pi}
 \le
 \liminf_{H\to0}\kappa_H
 \le
 \limsup_{H\to0}\kappa_H
 \le
 \frac14.
 \label{eq:kappa-H-zero}
\end{equation}
\end{corollary}
\begin{proof}
Put $\varepsilon:=1-H$. Since $\varepsilon\log\varepsilon\to0$,
we have
\begin{equation*}
 C_{\rm low}(1-\varepsilon)
 =\frac{m_1(2-\varepsilon)}2
   \exp\left[\frac{\varepsilon}{2-\varepsilon}
                  \log\left(\frac{m_1\varepsilon}{4}\right)\right]
 \xrightarrow[\varepsilon\to0]{} m_1,
\end{equation*}
and
\begin{equation*}
 \begin{aligned}
 C_{\rm ma}(1-\varepsilon)
 &=\frac{m_1(2-\varepsilon)}{2(1-\varepsilon)}
   \exp\left[-\frac{\varepsilon}{2-\varepsilon}
     \log\left(\frac{2m_1\varepsilon(2-\varepsilon)}{1-\varepsilon}\right)\right]\\
 &\xrightarrow[\varepsilon\to0]{} m_1.
 \end{aligned}
\end{equation*}
At the other endpoint,
\[
 C_{\rm low}(H)\xrightarrow[H\to0]{}\frac{m_1}{2}\frac{m_1}{4}
 =\frac1{4\pi},
\]
and
\[
 \log C_{\rm ma}(H)
 =-\log(4(1-H))+\frac{2H}{1+H}
   \log\left(\frac{2m_1(1-H)(1+H)}H\right)
 \xrightarrow[H\to0]{}-\log4,
\]
since $H\log H\xrightarrow[H\to0]{}0$. The bounds in \eqref{eq:analytic-bracket} give
both conclusions.
\end{proof}
The moving-average upper bound extends to linear combinations of
windows. The resulting bound, its scale-free form, and the optimal
weights for fixed windows are given in \cref{app:multiwindow-bounds}.
\section{Finite-dimensional Markov approximations}
\label{sec:factor-convergence}
Assume $H\in(0,1/2)$ and put $\beta:=\frac12-H\in(0,1/2)$.
\emph{Euler's gamma function} is
\[
 \appdef{\Gamma}{(0,\infty)}{(0,\infty)}{z}
 {\int_0^\infty r^{z-1}e^{-r}\dd r}.
\]
Define the \emph{fractional kernel} by
\begin{equation}
 \appdef{k_H}{(0,\infty)}{(0,\infty)}{t}
 {\sqrt{2H}\,t^{H-1/2}
  =\frac{\sqrt{2H}}{\Gamma(\beta)}
    \int_0^\infty e^{-xt}x^{\beta-1}\dd x}.
 \label{eq:fractional-laplace}
\end{equation}
The substitution $r=xt$ gives
\[
 \int_0^\infty e^{-xt}x^{\beta-1}\dd x
 =t^{-\beta}\Gamma(\beta),\qquad
 \Gamma(\beta+1)=\beta\Gamma(\beta).
\]
For a function $k$ on $(0,\infty)$, define its \emph{causal extension} by
\[
 \appdef{k^\circ}{\mR}{\mR}{r}
 {\begin{cases}k(r),&r>0,\\0,&r\le0.\end{cases}}
\]
Following the Ornstein--Uhlenbeck mixture approach
\cite{carmona2000,harms2021}, we approximate this kernel by positive
exponential sums. For $n\in\Nn$, let $N_n\in\Nn$ and define
\begin{equation}
 \appdef{k_n}{(0,\infty)}{(0,\infty)}{t}
 {\sum_{i=1}^{N_n}w_{i,n}e^{-x_{i,n}t}},
 \qquad w_{i,n},x_{i,n}>0.
 \label{eq:positive-exponential-kernel}
\end{equation}
On a two-sided Brownian space, define
\[
 A_H:=\sqrt{2H\Gamma_H},
 \qquad
 Y_t^{i,n}:=\int_{-\infty}^te^{-x_{i,n}(t-s)}\dd W_s,
\]
and
\begin{equation}
 \begin{aligned}
 Z_t^n&:=\frac1{A_H}\sum_{i=1}^{N_n}w_{i,n}(Y_t^{i,n}-Y_0^{i,n})\\
 &=\frac1{A_H}\int_{\mR}
       \{k_n^\circ(t-s)-k_n^\circ(-s)\}\dd W_s.
 \end{aligned}
 \label{eq:factor-target}
\end{equation}
The coupled two-sided fBm is
\begin{equation}
 B_t^H
 :=\frac1{A_H}\int_{\mR}
 \{k_H^\circ(t-s)-k_H^\circ(-s)\}\dd W_s,
 \qquad t\in\mR.
 \label{eq:coupled-fbm-kernel}
\end{equation}
Let $W_0=0$ and let
$Y_t^n:=(Y_t^{1,n},\ldots,Y_t^{N_n,n})\in\mR^{N_n}$.  Define the
completed full-past filtrations
\[
 \mathcal F_t^W:=\overline{\sigma(W_s:s\le t)}^{\Pp},\qquad
 \mathcal F_t^{Y^n}:=\overline{\sigma(Y_s^n:s\le t)}^{\Pp},
 \quad t\in\mR,
\]
and write
$\mathbf F^W=(\mathcal F_t^W)_{t\ge0}$ and
$\mathbf F^{Y^n}=(\mathcal F_t^{Y^n})_{t\ge0}$ for control times.
The process at negative times represents the past observed at time zero.
\begin{lemma}[factor and Brownian filtrations are equal]
\label{lem:redundant-brownian-past}
For every $t\ge0$,
\begin{equation}
 \mathcal F_t^{Y^n}=\mathcal F_t^W.
 \label{eq:factor-brownian-filtrations}
\end{equation}
Therefore, for every $L>0$,
\begin{equation}
 \mathcal U_L(Z^n;\mathbf F^W)
 =
 \mathcal U_L(Z^n;\mathbf F^{Y^n}).
 \label{eq:redundant-brownian-past}
\end{equation}
\end{lemma}
\begin{proof}
Each factor is a Wiener integral of the Brownian past, so
$\mathcal F_t^{Y^n}\subset\mathcal F_t^W$ for $t\ge0$.
Conversely, fix any component $i$. Its OU equation is
\[
 \dd Y_s^{i,n}=-x_{i,n}Y_s^{i,n}\dd s+\dd W_s.
\]
Integrating from zero to $s\le t$ gives
\[
 W_s=Y_s^{i,n}-Y_0^{i,n}+x_{i,n}\int_0^sY_r^{i,n}\dd r.
\]
Because $t\ge0$, the observed factor past contains the whole interval
between $s$ and zero, including $Y_0^{i,n}$. Thus it recovers every
$W_s$, $s\le t$, and the completed filtrations coincide. The control
classes, and hence their reset values, are identical.
\end{proof}
For $L>0$, define
\begin{equation}
 \appdef{u_n}{(0,\infty)}{[0,\infty)}{L}
 {\mathcal U_L(Z^n;\mathbf F^W)}.
 \label{eq:factor-cell-value}
\end{equation}
\begin{theorem}[quantitative convergence of factor cell values]
\label{thm:factor-value-convergence}
For $n\in\Nn$ and $L\ge1$, measure the approximation error by
\[
 A_{n,L}:=\E\left[\int_0^L|Z_t^n-B_t^H|^2\dd t\right],\qquad
 B_{n,L}:=\E[|Z_L^n-B_L^H|].
\]
There is a finite constant $C$, independent of $n,L$, such that
\begin{equation}
 \left|\frac{u_n(L)}L-\kappa_H\right|
 \le\frac{C}{L}+\frac{B_{n,L}}L
       +\frac{A_{n,L}}{2L}
       +C\sqrt{\frac{A_{n,L}}L}.
 \label{eq:factor-quantitative-value-bound}
\end{equation}
In particular, if $n_k,L_k\xrightarrow[k\to\infty]{}\infty$ and
\begin{equation}
 \frac1{L_k}\E\left[\int_0^{L_k}|Z_t^{n_k}-B_t^H|^2\dd t\right]
 +\frac1{L_k^2}\E[|Z_{L_k}^{n_k}-B_{L_k}^H|^2]\xrightarrow[k\to\infty]{}0,
 \label{eq:factor-value-criterion}
\end{equation}
then
\begin{equation}
 \frac{u_{n_k}(L_k)}{L_k}\xrightarrow[k\to\infty]{}\kappa_H.
 \label{eq:factor-value-limit}
\end{equation}
\end{theorem}
\begin{proof}
By \cref{lem:mvn-filtration}, the fBm reset value in $\mathbf F^W$ is
$u_H^{\past}(L)$. Stability in this common filtration gives
\[
 |u_n(L)-u_H^{\past}(L)|
 \le B_{n,L}+\frac{A_{n,L}}2+
             \sqrt{2u_H^{\past}(L)A_{n,L}}.
\]
By \cref{thm:past-finite-horizon}, for $L\ge1$,
\begin{equation*}
 0\le u_H^{\past}(L)-\kappa_HL\le2\sqrt{2\kappa_H},
\end{equation*}
so
\begin{equation*}
 u_H^{\past}(L)\le(\kappa_H+2\sqrt{2\kappa_H})L.
\end{equation*}
The triangle inequality combines these bounds into
\begin{align*}
 \left|\frac{u_n(L)}L-\kappa_H\right|
 \le
   \frac{2\sqrt{2\kappa_H}}L+\frac{B_{n,L}}L+\frac{A_{n,L}}{2L}+\sqrt{2(\kappa_H+2\sqrt{2\kappa_H})}
                     \sqrt{\frac{A_{n,L}}L}.
\end{align*}
This is \eqref{eq:factor-quantitative-value-bound} with a constant
independent of $n,L$.
Under \eqref{eq:factor-value-criterion}, Cauchy--Schwarz also gives
\begin{equation*}
 \frac{B_{n_k,L_k}}{L_k}
 \le
 \left(\frac{\E[|Z_{L_k}^{n_k}-B_{L_k}^H|^2]}{L_k^2}\right)^{1/2}
 \xrightarrow[k\to\infty]{}0,
\end{equation*}
and
\begin{equation*}
 \frac{A_{n_k,L_k}}{L_k}\xrightarrow[k\to\infty]{}0.
\end{equation*}
Every term in the bound vanishes.
\end{proof}
\begin{corollary}[kernel criterion]
\label{cor:factor-kernel-criterion}
For $n\in\Nn$, put $\delta k_n:=k_n-k_H$ on $(0,\infty)$ and, for $t\ge0$,
define
\begin{equation}
 \begin{aligned}
 \varepsilon_n(t)^2
 :=\frac1{A_H^2}\int_0^t|\delta k_n(r)|^2\dd r+\frac1{A_H^2}\int_0^\infty
       |\delta k_n(t+r)-\delta k_n(r)|^2\dd r.
 \end{aligned}
 \label{eq:factor-kernel-error}
\end{equation}
Then
\begin{equation}
 \E[|Z_t^n-B_t^H|^2]=\varepsilon_n(t)^2.
 \label{eq:factor-kernel-isometry}
\end{equation}
Hence \eqref{eq:factor-value-limit} holds for any
$n_k\xrightarrow[k\to\infty]{}\infty$ and $L_k\xrightarrow[k\to\infty]{}\infty$ such that
\begin{equation}
 \frac1{L_k}\int_0^{L_k}\varepsilon_{n_k}(t)^2\dd t
 +\frac{\varepsilon_{n_k}(L_k)^2}{L_k^2}
 \xrightarrow[k\to\infty]{}0.
 \label{eq:deterministic-factor-criterion}
\end{equation}
\end{corollary}
\begin{proof}
The difference is the Wiener integral
\[
 Z_t^n-B_t^H=\frac1{A_H}\int_\mR
 \{\delta k_n^\circ(t-s)-\delta k_n^\circ(-s)\}\dd W_s.
\]
The integrand vanishes for $s\ge t$. On $0<s<t$ it equals
$\delta k_n(t-s)$, whereas on $s\le0$ it equals
$\delta k_n(t-s)-\delta k_n(-s)$. The Wiener isometry, followed by
the substitutions $r=t-s$ and $r=-s$ in the two integrals, therefore gives
\begin{align*}
 A_H^2\E\left[|Z_t^n-B_t^H|^2\right]
 &=
   \int_0^t|\delta k_n(t-s)|^2\dd s
   +\int_{-\infty}^0|\delta k_n(t-s)-\delta k_n(-s)|^2\dd s\\
 &=
   \int_0^t|\delta k_n(r)|^2\dd r
   +\int_0^\infty|\delta k_n(t+r)-\delta k_n(r)|^2\dd r\\
 &=A_H^2\varepsilon_n(t)^2.
\end{align*}
The second integral controls the stationary past. No global
$L^2(0,\infty)$ assumption on $\delta k_n$ is needed.
Substitution in \eqref{eq:factor-value-criterion} proves the claim.
\end{proof}
\begin{theorem}[a positive quadrature]
\label{cor:factor-diagonal}
Fix $H\in(0,1/2)$ and write $\beta=1/2-H$.
Given $0<\varepsilon<R$ and $0<\eta\le1$, choose a partition
$\varepsilon=x_0<\cdots<x_N=R$ into intervals $[x_j,x_{j+1}]$,
$j=0,\ldots,N-1$, with
$x_{j+1}/x_j\le1+\eta$. Define
\begin{equation}
 \appdef{k_{\varepsilon,R,\eta}}{(0,\infty)}{(0,\infty)}{t}
 {\sum_{j=0}^{N-1}w_j e^{-x_jt}},
 \label{eq:constructive-quadrature}
\end{equation}
where
\[
 w_j:=\frac{\sqrt{2H}}{\Gamma(\beta+1)}
           (x_{j+1}^{\beta}-x_j^{\beta})>0.
\]
Let $Z^{\varepsilon,R,\eta}$ be the target built from this kernel
using stationary factors, with initial value zero. There is $C_H<\infty$ such that, for $L\ge1$,
\begin{equation}
 \sup_{0\le t\le L}
 \|Z_t^{\varepsilon,R,\eta}-B_t^H\|_{L^2}
 \le C_H\{\sqrt L\,\varepsilon^\beta+R^{-H}+\eta L^H\}.
 \label{eq:constructive-kernel-bound}
\end{equation}
In particular, fix $\gamma>0$ and, for $L\ge2$, choose
\begin{equation}
 \varepsilon_L=L^{-(\gamma+1/2)/\beta},\qquad
 R_L=L^{\gamma/H},\qquad \eta_L=L^{-(H+\gamma)}.
 \label{eq:constructive-resolution}
\end{equation}
A geometric partition then uses
$N_L=O_{H,\gamma}(L^{H+\gamma}\log L)$ factors with positive weights, and its cell value
satisfies
\begin{equation}
 \left|\frac{u_{N_L}(L)}L-\kappa_H\right|
 \le C_{H,\gamma}\{L^{-\gamma}+L^{-1}\}.
 \label{eq:constructive-value-rate}
\end{equation}
Here $u_{N_L}$ denotes the cell value built from this quadrature with
$N_L$ factors.
\end{theorem}
\begin{proof}
See \cref{app:positive-quadrature}.
\end{proof}
\begin{proposition}[the small-cost limit for fixed factors]
\label{prop:fixed-factor-brownian}
Fix $n$, put $D_n:=A_H^{-1}\sum_{i=1}^{N_n}w_{i,n}>0$, and define
\[
 C_n:=\frac{m_1T}{\sqrt2 A_H}\sum_{i=1}^{N_n}w_{i,n}\sqrt{x_{i,n}}.
\]
Write $V_\lambda^{(n)}(T)$ for the tracking value of $Z^n$ with unit
weights and observation filtration $\mathbf F^{Y^n}$, without a terminal reset.
Write $V_\lambda^B(T;D_n)$ for the value of the Brownian target
$D_nW$ with unit weights and its natural filtration. Then
\begin{equation}
 |V_\lambda^{(n)}(T)-V_\lambda^B(T;D_n)|\le\lambda C_n.
 \label{eq:fixed-factor-BV-error}
\end{equation}
Therefore,
\begin{equation}
 V_\lambda^{(n)}(T)
 \underset{\lambda\to0}{\sim}\kappa_{1/2}T D_n^{4/3}\lambda^{2/3}.
 \label{eq:fixed-factor-brownian-rate}
\end{equation}
\end{proposition}
\begin{proof}
Integrating the OU equations and using $W_0=0$ gives
\[
 Y_t^{i,n}-Y_0^{i,n}=W_t-x_{i,n}\int_0^tY_s^{i,n}\dd s.
\]
Therefore
\[
 Z_t^n=D_nW_t+A_t^n,\qquad
 A_t^n:=-\frac1{A_H}\sum_{i=1}^{N_n}w_{i,n}x_{i,n}\int_0^tY_s^{i,n}\dd s.
\]
This adapted perturbation is absolutely continuous and starts at zero.
Stationarity and the Wiener isometry give
\begin{align*}
 \E[(Y_s^{i,n})^2]&=\int_0^\infty e^{-2x_{i,n}r}\dd r
 =\frac1{2x_{i,n}}.
\end{align*}
Then,
\begin{align*}
 \E\left[|Y_s^{i,n}|\right]&=\frac{m_1}{\sqrt{2x_{i,n}}}.
\end{align*}
The triangle inequality and these stationary Gaussian moments now yield
\begin{align*}
 \E\left[\TV_{(0,T]}(A^n)\right]
 &=\E\left[\int_0^T\left|\frac1{A_H}\sum_{i=1}^{N_n}w_{i,n}x_{i,n}Y_s^{i,n}\right|\dd s\right]\\
 &\le
   \frac1{A_H}\sum_{i=1}^{N_n}w_{i,n}x_{i,n}
                       \int_0^T\E\left[|Y_s^{i,n}|\right]\dd s\\
 &=
   \frac{m_1T}{\sqrt2A_H}\sum_{i=1}^{N_n}w_{i,n}\sqrt{x_{i,n}}=C_n.
\end{align*}
Apply \cref{prop:BV-perturbation} in the full-past filtration.
The negative Brownian path is independent of the future increments.
\Cref{lem:independent-randomization} removes it from the Brownian value.
Thus \eqref{eq:fixed-factor-BV-error} follows, and
\[
 \lambda^{-2/3}|V_\lambda^{(n)}(T)-V_\lambda^B(T;D_n)|
 \le C_n\lambda^{1/3}\xrightarrow[\lambda\to0]{}0.
\]
For fixed $n$, \cref{thm:cell,prop:brownian} give the asymptotic formula.
\end{proof}
\begin{remark}[different small-cost exponents]
For every fixed $n$, \cref{prop:fixed-factor-brownian} gives
\[
  \lim_{\lambda\to0}
  \frac{\log V_\lambda^{(n)}(T)}{\log\lambda}=\frac23,
\]
whereas \cref{thm:cell} gives exponent $2H/(1+H)$ for fBm.  The joint limit
in \cref{thm:factor-value-convergence} couples resolution with horizon
growth.
\end{remark}
\begin{proposition}[Markov representation of the factor cell]
\label{prop:factor-Markov-reduction}
Let
\begin{align*}
 \mathcal F_t^{n,+}
 &:=\overline{\sigma\bigl(Y_0^n,W_s-W_0:0\le s\le t\bigr)}^{\Pp}\\
 &=\overline{\sigma(Y_s^n:0\le s\le t)}^{\Pp},
 \qquad 0\le t\le L.
\end{align*}
The second equality follows from the factor SDE and the reconstruction
formula used in \cref{lem:redundant-brownian-past}.  Then
\begin{equation}
 \mathcal U_L(Z^n;\mathbf F^{Y^n})
 =\mathcal U_L(Z^n;\mathbf F^{n,+}).
 \label{eq:factor-Markov-reduction}
\end{equation}
Thus the reset value depends on the past before zero only through $Y_0^n$.
\end{proposition}
\begin{proof}
The factor solution is
\[
 Y_t^{i,n}=e^{-x_{i,n}t}Y_0^{i,n}
          +\int_0^te^{-x_{i,n}(t-s)}\dd W_s,\qquad t\ge0.
\]
Thus $Z^n|_{[0,L]}$ is adapted to $\mathbf F^{n,+}$.
Let $U:=(Y_s^n)_{s\le0}$, viewed in
$\mathcal C_{\loc}((-\infty,0];\mR^{N_n})$.
Future Brownian increments are independent of the whole initial past,
so
\[
 U\indep\mathcal F_L^{n,+}\mid Y_0^n,
 \qquad
 \mathcal F_t^{Y^n}=\overline{\mathcal F_t^{n,+}\vee\sigma(U)}^{\Pp}.
\]
Apply \cref{lem:independent-randomization} with the initial state
$Y=Y_0^n$ to obtain \eqref{eq:factor-Markov-reduction}.
\end{proof}
The controlled Markov state lies in $\mR\times\mR^{N_n}$.
Introduce the coefficient vector
\[
 a:=(w_{1,n}/A_H,\ldots,w_{N_n,n}/A_H)^\top\in(0,\infty)^{N_n}
\]
and the diagonal matrix
\[
 X_n:=\diag(x_{1,n},\ldots,x_{N_n,n})\in\mR^{N_n\times N_n}.
\]
Here $\mathbf1\in\mR^{N_n}$ is the all-ones vector and
$D_n=\mathbf1^\top a$. Since
$Z_t^n=a^\top(Y_t^n-Y_0^n)$, the error $e:=\psi-Z^n$ satisfies
\begin{equation}
 \dd e_t=\dd\psi_t+a^\top X_nY_t^n\dd t-D_n\dd W_t,
\end{equation}
and
\begin{equation}
 \dd Y_t^n=-X_nY_t^n\dd t+\mathbf1\dd W_t.
 \label{eq:factor-controlled-state}
\end{equation}
The initial vector $Y_0^n$ is centered Gaussian, and the Wiener
isometry gives its covariance:
\[
 \E[Y_0^{i,n}Y_0^{j,n}]
 =
 \int_{-\infty}^0e^{(x_{i,n}+x_{j,n})s}\dd s
 =\frac1{x_{i,n}+x_{j,n}}.
\]
For $f\in\mathcal C^2(\mR^{N_n+1};\mR)$, the generator between
interventions is
\begin{equation}
 \begin{aligned}
 \mathcal L_nf(e,y)
 &=(a^\top X_ny)\partial_ef-(X_ny)\cdot\nabla_yf+\frac12
       \left(-D_n\partial_e+\sum_{i=1}^{N_n}\partial_{y_i}\right)^2f\\
 &=(a^\top X_ny)\partial_ef-(X_ny)\cdot\nabla_yf
       +\frac{D_n^2}{2}\partial_{ee}f
 -D_n\sum_{i=1}^{N_n}\partial_{ey_i}f
       +\frac12\sum_{i=1}^{N_n}\sum_{j=1}^{N_n}\partial_{y_iy_j}f.
 \end{aligned}
 \label{eq:factor-generator}
\end{equation}
The running, variation, and terminal costs are $e^2/2$, one per unit
adjustment, and $|e|$.
\paragraph{A quadratic-hedging application.}
Under the assumptions of \cref{sec:hedging}, the
Galtchouk--Kunita--Watanabe decomposition identifies the penalized
mean-square hedging value, after subtracting the irreducible term,
with $V_\lambda(\Theta)$. The identity and its small-cost consequence
are given in \cref{sec:hedging}.
\section{Conclusion}
\label{sec:conclusion}
The increment scale determines the exponent, while a causal fractional
cell determines the constant. We use resets to obtain subadditivity and finite-memory approximation
to show that the initial past has no long-run effect on the value. Averaging over shifted partitions then yields a
bounded full-past correction and continuity in $H$. Coupling local paths with their filtrations lets us transfer the cell
value to
nonstationary targets. Homogeneity then determines the near-optimal cost split. Finally, quadrature with positive weights approximates the cell
when resolution grows with the horizon, whereas fixed factors retain a
Brownian small-cost limit.
\appendix
\section{Causal deconvolution of the Mandelbrot--Van Ness map}
\label{app:mvn-deconvolution}
\begin{lemma}[Brownian and fractional pasts generate the same sigma-field]
\label{lem:mvn-filtration}
Fix $H\in(0,1)$.  Let
$W=(W_t)_{t\in\mR}$ be a two-sided standard Brownian motion with $W_0=0$,
and define
\[
 \mathbf B_t^H:=\frac1{\sqrt{\Gamma_H}}
 \int_{\mR}\mathfrak g_t^H(s)\dd W_s,
 \qquad
 \mathfrak g_t^H(s):=(t-s)_+^{H-1/2}-(-s)_+^{H-1/2},
\]
with the order-zero convention from \eqref{eq:mvn-kernel}.  For $t\in[0,\infty)$, set
\[
 \mathcal H_t^W
 :=\overline{\Span}^{L^2(\Omega)}\{W_s:s\le t\},
 \qquad
 \mathcal H_t^{B}
 :=\overline{\Span}^{L^2(\Omega)}\{\mathbf B_s^H:s\le t\}.
\]
The Gaussian spans and completed pasts coincide:
\begin{equation}
 \mathcal H_t^W=\mathcal H_t^{B},
 \qquad
 \overline{\sigma(W_s:s\le t)}^{\Pp}
 =\overline{\sigma(\mathbf B_s^H:s\le t)}^{\Pp}.
 \label{eq:mvn-filtration-equality}
\end{equation}
The result also holds after shifting the time origin, subtracting the
value at that origin, and scaling time, for nonnegative relative control times.
\end{lemma}
\begin{proof}
To compare the two pasts, put $\beta:=H-1/2$ and introduce the \emph{oriented
Brownian interval indicator}, for $s\in\mR$, by
\[
 \iota_s(r):=
 \begin{cases}\ind_{(0,s]}(r),&s\ge0,\\-\ind_{(s,0]}(r),&s<0\end{cases}.
\]
Then,
\[
 W_s=\int_\mR\iota_s(r)\dd W_r.
\]
For $t\ge0$, these indicators with $s\le t$ span a dense subspace of
$L^2((-\infty,t])$, so the Wiener isometry gives
\[
 \mathcal H_t^W
 =\left\{\int_\mR f(r)\dd W_r:
     f\in L^2(\mR),\ f=0\ \text{a.e. on }(t,\infty)\right\}.
\]
Since $\mathfrak g_s^H$ vanishes on $(s\vee0,\infty)$,
$\mathbf B_s^H\in\mathcal H_t^W$ for every $s\le t$, hence
$\mathcal H_t^B\subset\mathcal H_t^W$. For the converse, we work in the deterministic Wiener-integrand space.
Starting from a step function $f=\sum_{j=0}^{M-1} a_j\ind_{(r_j,r_{j+1}]}$,
where $M\in\Nn$ and $r_0<\cdots<r_M$, define
\[
 I_B(f):=\sum_{j=0}^{M-1} a_j(\mathbf B_{r_{j+1}}^H-\mathbf B_{r_j}^H),
 \qquad \|f\|_{\mathfrak H_H}^2:=\E[I_B(f)^2].
\]
Let $\mathfrak H_H$ be the Hilbert completion after identifying functions
of zero seminorm. The map $I_B:\mathfrak H_H\to L^2(\Omega)$ is its isometric
extension. The Wiener-integral deconvolution identity of
\cite{pipiras2002}, in the present normalization, is
\begin{equation}
 W_s=\frac{\sqrt{\Gamma_H}}{\Gamma(1+\beta)\Gamma(1-\beta)}
 \int_\mR\bigl\{(s-r)_+^{-\beta}-(-r)_+^{-\beta}\bigr\}
                         \dd\mathbf B_r^H.
 \label{eq:mvn-inverse}
\end{equation}
This is an $L^2$ Wiener integral. The inverse kernel is supported on
$(-\infty,s\vee0]$, and its step approximations $f_m$ use only that
half-line. By the Wiener isometry,
\[
\left\|W_s-\frac{\sqrt{\Gamma_H}\,I_B(f_m)}{\Gamma(1+\beta)\Gamma(1-\beta)}\right\|_{L^2}^2
=\frac{\Gamma_H}{\Gamma(1+\beta)^2\Gamma(1-\beta)^2}
\underbrace{\left\|(s-\cdot)_+^{-\beta}-(-\cdot)_+^{-\beta}-f_m\right\|_{\mathfrak H_H}^2}_{\displaystyle\xrightarrow[m\to\infty]{}0}.
\]
For $s\le t$ and $t\ge0$, each $I_B(f_m)$ is a finite linear combination
of fractional observations at times at most $t$. Thus
$W_s\in\mathcal H_t^B$. The normalization in
\eqref{eq:mvn-inverse} comes from composing the right-sided fractional
operators of orders $\beta$ and $-\beta$, whose kernel factors are
$\Gamma(1+\beta)$ and $\Gamma(1-\beta)$. At $\beta=0$,
$\Gamma_{1/2}=1$ and the identity reduces to the Brownian interval integral.
The two inclusions show that the closed Gaussian spans coincide.
Since an $L^2$ limit has an almost-surely convergent subsequence,
completing a Gaussian linear span does not change its completed
generated sigma-field. This proves
both equalities in \eqref{eq:mvn-filtration-equality}.
The same support argument applies after shifting the time origin,
subtracting the value at that origin, and rescaling, for
nonnegative relative times.
\end{proof}
\section{Conditional removal of auxiliary variables}
\label{app:conditional-randomization}
\begin{proof}[Proof of \cref{lem:independent-randomization}]
\emph{Conditional averaging.}
Let $\psi\in\cA_L(\mathbf G^U)$ have finite cost. To average out the
auxiliary variable, represent the countable base generators by a
standard Borel variable $Z$ and choose a Borel path map
$\psi=\Psi(Z,U)$. Conditional independence then gives
$\mathcal L(U\mid\mathcal G_L)=\eta_Y$, where
$\eta_y:=\mathcal L(U\mid Y=y)$ is a regular conditional distribution
\cite{kallenberg2021}. We therefore define
\[
 \bar\psi_t:=\int_S\Psi(Z,u)_t\,\eta_Y(\dd u).
\]
As $\psi_{0-}=0$,
\[
 \sup_{t\le L}|\psi_t|\le\CTV_{[0,L]}(\psi),
\qquad\text{and}\qquad
 \E\left[\CTV_{[0,L]}(\psi)\right]<\infty.
\]
Conditional Fubini supplies one integrable pathwise dominating function
for the average. Dominated convergence gives right-continuity and left
limits of $\bar\psi$. For every finite partition $\pi$ of $[0,L]$
whose interior points are rational,
\[
 |\bar\psi_0|+\sum_{[s,t]\in\pi}|\bar\psi_t-\bar\psi_s|
 \le\int_S\left(|\Psi(Z,u)_0|
       +\sum_{[s,t]\in\pi}|\Psi(Z,u)_t-\Psi(Z,u)_s|\right)\eta_Y(\dd u).
\]
Taking the supremum over these partitions yields
\[
 \CTV_{[0,L]}(\bar\psi)
 \le\int_S\CTV_{[0,L]}(\Psi(Z,u))\eta_Y(\dd u).
\]
The same argument on subintervals gives domination of the variation
measures and hence of their integrals against deterministic nonnegative
Borel weights. 
Set $Q_L:=(\Qq\cap[0,L))\cup\{L\}$.
For $q\in Q_L$, adaptedness represents $\Psi(Z,U)_q$ using only
$\mathcal G_q$ and $U$. Since $Y$ is $\mathcal G_0$-measurable,
integration against $\eta_Y$ makes $\bar\psi_q$ $\mathcal G_q$-measurable.
A common null set suffices for all $q\in Q_L$.
Right-continuity of paths and filtration extends adaptedness to
every $t<L$, while the terminal time is included in $Q_L$.
\begin{samepage}
Jensen also gives
\[
 \E\left[\int_0^L\bar\psi_t^2\dd t\right]\le\E\left[\int_0^L\psi_t^2\dd t\right],
\]
and
\[
 \E\left[|\bar\psi_L|\right]\le\E\left[|\psi_L|\right].
\]
\end{samepage}
Thus $\bar\psi\in\cA_L(\mathbf G)$.
The target is $\mathcal G_L$-measurable. Conditional Jensen for the
quadratic and absolute-value terms, together with the variation bound,
yields
\[
 \mathfrak C_L^X(\bar\psi)
 \le\int_S\mathfrak C_L^X(\Psi(Z,u))\eta_Y(\dd u)
 =\E[\mathfrak C_L^X(\psi)\mid\mathcal G_L].
\]
Taking expectations proves the reverse inequality in
\eqref{eq:conditional-seed-value}. Conditioning on $Y$ proves the reverse
inequality in \eqref{eq:conditional-seed-value-given-state}.
The other inequalities follow from inclusion of the control classes.\\
\emph{Initial conditional values.}
Let $V$ be the base conditional value given $\mathcal G_0$, and $W$ the
enlarged conditional value given $\mathcal G_0^U$.
For integrable $\mathcal G_L$-measurable $C$,
\[
 \E[C\mid\mathcal G_0^U]=\E[C\mid\mathcal G_0],
\]
so testing enlarged values with base controls gives $W\le V$.
The conditional cost families are downward directed: paste two controls
on the event, known initially, where the first conditional cost is smaller.
They therefore admit decreasing sequences approaching their essential
infima, bounded above by the conditional cost of a fixed finite-cost
control. Dominated convergence then gives
\[
 \E[W]=\mathcal U_L(X;\mathbf G^U)
       =\mathcal U_L(X;\mathbf G)=\E[V].
\]
Since $0\le W\le V$ and both are integrable, $W=V$ almost surely.
For a $\mathcal G_0$-measurable pre-initial position $x_0$, average
$\psi-x_0$ and add $x_0$ back. The same convexity arguments apply to the
weighted costs whenever a finite-cost admissible control exists.
\end{proof}
\section{Path-space rules and concatenation}
\label{app:canonical}
Fix $L>0$. As above, write
$\mathcal C_L=\mathcal C([0,L];\mR)$ and
$\mathcal D_L=\mathcal D([0,L];\mR)$, and put
$Q_L:=(\Qq\cap[0,L))\cup\{L\}$.
These path spaces are Polish under the uniform and Skorokhod $J_1$
topologies, respectively \cite[Chapters~2--3]{billingsley1999}.
For $\mu\in\mathcal P(\mathcal C_L)$, introduce
the coordinate maps
\[
 \appdef{X_t}{\mathcal C_L}{\mR}{x}{x(t)},
 \qquad 0\le t\le L.
\]
The coordinate maps then induce the raw and completed filtrations
\[
 \mathcal G_t^0:=\sigma(X_s:0\le s\le t),\qquad
 \mathcal G_t^\mu:=\mathcal G_t^0\vee\mathcal N^\mu,
\]
where
\[
 \mathcal N^\mu
 :=\left\{N\subset \mathcal C_L:\exists N_0\in\mathcal G_L^0,\ N\subset N_0,\ \mu(N_0)=0\right\}.
\]
We assume that $(\mathcal G_t^\mu)_{0\le t\le L}$ is right-continuous.
For fBm this follows from its causal Volterra representation and inverse
\cite{decreusefond1999,pipiras2002}. On $\mathcal D_L$, the coordinate
evaluations at times in $Q_L$ generate $\mathcal B(\mathcal D_L)$
\cite[Chapter~3]{billingsley1999}. Let $\Pi(Q_L)$ be the countable set of finite
partitions $0=q_0<\cdots<q_n=L$ with $q_j\in Q_L$. Right-continuity gives
\[
 \TV_{(0,L]}(z)
 =\sup_{\pi\in\Pi(Q_L)}\sum_{j=1}^{n(\pi)}|z_{q_j}-z_{q_{j-1}}|.
\]
Because the supremum is countable, the variation is Borel. Consequently,
\[
 \{z\in\mathcal D_L:\TV_{(0,L]}(z)<\infty\}
 =\bigcup_{m\in\Nn}\{z:\TV_{(0,L]}(z)\le m\}
 \in\mathcal B(\mathcal D_L).
\]
\begin{theorem}[path-space representation and law transport]
\label{thm:canonical-transport}
Let $X$ be the canonical coordinate just defined and let
$\psi=(\psi_t)_{0\le t\le L}$ be a process on $\mathcal C_L$ that is
adapted to $(\mathcal G_t^\mu)$, right-continuous, of finite variation
$\mu$-a.s., and integrable for the relevant cost.  There is a Borel map
$\Phi:\mathcal C_L\to \mathcal D_L$ such that
\[
  \mu\bigl(\Phi(X)_t=\psi_t\text{ for all }t\in[0,L]\bigr)=1.
\]
There is a $\mu$-null set outside which the map is causal at every
$q\in Q_L$: for each $q\in Q_L$, there is a Borel map
$f_q:\mathcal C([0,q];\mR)\to\mR$ satisfying
\begin{equation}
 \Phi(X)_q=f_q(X|_{[0,q]})
 \qquad\mu\text{-a.s.}
 \label{eq:canonical-nonanticipation}
\end{equation}
If $Y$ is any continuous process with law $\mu$, defined on any complete
probability space and equipped with its usual natural filtration, then
$\Phi(Y)$ is adapted and
\begin{equation}
 (Y,\Phi(Y))\overset d=(X,\Phi(X)).
 \label{eq:law-transport}
\end{equation}
Conversely, every adapted right-continuous finite-variation control on the
natural realization of $Y$ has such a Borel representative on
$\mathcal C_L$.
\end{theorem}
\begin{proof}
Rational coordinate measurability makes $\psi$ a $\mathcal D_L$-valued
variable measurable for the completion of $\sigma(X)$.
Standard Borel factorization \cite{kallenberg2021} gives a Borel map
$\Phi:\mathcal C_L\to\mathcal D_L$ with
\[
 \psi=\Phi(X)\quad\mu\text{-a.s. as path-valued variables}.
\]
For $q\in Q_L$, adaptedness and Doob--Dynkin give a Borel map
$f_q:\mathcal C([0,q];\mR)\to\mR$ such that
\[
 \Phi(X)_q=f_q(X|_{[0,q]})\quad\mu\text{-a.s.}
\]
We can combine the exceptional sets into one null set because $Q_L$ is
countable. Equality in law transfers all these identities to $Y$.
For $t<L$, choose $q_n\underset{n\to\infty}{\downarrow} t$ in $Q_L$. Then
\[
 \Phi(Y)_t=\lim_{n\to\infty}\Phi(Y)_{q_n},\qquad
 \bigcap_{n=1}^{\infty}\mathcal F_{q_n}^Y=\mathcal F_t^Y.
\]
The limit is measurable for each $\mathcal F_{q_n}^Y$, hence for
$\mathcal F_t^Y$. Thus $\Phi(Y)$ is adapted. The terminal time is covered because
$L\in Q_L$. The map $x\mapsto(x,\Phi(x))$ from $\mathcal C_L$ to
$\mathcal C_L\times\mathcal D_L$ is Borel. For any Borel
$A\subset\mathcal C_L\times\mathcal D_L$,
\[
 \Pp((Y,\Phi(Y))\in A)
 =\int_{\mathcal C_L}\ind_A(x,\Phi(x))\mu(\dd x)
 =\mu((X,\Phi(X))\in A).
\]
This proves \eqref{eq:law-transport} and preserves admissibility and
costs. Conversely, on a natural realization of $Y$, factor any admissible
$\mathcal D_L$-valued control through $Y$ and repeat the rational-time
argument on the canonical space.
\end{proof}
\begin{corollary}[equality of natural-filtration values]
\label{cor:natural-value-transport}
Let $Y$ be a continuous process with law $\mu$, equipped with its usual
natural filtration.  Let $J:\mathcal C_L\times \mathcal D_L\to[0,\infty]$ be Borel.
Then the infimum of $\E[J(Y,\psi)]$ over adapted right-continuous
finite-variation controls equals the corresponding infimum on the
path space.
\end{corollary}
\begin{proof}
By \cref{thm:canonical-transport}, every canonical rule transports to an
admissible rule for $Y$, and every admissible rule for $Y$ has a canonical
representative. In either direction,
\[
 \begin{aligned}
 \E[J(Y,\Phi(Y))]
 =
 \int_{\mathcal C_L}J(x,\Phi(x))\mu(\dd x)=\E_{\mu}[J(X,\Phi(X))].
 \end{aligned}
\]
Taking infima gives equality.
\end{proof}
\begin{lemma}[path-space concatenation with reset]
\label{lem:canonical-concatenation}
For $S,T>0$, let
$\Phi^S:\mathcal C([0,S];\mR)\to \mathcal D([0,S];\mR)$ and
$\Phi^T:\mathcal C([0,T];\mR)\to \mathcal D([0,T];\mR)$ be Borel representatives given by
\cref{thm:canonical-transport}.  For $x\in\mathcal C_{S+T}$, write
$x^S:=x|_{[0,S]}$ and $y_t:=x_{S+t}-x_S$. Define
\[
 \appdef{\Psi}{\mathcal C_{S+T}}{\mathcal D_{S+T}}{x}{
 \left(t\mapsto
 \begin{cases}
 \Phi^S(x^S)_t,&0\le t<S,\\
 x_S+\Phi^T(y)_{t-S},&S\le t\le S+T
 \end{cases}\right)}.
\]
The resulting $\Psi$ is Borel and adapted under fBm law. Its running
cost splits across the blocks, while its variation satisfies
\begin{equation}
 \CTV_{[0,S+T]}(\Psi)
 \le\CTV_{[0,S]}(\Phi^S(x^S))
 +|\Phi^S(x^S)_S-x_S|
 +\CTV_{[0,T]}(\Phi^T(y)),
 \label{eq:concatenation-variation}
\end{equation}
Moreover, the second rule alone determines the terminal error:
\begin{equation}
 |\Psi_{S+T}-x_{S+T}|
 =|\Phi^T(y)_T-y_T|.
 \label{eq:concatenation-terminal}
\end{equation}
\end{lemma}
\begin{proof}
Restriction and shifted-increment maps are continuous on the input
spaces. Every rational evaluation of $\Psi$ is Borel, and the
concatenation is c\`adl\`ag, so $\Psi$ is a Borel $\mathcal D_{S+T}$-valued
map. Under fBm law, $Y_r:=X_{S+r}-X_S$ has law $\mu_H^T$ and
\[
 \sigma(Y_v:0\le v\le r)\subset\mathcal F_{S+r}^{B^H}.
\]
Law transport therefore makes the second rule adapted.
At the interface, the triangle inequality gives
\begin{align*}
 |x_S+\Phi^T(y)_0-\Phi^S(x^S)_{S-}|
 \le |\Phi^S(x^S)_S-\Phi^S(x^S)_{S-}|
       +|\Phi^S(x^S)_S-x_S|+|\Phi^T(y)_0|.
\end{align*}
The first term belongs to the first variation, the second is its reset,
and the third is the second initial atom. All other increments lie in
disjoint blocks, proving \eqref{eq:concatenation-variation}.
To split the running cost, observe that
\[
 \Psi_t-x_t=\Phi^S(x^S)_t-x_t \qquad(0\le t<S),
\]
and
\[
 \Psi_{S+r}-x_{S+r}=\Phi^T(y)_r-y_r \qquad(0\le r\le T).
\]
These identities split the running integral. Taking $r=T$ also gives
\eqref{eq:concatenation-terminal}.
\end{proof}
The same factorization and concatenation arguments apply to full-past
rules on $\mathcal C_0(\mR)$, using rational coordinates on the whole
line. The map $x\mapsto(x_{S+r}-x_S)_{r\in\mR}$ preserves
the two-sided fBm law. They also apply to finite-window rules on
$[-M,L]$.
\section{Finite-memory approximation}
\label{app:finite-memory}
On the two-sided canonical fBm space, compare the finite-window and
full-past observation fields, defined for $0\le t\le L$ by
\[
 \mathcal H_t^{(M)}
 =\overline{\sigma(\mathbf B_s^H:-M\le s\le t)}^{\Pp},\qquad
 \mathcal H_t^{(-\infty)}
 =\overline{\sigma(\mathbf B_s^H:s\le t)}^{\Pp}.
\]
The full-past field is $\mathcal F_t^{\mathbf B}$ by
\cref{lem:mvn-filtration}. These finite-window fields are right-continuous after
completion for $t\ge0$. To see this, set
$Y_r:=\mathbf B_{-M+r}^H-\mathbf B_{-M}^H$ for $0\le r\le M+t$.
Since $\mathbf B_0^H=0$, we have $Y_M=-\mathbf B_{-M}^H$, so
\[
 \mathbf B_s^H=Y_{s+M}-Y_M\quad(-M\le s\le t).
\]
The two observed paths determine each other. The process $Y$ has the
one-sided fBm law, whose completed natural filtration is right-continuous.
\begin{proof}[Proof of \cref{lem:finite-memory-density}]
Write $J_L$ for the full-past reset cost. Given $\eta>0$, choose
$\psi$ with $J_L(\psi)\le u_H^{\past}(L)+\eta$. We approximate this rule by restricting its values to $[-R,R]$, discretizing time and then
restricting the observation window.\\
\emph{Restricting the control values.}
First set $\psi^R:=(-R)\vee\psi\wedge R$. Since the map
$x\mapsto(-R)\vee x\wedge R$ is $1$-Lipschitz and fixes zero,
\[
 \CTV_{[0,L]}(\psi^R)\le\CTV_{[0,L]}(\psi).
\]
Finite cost gives $\E\left[\int_0^L\psi_t^2\dd t\right]<\infty$ and
$\E\left[|\psi_L|\right]<\infty$. The bounds
\[
 |\psi_t^R-\mathbf B_t^H|^2\le2|\psi_t|^2+2|\mathbf B_t^H|^2,
\]
and
\[
 |\psi_L^R-\mathbf B_L^H|\le|\psi_L|+|\mathbf B_L^H|
\]
let us apply dominated convergence to the tracking and reset costs. Thus
$\limsup_RJ_L(\psi^R)\le J_L(\psi)$. Choose $R$ with increase at most
$\eta$ and relabel this bounded rule as $\psi$.\\
\emph{Time discretization.}
Next, on a dyadic partition $0=t_0<\cdots<t_K=L$, define
\[
 \psi_t^\pi:=\psi_{t_k}\quad(t_k\le t<t_{k+1}),\qquad
 \psi_L^\pi:=\psi_L.
\]
Ordered sampling gives
\[
 \CTV_{[0,L]}(\psi^\pi)
 =|\psi_0|+\sum_{k=1}^K|\psi_{t_k}-\psi_{t_{k-1}}|
 \le\CTV_{[0,L]}(\psi).
\]
Outside the countable pathwise jump set, $\psi_t^\pi\to\psi_t$ as $\operatorname{mesh}(\pi)\to0$. Boundedness and dominated convergence give convergence of the tracking
cost, while the reset is unchanged. Choose a grid that increases the cost by at
most $\eta$, and write $\xi_k:=\psi_{t_k}$.\\
\emph{Finite observation window.}
Finally, replace each coordinate by
$\xi_k^{(M)}:=\E[\xi_k\mid\mathcal H_{t_k}^{(M)}]$ and measure the
resulting error by $d_M:=\max_{0\le k\le K}\|\xi_k^{(M)}-\xi_k\|_{L^2}$.
Since, after completion,
\[
 \bigvee_{m\in\Nn}\mathcal H_{t_k}^{(m)}
 =\mathcal H_{t_k}^{(-\infty)},
\]
martingale convergence gives $d_M\xrightarrow[M\to\infty]{}0$. The coordinates remain bounded
by $R$. Their step interpolation $\psi^{\pi,M}$, including terminal
coordinate $\xi_K^{(M)}$, is a Borel finite-window adapted rule.
Since
$\E[\int_0^L(\psi_t^{\pi,M}-\psi_t^\pi)^2\dd t]\le Ld_M^2$,
expansion of the square and Cauchy--Schwarz give
\begin{align*}
 \left|\E\left[\int_0^L
 \bigl[(\psi_t^{\pi,M}-\mathbf B_t^H)^2
       -(\psi_t^\pi-\mathbf B_t^H)^2\bigr]\dd t\right]\right|\le2\sqrt L\,d_M
       \left(\E\left[\int_0^L(\psi_t^\pi-\mathbf B_t^H)^2\dd t\right]\right)^{1/2}
       +Ld_M^2\xrightarrow[M\to\infty]{}0.
\end{align*}
For each mesh jump,
\[
 \bigl||\xi_k^{(M)}-\xi_{k-1}^{(M)}|-|\xi_k-\xi_{k-1}|\bigr|
 \le|\xi_k^{(M)}-\xi_k|+|\xi_{k-1}^{(M)}-\xi_{k-1}|.
\]
Including the initial atom and terminal reset, the expected absolute
change in these costs is at most $(2K+2)d_M\xrightarrow[M\to\infty]{}0$.
For large $M$, the resulting rule has cost at most
$u_H^{\past}(L)+4\eta$. The finite-memory values decrease in $M$ and
are bounded below by $u_H^{\past}(L)$. Let $\eta\to0$.
Removing the terminal penalty gives the same proof for $v_H^{\past}$.
\end{proof}
\section{A continuous-time subadditive lemma}
\label{app:fekete}
Let $f:(0,\infty)\to[0,\infty)$ be subadditive and locally bounded.
Given $q>0$, write $L=nq+r$, where $n:=\lfloor L/q\rfloor\ge1$ and
$0\le r<q$. Local boundedness gives
$C_q:=\sup_{q\le s\le2q}f(s)<\infty$. Instead of using the remainder
$r$ alone, decompose $L=(n-1)q+(q+r)$ so that subadditivity only calls
on values away from zero:
\[
 f(L)\le(n-1)f(q)+f(q+r)\le nf(q)+C_q.
\]
As $n/L\xrightarrow[L\to\infty]{}1/q$,
\[
 \limsup_{L\to\infty}\frac{f(L)}L\le\frac{f(q)}q.
\]
Taking the infimum over $q>0$ gives
\[
 \begin{aligned}
 \inf_{q>0}\frac{f(q)}q
 \le\liminf_{L\to\infty}\frac{f(L)}L
 \le\limsup_{L\to\infty}\frac{f(L)}L
 \le\inf_{q>0}\frac{f(q)}q,
 \end{aligned}
\]
which proves the limit and infimum identity.
\section{The Brownian benchmark}
\label{app:brownian-verification}
\begin{proof}[Proof of \cref{prop:brownian}]
\emph{Verification lower bound.}
Set $d:=(3/2)^{1/3}$ and $\kappa:=d^2/2$. Define
\[
 \appdef{v}{\mR}{\mR}{x}{
 \begin{cases}
 \displaystyle\frac{d^2x^2}{2}-\frac{x^4}{12},&|x|\le d,\\[4pt]
 \displaystyle\frac{5d^4}{12}+|x|-d,&|x|>d.
 \end{cases}}
\]
Its derivatives satisfy
\begin{equation*}
 v'(x)=
 \begin{cases}
 d^2x-x^3/3,&|x|\le d,\\
 \sign(x),&|x|>d,
 \end{cases}
\end{equation*}
and
\begin{equation*}
 v''(x)=
 \begin{cases}
 d^2-x^2,&|x|\le d,\\
 0,&|x|>d.
 \end{cases}
\end{equation*}
Because $2d^3/3=1$, both derivatives match at $\pm d$.
Thus $v\in\mathcal C^2(\mR)$, $|v'|\le1$, $0\le v(x)\le|x|$, and
\[
 \frac12v''(x)+\frac12x^2-\kappa
 =\begin{cases}0,&|x|\le d,\\(x^2-d^2)/2,&|x|>d\end{cases}\ge0.
\]
For an admissible $\psi$, put $e:=\psi-B$ and $e_{0-}:=0$.
It\^o's formula, including the initial jump, and the Lipschitz bound give
\begin{align*}
 v(e_L)
 &=
 v(e_0)-v(0)+\int_{(0,L]}v'(e_{t-})\dd\psi_t^c\\
 &\quad+\sum_{\substack{0<t\le L\\ \Delta e_t\ne0}}\{v(e_t)-v(e_{t-})\}
       -\int_0^Lv'(e_{t-})\dd B_t+\frac12\int_0^Lv''(e_{t-})\dd t\\
 &\ge
 -\CTV_{[0,L]}(\psi)-\int_0^Lv'(e_{t-})\dd B_t+\frac12\int_0^Lv''(e_{t-})\dd t.
\end{align*}
The stochastic integral has mean zero because $v'$ is bounded.
Taking expectations gives
\begin{align*}
 \E\left[\frac12\int_0^Le_t^2\dd t+\CTV_{[0,L]}(\psi)\right]\ge\E\left[\int_0^L\frac{e_t^2+v''(e_t)}2\dd t\right]
      -\E[v(e_L)] \ge
 \kappa L-\E[v(e_L)].
\end{align*}
After adding the reset,
\begin{align*}
 \E\left[\frac12\int_0^Le_t^2\dd t+\CTV_{[0,L]}(\psi)+|e_L|\right]
 \ge\kappa L+\E[|e_L|-v(e_L)]
 \ge\kappa L.
\end{align*}
This is the singular-control verification argument
\cite{karatzas1983,harrison1983,taksar1985,protter2005}.\\
\emph{Reflected upper bound.}
For any $d>0$, let $e=-B+K^--K^+$ be reflected on $[-d,d]$, starting
at zero. The continuous nondecreasing regulators $K^\pm$ start at zero,
with $\dd K^-$ supported on $\{e=-d\}$ and $\dd K^+$ on $\{e=d\}$.
The control $\psi:=K^--K^+$ has augmented variation
$\CTV_{[0,L]}(\psi)=K_L^-+K_L^+$.
It\^o's formula gives
\begin{align*}
 e_L^2&=-2\int_0^Le_t\dd B_t+L-2d(K_L^-+K_L^+).
\end{align*}
Then,
\begin{align*}
 \E[K_L^-+K_L^+]
 &=\frac{L-\E[e_L^2]}{2d}
 =\frac{L}{2d}+O(1).
\end{align*}
To compute the running cost, use
\[
 \appdef{f}{[-d,d]}{\mR}{x}{\frac{x^4}{6}-\frac{d^2x^2}{3}}.
\]
Its derivatives satisfy
\[
 f'(x)=\frac{2x(x^2-d^2)}3,\qquad f'(\pm d)=0,\qquad
 \frac12f''(x)=x^2-\frac{d^2}{3}.
\]
\begin{samepage}
The regulator terms vanish, so
\begin{align*}
 \E[f(e_L)]-f(0)
 &=\E\left[\int_0^L\left(e_t^2-\frac{d^2}{3}\right)\dd t\right].
\end{align*}
Then,
\begin{align*}
 \E\left[\int_0^Le_t^2\dd t\right]
 =\frac{d^2L}{3}+\E[f(e_L)]-f(0)
 =\frac{d^2L}{3}+O(1).
\end{align*}
\end{samepage}
As $|e_L|\le d$, the reset contributes $O(1)$. The cost rate is
\[
 \appdef{c}{(0,\infty)}{(0,\infty)}{d}{\frac{d^2}{6}+\frac1{2d}}.
\]
Since
\[
 c'(d)=\frac d3-\frac1{2d^2},\qquad
 c''(d)=\frac13+\frac1{d^3}>0,
\]
the unique minimum occurs at $d^3=3/2$, and
\[
 c(d_*)=\frac{d_*^2}{6}+\frac1{2d_*}
 =\frac{d_*^2}{6}+\frac{d_*^2}{3}
 =\frac{d_*^2}{2}=\kappa.
\]
\end{proof}
\section{Prediction states}
\label{app:prediction-states}
\subsection{The prediction curve as a state variable}
Let $W=(W_t)_{t\ge0}$ be standard Brownian motion with augmented natural
filtration $\mathbf F^W$. Choose a jointly Borel,
filtration-preserving Volterra kernel
$\mathcal K_H:\{(t,s)\in[0,\infty)^2:s\le t\}\to\mR$
\cite{decreusefond1999} such that
\begin{equation}
 B_t^H=\int_0^t\mathcal K_H(t,s)\dd W_s,
 \qquad \mathcal F_t^{B^H}=\mathcal F_t^W,
 \label{eq:volterra-fbm-representation}
\end{equation}
with $\mathcal K_H(t,\cdot)\in L^2([0,t])$. The state space
$\cE:=\mathcal C_{\loc}([0,\infty);\mR)$ is Polish under uniform
convergence on compact intervals.
Write $\mathcal N_{\cE}(g,Q)$ for the
Gaussian probability law on $\cE$ with mean $g\in\cE$ and covariance
kernel $Q$. For $h\ge0$, the shift
\[
 \appdef{\tau_h}{\cE}{\cE}{g}{\bigl(r\mapsto g(h+r)\bigr)}
\]
is continuous. The \emph{prediction curve} is defined by
\begin{equation}
 G_t(r):=\E[B_{t+r}^H\mid\mathcal F_t^{B^H}]
 =\int_0^t\mathcal K_H(t+r,s)\dd W_s,
 \qquad r\ge0.
 \label{eq:prediction-curve}
\end{equation}
The $L^2$ contraction of conditional expectation, followed by Gaussian
moment equivalence, gives
\begin{align*}
 \E\left[|G_t(r)-G_t(q)|^2\right]
 \le
 \E[|B_{t+r}^H-B_{t+q}^H|^2]
 \overset{\eqref{eq:fbm-covariance}}=|r-q|^{2H}.
\end{align*}
Then,
\begin{align*}
 \E\left[|G_t(r)-G_t(q)|^p\right]
 &\le
 C_p|r-q|^{pH},\qquad p\ge2.
\end{align*}
For $p>1/H$, Kolmogorov's criterion \cite[Chapter~4]{kallenberg2021}
gives a continuous version on each
compact $r$-interval. Start from the jointly measurable parameterized
Wiener integrals for $r\in\Qq_+$ and extend by continuity. Separability
of $\cE$ gives a version with $G_t(0)=B_t^H$ such that
\[
 \left|\begin{array}{rcl}
 [0,\infty)\times\Omega&\longrightarrow&\cE\\
 (t,\omega)&\longmapsto&G_t(\omega)
 \end{array}\right.
\]
is $\mathcal B([0,\infty))\otimes\mathcal F$--$\mathcal B(\cE)$
measurable.
\begin{theorem}[Markov property of the prediction state]
\label{thm:prediction-state}
The $[0,\infty)\times\cE$-valued process
$((t,G_t))_{t\ge0}$ is Markov. For $h,r\ge0$,
\begin{equation}
 G_{t+h}(r)=G_t(h+r)
 +\int_t^{t+h}\mathcal K_H(t+h+r,s)\dd W_s,
 \label{eq:prediction-update}
\end{equation}
and, conditionally on $\mathcal F_t^{B^H}$,
\begin{equation}
 B_{t+r}^H
 =G_t(r)+\int_t^{t+r}\mathcal K_H(t+r,s)\dd W_s.
 \label{eq:future-decomposition}
\end{equation}
With $Q_t$ defined in \eqref{eq:future-covariance},
\[
  \mathcal L\left((B_{t+r}^H)_{r\ge0}\mid\mathcal F_t^{B^H}\right)
  =\mathcal N_{\cE}(G_t,Q_t).
\]
\end{theorem}
\begin{proof}
Splitting the Wiener integral at $t$ gives
\begin{align*}
 B_{t+r}^H
 &=\int_0^t\mathcal K_H(t+r,s)\dd W_s
   +\int_t^{t+r}\mathcal K_H(t+r,s)\dd W_s\\
 &=G_t(r)+\int_t^{t+r}\mathcal K_H(t+r,s)\dd W_s.
\end{align*}
The second term is centered Gaussian and independent of
$\mathcal F_t^{B^H}=\mathcal F_t^W$. Its covariance is
\begin{equation}
 Q_t(r,q):=\int_t^{t+(r\wedge q)}
 \mathcal K_H(t+r,s)\mathcal K_H(t+q,s)\dd s.
 \label{eq:future-covariance}
\end{equation}
Hence
\[
 \mathcal L\bigl((B_{t+r}^H)_{r\ge0}\mid\mathcal F_t^{B^H}\bigr)
 =\mathcal N_{\cE}(G_t,Q_t).
\]
Likewise,
\begin{align*}
 G_{t+h}(r)
 &=\int_0^t\mathcal K_H(t+h+r,s)\dd W_s
   +\int_t^{t+h}\mathcal K_H(t+h+r,s)\dd W_s\\
 &=G_t(h+r)+\int_t^{t+h}\mathcal K_H(t+h+r,s)\dd W_s.
\end{align*}
Denote the last integral, as a function of $r$, by $I_{t,h}$.
The Gaussian increment estimate gives a continuous version
$I_{t,h}\in L^0(\Omega;\cE)$.
It is independent of $\mathcal F_t^{B^H}$, centered Gaussian, and has covariance
\[
 \E[I_{t,h}(r)I_{t,h}(q)]
 =\int_t^{t+h}\mathcal K_H(t+h+r,s)\mathcal K_H(t+h+q,s)\dd s.
\]
Define the transition map
\[
 \appdef{P_{t,h}}{\cE}{\mathcal P(\cE)}
 {g}{\mathcal L(\tau_hg+I_{t,h})}.
\]
For every $A\in\mathcal B(\cE)$,
\[
 P_{t,h}(g,A)=\E[\ind_A(\tau_hg+I_{t,h})].
\]
Since addition and $\tau_h$ are continuous, $P_{t,h}$ is a Borel kernel.
Independence gives
\[
 \mathcal L(G_{t+h}\mid\mathcal F_t^{B^H})=P_{t,h}(G_t,\cdot).
\]
Together with the deterministic update $t\mapsto t+h$, this proves the
Markov property of $(t,G_t)$.
\end{proof}
\begin{corollary}[prediction state for continuation values]
\label{cor:prediction-control-state}
Fix $t\ge0$, $\tau>0$, and $\lambda>0$, and write
$\E_t[\cdot]:=\E[\cdot\mid\mathcal F_t^{B^H}]$.  Let
$\bar a,\bar b:[0,\tau]\to(0,\infty)$ be bounded Borel functions bounded
away from zero.  For $x\in\mR$, let $\mathfrak A_{t,t+\tau}(x)$ be the set of
$\mathbf F^{B^H}$-adapted right-continuous finite-variation controls $\xi$ on
$[t,t+\tau]$ with pre-initial value $\xi_{t-}=x$ and finite expected cost.
The variation on this interval includes the atom $|\xi_t-x|$.  Define
\begin{equation}
 \begin{aligned}
 \mathcal V_{t,\tau}(x)
 :=\operatorname*{ess\,inf}_{\xi\in\mathfrak A_{t,t+\tau}(x)}
 \E_t\biggl[
 \frac12\int_t^{t+\tau}\bar a_{s-t}(\xi_s-B_s^H)^2\dd s
 +\lambda\int_{[t,t+\tau]}\bar b_{s-t}\dd|\xi|_s
 \biggr].
 \end{aligned}
 \label{eq:conditional-continuation-value}
\end{equation}
There is a Borel function
$\mathsf V_{t,\lambda,\tau}^{\bar a,\bar b}$ on
$\cE\times\mR$, with values in $[0,\infty)$,
such that, for every $x\in\mR$,
\begin{equation}
 \mathcal V_{t,\tau}(x)
 =\mathsf V_{t,\lambda,\tau}^{\bar a,\bar b}(G_t,x)
 \qquad\Pp\text{-a.s.}
 \label{eq:prediction-sufficient-state}
\end{equation}
Since $G_t(0)=B_t^H$, the same value can be written in the variables
$(G_t,e)$ with $e:=x-G_t(0)$ and $x=e+G_t(0)$.
\end{corollary}
\begin{proof}
Introduce the future innovation $I^t:=(B_{t+s}^H-G_t(s))_{0\le s\le\tau}$
and the observed past $P_t:=(B_r^H)_{0\le r\le t}$. These are path-valued
random variables in $\mathcal C([0,\tau];\mR)$ and
$\mathcal C([0,t];\mR)$, respectively.
By \cref{thm:prediction-state}, $I^t$ is independent of
$\mathcal F_t^{B^H}$ and has covariance $Q_t$. Let $\mathbf H^t$ be
the usual augmentation of
$\mathcal H_s^{t,0}:=\sigma(G_t,I_r^t:0\le r\le s)$.
Brownian independence and the zero--one law imply
\[
 \mathcal H_0^t=\overline{\sigma(G_t)}^{\Pp},\qquad
 \mathcal F_{t+s}^{B^H}
 =\overline{\mathcal H_s^t\vee\sigma(P_t)}^{\Pp},\qquad
 P_t\indep\mathcal H_\tau^t\mid G_t.
\]
To apply \eqref{eq:conditional-seed-initial-value}, center the target
and initial position at $G_t(0)$ by setting
$X_s:=G_t(s)-G_t(0)+I_s^t$ and $x_0:=x-G_t(0)$.
Then $X_0=0$ and $x_0\in L^0(\Omega,\mathcal H_0^t,\Pp)$.
The weighted continuation cost is unchanged. The conditional
initial-information identity, applied with the auxiliary variable $P_t$,
shows that $\mathcal V_{t,\tau}(x)$ is measurable with respect to $\sigma(G_t)$.
The Doob--Dynkin lemma \cite{kallenberg2021} gives a Borel function
$V_x$ on $\cE$, with values in
$[0,\infty)$, satisfying
$\mathcal V_{t,\tau}(x)=V_x(G_t)$ almost surely. 
For any $x,y\in\mR$, the same control is admissible for both positions,
and its initial atoms differ by at most $|x-y|$. Thus
\[
 |\mathcal V_{t,\tau}(x)-\mathcal V_{t,\tau}(y)|
 \le\lambda\bar b_0|x-y|\quad\text{a.s.}
\]
Let $E_0\subset\cE$ be the Borel full-$\mathcal L(G_t)$-measure set
on which all rational values $V_q$ are finite and obey this bound.
Define
\[
 \appdef{\mathsf V_{t,\lambda,\tau}^{\bar a,\bar b}}
 {\cE\times\mR}{[0,\infty)}{(g,x)}{
 \begin{cases}
 \displaystyle\inf_{q\in\Qq}\{V_q(g)+\lambda\bar b_0|x-q|\},&g\in E_0,\\
 0,&g\notin E_0.
 \end{cases}}
\]
This countable infimum is Borel and finite. For $g\in E_0$, it is the
unique $\lambda\bar b_0$-Lipschitz extension of the rational values.
Approximating any fixed $x$ by rationals proves
\eqref{eq:prediction-sufficient-state}.
\end{proof}
\begin{remark}
The state is generally infinite-dimensional. See also
\cite{carmona1998,viens2019,cuchiero2020,kirchner2025}.
\end{remark}
\subsection{Why the scalar error is not a general control state}
The next one-step problem shows that a quadratic tracking decision can use
information beyond the current error.
For \emph{fractional Gaussian noise} $\Delta_k:=B_k^H-B_{k-1}^H$, the
increment covariances are
\[
 \gamma_k:=\Cov(\Delta_1,\Delta_{k+1})
 =\tfrac12\bigl((k+1)^{2H}-2k^{2H}+|k-1|^{2H}\bigr).
\]
\begin{proposition}[a one-step test of scalar-state sufficiency]
\label{prop:memory}
Fix $a,b,\lambda>0$ and $u\in\mR$. For each observation
$(\Delta_1,\Delta_2)$, let $v^*$ be the unique minimizer over $v\in\mR$
of the conditional objective
\begin{equation}
 \frac a2\E\left[(v-B_3^H)^2\mid\Delta_1,\Delta_2\right]
 +\lambda b|v-u|.
 \label{eq:memory-subproblem}
\end{equation}
The minimizer satisfies
$v^*\in L^0(\Omega,\sigma(\Delta_1,\Delta_2),\Pp)$.
If $H\ne1/2$, $v^*$ is not $\sigma(B_2^H)$-measurable and therefore cannot
be a function of $u-B_2^H$ alone.
\end{proposition}
\begin{proof}
Let $\widehat B_3:=\E[B_3^H\mid\Delta_1,\Delta_2]$. The identity
\[
 \E[(v-B_3^H)^2\mid\Delta_1,\Delta_2]
 =(v-\widehat B_3)^2+\Var(B_3^H\mid\Delta_1,\Delta_2)
\]
reduces the problem to the strictly convex function
$a(v-\widehat B_3)^2/2+\lambda b|v-u|$.
Its subgradient condition is
\[
 0\in a(v-\widehat B_3)+\lambda b\,\partial|v-u|,
 \qquad
 \partial|v-u|=
 \begin{cases}\{-1\},&v<u,\\[-2pt][-1,1],&v=u,\\\{1\},&v>u.\end{cases}
\]
Solving on $v<u$, $v=u$, and $v>u$ gives
\begin{equation}
 v^*=u+\sign(\widehat B_3-u)
 \left(|\widehat B_3-u|-\frac{\lambda b}{a}\right)_+.
 \label{eq:memory-action}
\end{equation}
Put $r:=\gamma_1=2^{2H-1}-1\in(-1/2,1)$ and $s:=\gamma_2$.
Stationarity identifies the observed covariance matrix and its
cross-covariance with the next increment:
\[
 \Cov((\Delta_1,\Delta_2))=
 \begin{pmatrix}1&r\\r&1\end{pmatrix},\qquad
 \Cov(\Delta_3,(\Delta_1,\Delta_2))=(s,r).
\]
Gaussian regression therefore yields
\begin{align*}
 \E[\Delta_3\mid\Delta_1,\Delta_2]
 &=\frac1{1-r^2}(s,r)
 \begin{pmatrix}1&-r\\-r&1\end{pmatrix}
 \binom{\Delta_1}{\Delta_2}\\
 &=\frac{s-r^2}{1-r^2}\Delta_1
   +\frac{r(1-s)}{1-r^2}\Delta_2.
\end{align*}
Using $B_2^H=\Delta_1+\Delta_2$ and
\[
 \frac{s-r^2-r(1-s)}{1-r^2}
 =\frac{(s-r)(1+r)}{(1-r)(1+r)}
 =\frac{s-r}{1-r},
\]
we obtain
\begin{equation}
 \widehat B_3=
 \left(1+\frac{r(1-s)}{1-r^2}\right)B_2^H
 +\frac{s-r}{1-r}\Delta_1.
 \label{eq:prediction-two-increments}
\end{equation}
For $\alpha:=2H\in(0,2)$, define
\[
 \appdef{d}{[0,\infty)}{\mR}{x}{(x+1)^\alpha-x^\alpha}.
\]
For $k\ge1$ and $x>0$, its increments and curvature satisfy
\begin{align*}
 2\gamma_k=d(k)-d(k-1),\qquad
 d''(x)&=\alpha(\alpha-1)
 \bigl((x+1)^{\alpha-2}-x^{\alpha-2}\bigr).
\end{align*}
The bracket is negative. Hence $d$ is strictly convex for $\alpha<1$
and strictly concave for $\alpha>1$, so
\[
 \gamma_2>\gamma_1\quad(H<1/2),\qquad
 \gamma_2<\gamma_1\quad(H>1/2).
\]
Thus $s-r\ne0$ for $H\ne1/2$. Finally,
\begin{align*}
 \Var(B_2^H)&=2(1+r),\qquad
 \Cov(\Delta_1,B_2^H)=1+r.
\end{align*}
Then,
\begin{align*}
 \Var(\Delta_1\mid B_2^H)
 &=
 1-\frac{(1+r)^2}{2(1+r)}=\frac{1-r}{2}.
\end{align*}
Then,
\begin{align*}
 \Var(\widehat B_3\mid B_2^H)
 \overset{\eqref{eq:prediction-two-increments}}=
 \left(\frac{s-r}{1-r}\right)^2\frac{1-r}{2}
 =\frac{(s-r)^2}{2(1-r)}>0.
\end{align*}
The conditional Gaussian law of $\widehat B_3$ therefore has full support.
The soft-threshold map in \eqref{eq:memory-action} has slope one on
both tails, each of positive conditional probability. Thus
$\Var(v^*\mid B_2^H)>0$, excluding $\sigma(B_2^H)$-measurability and,
for fixed $u$, any representation through $u-B_2^H$ alone.
\end{proof}
\section{Technical estimates for local homogenization}
\label{app:local-estimates}
\subsection{Stability of reset values}
\label{app:reset-stability}
\begin{proof}[Proof of \cref{lem:value-stability}]
The same control $\psi\in\cA_L(\mathbf G)$ is admissible for both targets.
The variation terms cancel, so we can expand the squares and apply the
triangle inequality. Cauchy--Schwarz then controls the mixed term.
\begin{align*}
 J_L^X(\psi)-J_L^Y(\psi)
 &=\E\left[\int_0^L\left\{
       -(\psi_t-Y_t)(X_t-Y_t)+\frac12|X_t-Y_t|^2
       \right\}\dd t\right]\\
 &\qquad{}+\E[|\psi_L-X_L|-|\psi_L-Y_L|]\\
 &\quad\le
   \E\left[\int_0^L|\psi_t-Y_t|\,|X_t-Y_t|\dd t\right]
   +\frac{A_L}{2}+B_L\\
 &\quad\le
   \left(\E\left[\int_0^L(\psi_t-Y_t)^2\dd t\right]\right)^{1/2}
   \sqrt{A_L}+\frac{A_L}{2}+B_L\\
 &\quad\le
   \sqrt{2J_L^Y(\psi)A_L}+\frac{A_L}{2}+B_L.
\end{align*}
Taking a minimizing sequence for $Y$ gives
\[
 U_X-U_Y\le\sqrt{2U_YA_L}+A_L/2+B_L.
\]
Interchange $X,Y$ and apply the bound from the smaller value to the
larger to obtain \eqref{eq:value-stability-quantitative}.
Moreover, Young's inequality gives
\[
 \sqrt{2U_YA_L}\le
 \delta U_Y+\frac{A_L}{2\delta}
\]
and hence \eqref{eq:value-stability}. Finally, if $U_X,U_Y\le CL$,
divide by $L$ to get
\[
 \frac{|U_X-U_Y|}{L}
 \le\frac{B_L}{L}+\frac{A_L}{2L}+\sqrt{2C\frac{A_L}{L}}\xrightarrow[L\to\infty]{}0.
\]
\end{proof}
\subsection{Uniform growing-cell approximation}
\label{app:growing-cells}
\begin{proof}[Proof of \cref{lem:growing-cell-transfer}]
We first choose a horizon on which the tangent error remains uniform.
Put $d_J:=T-\sup J>0$. For integers $m\ge2$ and
$0<h<\min\{h_*,d_J/(c_+m^2)\}$, define
\[
 \appdef{\varepsilon_m}{(0,\min\{h_*,d_J/(c_+m^2)\})}{[0,\infty)}{h}
 {\displaystyle\sup_{\substack{t\in J,\ c\in[c_-,c_+]\\0\le u\le m^2}}
 \|X_u^{t,h,c}-\mathbf B_u^{t,h,c}\|_{L^2(\overline\Pp)}}.
\]
By \textnormal{(FT2)}, $\varepsilon_m(h)\xrightarrow[h\to0]{}0$ for each fixed $m$.
Starting with any $h_1>0$, choose recursively
\[
 0<h_m<\min\left\{\frac{h_{m-1}}2,\frac{h_*}2,\frac1{m^2},
                    \frac{d_J}{2c_+m^2}\right\},
 \qquad
 \sup_{0<h\le h_m}\varepsilon_m(h)\le m^{-2}.
\]
Set $L_h:=m$ for $h\in(h_{m+1},h_m]$ and assign it any positive value
for larger $h$. Then
\[
 L_h\underset{h\to0}{\uparrow}\infty,
\quad\text{and}\quad
 hL_h\le L_h^{-1}\xrightarrow[h\to0]{}0.
\]
Moreover,
\[
 \sup_{\substack{t\in J,\ c\in[c_-,c_+]}}
 \sup_{0\le u\le L_h^2}
 \|X_u^{t,h,c}-\mathbf B_u^{t,h,c}\|_{L^2}
 \le L_h^{-2}.
\]
In particular, this controls $[0,2L_h]$ for large $L_h$. Tonelli's
theorem and Cauchy--Schwarz therefore give
\begin{equation*}
 \E\left[\int_0^{L_h}|X_u^{t,h,c}-\mathbf B_u^{t,h,c}|^2\dd u\right]
 \le L_h\varepsilon_{L_h}(h)^2\le L_h^{-3},
\end{equation*}
and
\begin{equation*}
 \E\left[|X_{L_h}^{t,h,c}-\mathbf B_{L_h}^{t,h,c}|\right]
 \le \varepsilon_{L_h}(h)\le L_h^{-2}.
\end{equation*}
\textnormal{(FT1)} and \cref{lem:independent-randomization} give
\begin{equation*}
 \mathcal U_{L_h}(X^{t,h,c};(\mathcal F_{t+chu})_{u\le L_h})
 =\mathcal U_{L_h}(X^{t,h,c};\overline{\mathbf G}^{t,h,c}),
\end{equation*}
and
\begin{equation*}
 \mathcal U_{L_h}(\mathbf B^{t,h,c};\overline{\mathbf G}^{t,h,c})
 =u_H^{\past}(L_h).
\end{equation*}
For $0\le v\le1$ and $j=0,\ldots,L_h-1$, \textnormal{(FT3)} and
strict positivity of $\rho$ give
\[
\begin{aligned}
 \E\left[|X_{j+v}^{t,h,c}-X_j^{t,h,c}|^2\right]
 =\frac{\E\left[|\Theta_{t+ch(j+v)}-\Theta_{t+chj}|^2\right]}
        {\rho(t)^2(ch)^{2H}}
 \le
 \frac{C(chv)^{2H}}{\rho(t)^2(ch)^{2H}}
 \le C_Jv^{2H}.
\end{aligned}
\]
A rule with time step one therefore bounds both reset values by $C(1+L_h)$,
uniformly in $t,c,h$. The stability estimate yields
\[
 \frac{|\mathcal U_{L_h}(X^{t,h,c})-u_H^{\past}(L_h)|}{L_h}
 \overset{\eqref{eq:value-stability-quantitative}}\le
 L_h^{-3}+\tfrac12L_h^{-4}+CL_h^{-2}
 \xrightarrow[h\to0]{}0
\]
uniformly. Now use $u_H^{\past}(L_h)/L_h\xrightarrow[h\to0]{}\kappa_H$. 
Let $R_h:=qL_h+r_h$, with $q$ in a compact subset of $(0,\infty)$ and
$r_h=O(1)$ uniformly in all parameters. For small $h$, $0<R_h\le L_h^2$ and $R_h\asymp L_h$.
The same coupling gives, uniformly,
\begin{equation*}
 \E\left[\int_0^{R_h}|X_u^{t,h,c}-\mathbf B_u^{t,h,c}|^2\dd u\right]
 \le R_hL_h^{-4}=o(R_h),
\end{equation*}
and
\begin{equation*}
 \E\left[|X_{R_h}^{t,h,c}-\mathbf B_{R_h}^{t,h,c}|\right]
 \le L_h^{-2}=o(R_h).
\end{equation*}
Fix $0<q_-<q_+<\infty$. The values obtained with time step one are
$O(R_h)$ uniformly in $q\in[q_-,q_+]$, while
\[
 \sup_{q\in[q_-,q_+]}
 \left|
 \frac{u_H^{\past}(qL_h+r_h)}{qL_h+r_h}-\kappa_H
 \right|
 \xrightarrow[h\to0]{}0,
\]
since the horizons diverge uniformly. Stability then gives the transfer
statement with $L_h$ replaced by $R_h$.
\end{proof}
\subsection{Block estimates}
\label{app:block-estimates}
\begin{lemma}[block reset decomposition]
\label{lem:virtual-reset-decomposition}
Let $\theta\in\mathcal C([0,T];\mR)$ and let $\phi$ be right-continuous
with finite variation. For $0\le r_0<r_1<\cdots<r_N\le T$, the reset variations on the blocks satisfy
\begin{equation}
\begin{aligned}
 \sum_{i=1}^N\Bigl[
 |\phi_{r_{i-1}}-\theta_{r_{i-1}}|
 +\TV_{(r_{i-1},r_i]}(\phi)
 +|\phi_{r_i}-\theta_{r_i}|\Bigr]
 \le\TV_{(r_0,r_N]}(\phi)
 +2\sum_{j=0}^N|\phi_{r_j}-\theta_{r_j}|.
\end{aligned}
 \label{eq:virtual-reset-decomposition}
\end{equation}
If $b:[0,T]\to(0,\infty)$ is Borel and $b\le b_+$, then the weighted
reset variations satisfy
\begin{align}
 &\sum_{i=1}^N\Bigg[
 b_{r_{i-1}}|\phi_{r_{i-1}}-\theta_{r_{i-1}}|
 +\int_{(r_{i-1},r_i]}b_t\dd|\phi|_t
 +b_{r_i}|\phi_{r_i}-\theta_{r_i}|\Bigg]\notag\\
 &\qquad\le\int_{(r_0,r_N]}b_t\dd|\phi|_t
 +2b_+\sum_{j=0}^N|\phi_{r_j}-\theta_{r_j}|.
 \label{eq:weighted-virtual-reset-decomposition}
\end{align}
\end{lemma}
\begin{proof}
Additivity on the disjoint intervals $(r_{i-1},r_i]$ makes the
left-hand side of \eqref{eq:virtual-reset-decomposition} equal to
\begin{align*}
 &\TV_{(r_0,r_N]}(\phi)
 +|\phi_{r_0}-\theta_{r_0}|
 +2\sum_{j=1}^{N-1}|\phi_{r_j}-\theta_{r_j}|+|\phi_{r_N}-\theta_{r_N}|\\
 &\qquad\le\TV_{(r_0,r_N]}(\phi)+2\sum_{j=0}^N|\phi_{r_j}-\theta_{r_j}|.
\end{align*}
Likewise, additivity gives the first expression below for the weighted
left-hand side. The bound then follows from $b\le b_+$:
\begin{align*}
 &\int_{(r_0,r_N]}b_t\dd|\phi|_t
 +b_{r_0}|\phi_{r_0}-\theta_{r_0}|
 +2\sum_{j=1}^{N-1}b_{r_j}|\phi_{r_j}-\theta_{r_j}|
 +b_{r_N}|\phi_{r_N}-\theta_{r_N}|\\
 &\qquad\le\int_{(r_0,r_N]}b_t\dd|\phi|_t
 +2b_+\sum_{j=0}^N|\phi_{r_j}-\theta_{r_j}|.
\end{align*}
\end{proof}
\begin{lemma}[separator points]
\label{lem:separator-points}
Let $e$ be a measurable process with
$\E[\int_0^T e_t^2\dd t]<\infty$.  Let
$J_1,\ldots,J_N$ be disjoint intervals of the same length
$h>0$.  There are deterministic points $r_i\in J_i$ such that
\begin{equation}
 \E[|e_{r_i}|^2]
 \le
 \frac1h\int_{J_i}\E[|e_t|^2]\dd t.
 \label{eq:separator-pointwise}
\end{equation}
Then,
\begin{equation}
 \sum_{i=1}^N\E[|e_{r_i}|]
 \le
 \left(
 \frac Nh
 \E[\int_{\cup_iJ_i}e_t^2\dd t]
 \right)^{1/2}.
 \label{eq:separator-sum}
\end{equation}
\end{lemma}
\begin{proof}
Choose $r_i\in J_i$ with
$\E\left[|e_{r_i}|^2\right]\le h^{-1}\int_{J_i}\E\left[|e_t|^2\right]\dd t$. Such a point
exists by the integral average. Two applications of Cauchy--Schwarz give
\begin{align*}
 \sum_{i=1}^N\E\left[|e_{r_i}|\right]
 &\le
 \sum_{i=1}^N(\E\left[|e_{r_i}|^2\right])^{1/2}\\
 &\le
 \left(N\sum_{i=1}^N\E\left[|e_{r_i}|^2\right]\right)^{1/2}\\
 &\overset{\eqref{eq:separator-pointwise}}\le
 \left(\frac Nh\sum_{i=1}^N\int_{J_i}\E\left[|e_t|^2\right]\dd t\right)^{1/2}\\
 &=
 \left(\frac Nh\E\left[\int_{\cup_iJ_i}e_t^2\dd t\right]\right)^{1/2}.
\end{align*}
\end{proof}
\begin{lemma}[short-interval control]
\label{lem:short-interval-rule}
Assume \eqref{eq:strong-tangent-global-bound}, and let $a,b$ be bounded
above.  Put $h_\lambda:=\lambda^{1/(1+H)}$ and
$\alpha_H:=2H/(1+H)$.  On every interval $I=[s,s+r]\subset[0,T]$, there
is a control that starts from $\Theta_s$, resets to
$\Theta_{s+r}$ at the end, and has expected cost at most
\begin{equation}
 C\lambda^{\alpha_H}(r+h_\lambda),
 \label{eq:short-interval-cost}
\end{equation}
where $C$ is independent of $s,r,\lambda$, for all sufficiently small
$\lambda>0$.
\end{lemma}
\begin{proof}
Take $h_\lambda<h_0$. For $r>0$, set
$N:=\lceil r/h_\lambda\rceil$ and
$t_k:=s+\min(kh_\lambda,r)$, $0\le k\le N$.
Use $\phi_t:=\Theta_{t_k}$ on $[t_k,t_{k+1})$, with
$\phi_{s+r}:=\Theta_{s+r}$. The initial adjustment is zero, and the
increment bound controls the tracking cost. For the variation cost,
Cauchy--Schwarz and $N\le r/h_\lambda+1$ give the remaining bounds:
\begin{align*}
 \E\left[\int_I|\phi_t-\Theta_t|^2\dd t\right]
 \overset{\eqref{eq:strong-tangent-global-bound}}\le
 C\sum_{k=0}^{N-1}\int_0^{t_{k+1}-t_k}u^{2H}\dd u
 \le Crh_\lambda^{2H},
\end{align*}
and
\begin{align*}
 \E\left[\CTV_{[s,s+r]}(\phi)\right]
 =\sum_{k=0}^{N-1}\E\left[|\Theta_{t_{k+1}}-\Theta_{t_k}|\right]
 \le CNh_\lambda^H
 \le
 C(rh_\lambda^{H-1}+h_\lambda^H).
\end{align*}
Since $\lambda=h_\lambda^{1+H}$ and
$\lambda^{\alpha_H}=h_\lambda^{2H}$, the weighted cost is at most
\[
 C\{rh_\lambda^{2H}+\lambda rh_\lambda^{H-1}+\lambda h_\lambda^H\}
 \le C\lambda^{\alpha_H}(r+h_\lambda).
\]
The case $r=0$ is immediate.
\end{proof}
\subsection{The Gaussian Volterra tangent}
\label{app:volterra-verification}
\begin{proof}[Proof of \cref{prop:volterra-strong-tangent}]
Fix a compact interval $J\subset(0,T)$, $0<c_-<c_+<\infty$, and
$R<\infty$. The limits below are uniform in $t\in J$,
$c\in[c_-,c_+]$, and $u\in[0,R]$.\\
We first construct the tangent and its observation filtration.
On an extension $(\overline\Omega,\overline{\mathcal F},\overline\Pp)$,
extend $W$ to a two-sided Brownian motion $\widetilde W$ by an independent
past and, if needed, an independent continuation after $T$. Set $\beta:=H-1/2\in(-1/2,1/2)$ and
\[
 \appdef{\widetilde W^{t,h}}{\mR}{L^2(\overline\Omega,\overline\Pp)}{z}
 {h^{-1/2}(\widetilde W_{t+hz}-\widetilde W_t)}.
\]
This is two-sided standard Brownian motion. Use the kernel $\mathfrak g_u^H$ from \eqref{eq:mvn-kernel}.
For general $H$, $\mathfrak g_u^H\in L^2(\mR)$: its local squared
singularities have exponent $2\beta=2H-1>-1$, while
\[
 \mathfrak g_u^H(z)=O_u(|z|^{H-3/2})\quad(z\to-\infty),
 \qquad 2H-3<-1.
\]
\emph{Wiener integration} is the linear isometry from $L^2(\mR)$ into
$L^2(\overline\Omega,\overline\Pp)$ satisfying
\[
 \E\left[\left(\int_{\mR}f\dd\widetilde W\right)\left(\int_{\mR}g\dd\widetilde W\right)\right]
 =\int_\mR f(z)g(z)\dd z,\qquad f,g\in L^2(\mR).
\]
The kernel identities
\begin{equation*}
 \mathfrak g_{au}^H(az)
 =a^{H-1/2}\mathfrak g_u^H(z)\qquad(a>0),
\end{equation*}
and
\begin{equation*}
 \mathfrak g_u^H(z)-\mathfrak g_v^H(z)
 =\mathfrak g_{u-v}^H(z-v)
\end{equation*}
give, by scaling, translation, and $\mathfrak g_{-u}^H(z)=-\mathfrak g_u^H(z+u)$,
\begin{align*}
 \|\mathfrak g_1^H\|_2^2
 =\int_0^1(1-z)^{2H-1}\dd z+\int_0^\infty\bigl((1+z)^{H-1/2}-z^{H-1/2}\bigr)^2\dd z
 =\Gamma_H.
\end{align*}
Then,
\[
 \|\mathfrak g_u^H\|_2^2=\Gamma_H|u|^{2H},
\]
and
\[
 \|\mathfrak g_u^H-\mathfrak g_v^H\|_2^2
 =\Gamma_H|u-v|^{2H}.
\]
Polarization therefore identifies the covariance:
\begin{align*}
 \langle\mathfrak g_u^H,\mathfrak g_v^H\rangle_{L^2(\mR)}=
 \frac{\Gamma_H}{2}(|u|^{2H}+|v|^{2H}-|u-v|^{2H})=\Gamma_HR_H(u,v).
\end{align*}
Consequently,
\begin{equation}
 \appdef{\mathbf B^{t,h}}{\mR}{L^2(\overline\Omega,\overline\Pp)}{u}
 {\displaystyle\frac{\sign(c_Kq(t))}{\sqrt{\Gamma_H}}
 \int_\mR\mathfrak g_u^H(z)\dd\widetilde W_z^{t,h}}
 \label{eq:local-two-sided-fbm}
\end{equation}
has the standard two-sided fBm law, using its continuous Gaussian version.
Let $U$ be the added Brownian path on $(-\infty,0]$, an independent
$\mathcal C_{\loc}((-\infty,0];\mR)$-valued variable.
For $u\ge0$ with $t+hu\le T$,
\begin{align*}
 \overline{\mathcal F_{t+hu}^W\vee\sigma(U)}
 =
 \overline{\sigma(\widetilde W_z^{t,h}:z\le u)}
 =
 \overline{\sigma(\mathbf B_v^{t,h}:v\le u)}.
\end{align*}
The first equality follows from $W_0=0$:
$\widetilde W_{-t/h}^{t,h}=-h^{-1/2}W_t$, so the local past recovers
$W_t$ and then the original Brownian path. The second equality is
\cref{lem:mvn-filtration}. Taking
$\mathbf B^{t,h,c}:=\mathbf B^{t,ch}$ and $U^{t,h,c}:=U$ proves
\textnormal{(FT1)}. 
To establish the uniform $L^2$ approximation, split the increment at $t$
and substitute $s=t+hz$:
\begin{equation}
 h^{-H}\{\Theta_{t+hu}-\Theta_t-F(t+hu)+F(t)\}
 =\int_\mR\mathfrak g_{h,t,u}(z)\dd\widetilde W_z^{t,h},
 \label{eq:volterra-scaled-kernel}
\end{equation}
where $\mathfrak g_{h,t,u}\in L^2(\mR)$ is
\[
 \appdef{\mathfrak g_{h,t,u}}{\mR}{\mR}{z}
 {h^{-\beta}\!\begin{cases}
 \begin{aligned}[t]
 &[K(h(u-z))-K(-hz)]\\[-2pt]
 &\qquad{}\times q(t+hz),
 \end{aligned}&-t/h\le z<0,\\
 K(h(u-z))q(t+hz),&0\le z\le u,\\
 0,&\text{otherwise}.
 \end{cases}}
\]
Endpoint kernel values may be chosen arbitrarily. We prove
\begin{equation}
 \sup_{t\in J}\sup_{0\le u\le R}
 \|\mathfrak g_{h,t,u}-c_Kq(t)\mathfrak g_u^H\|_{L^2(\mR)}
 \xrightarrow[h\to0]{}0.
 \label{eq:volterra-kernel-convergence}
\end{equation}
Fix $A>1+R$ and take $hA<\inf J$. Define the relative kernel error
\[
 \appdef{\zeta_K}{(0,T]}{[0,\infty)}{v}
 {\displaystyle\sup_{0<r\le v}\left|\dfrac{K(r)}{c_Kr^\beta}-1\right|}.
\]
By \eqref{eq:kernel-asymptotics}, $\zeta_K(v)\to0$ as $v\to0$.
For $z\in[-A,u]$ and positive arguments $x\le A+R$,
\begin{align*}
 &|h^{-\beta}K(hx)q(t+hz)-c_Kq(t)x^\beta|\\
 &\quad\le |c_K|x^\beta\bigl(q_+\zeta_K(h(A+R))+\omega_q(hA)\bigr).
\end{align*}
Apply this to $x=u-z$ and, on $z<0$, to $x=-z$.
Consequently,
\begin{align*}
 &\int_{-A}^R|\mathfrak g_{h,t,u}(z)-c_Kq(t)\mathfrak g_u^H(z)|^2\dd z\\
 &\quad\le
 C\bigl(\zeta_K(h(A+R))+\omega_q(hA)\bigr)^2
 \times
 \left(\int_{-A}^u(u-z)^{2\beta}\dd z
       +\int_{-A}^0(-z)^{2\beta}\dd z\right)\\
 &\quad=
 \frac{C}{2H}\bigl(\zeta_K(h(A+R))+\omega_q(hA)\bigr)^2
 \bigl((A+u)^{2H}+A^{2H}\bigr)\\
 &\quad\le
 C_{A,R}\bigl(\zeta_K(h(A+R))+\omega_q(hA)\bigr)^2
 \xrightarrow[h\to0]{}0.
\end{align*}
For the remaining tail $-t/h\le z<-A$, the mean-value theorem and the
derivative bound give
\begin{align*}
 h^{-\beta}|K(h(u-z))-K(-hz)|
 &\le
 h^{-\beta}(hu)\sup_{0\le\theta\le u}|K'(h(-z+\theta))|\\
 &\overset{\eqref{eq:kernel-bounds}}\le
 C_Rh^{-\beta+1+H-3/2}|z|^{H-3/2}
 =C_R|z|^{H-3/2}.
\end{align*}
The limiting kernel has the same tail bound. Therefore
\begin{align*}
 \sup_{t,u}\int_{-\infty}^{-A}
 |\mathfrak g_{h,t,u}-c_Kq(t)\mathfrak g_u^H|^2\dd z
 \le C_R\int_A^\infty z^{2H-3}\dd z
 =\frac{C_R}{2-2H}A^{2H-2}\xrightarrow[A\to\infty]{}0.
\end{align*}
The omitted original past is also negligible:
\[
 \sup_{t\in J,u\le R}\int_{-\infty}^{-t/h}|\mathfrak g_u^H(z)|^2\dd z
 \le C_R\left(\frac{h}{\inf J}\right)^{2-2H}\xrightarrow[h\to0]{}0.
\]
First take $h\to0$ on the compact region, then $A\to\infty$.
This proves \eqref{eq:volterra-kernel-convergence}. All bounds remain
uniform with $h$ replaced by $ch$.
The Wiener isometry and the normalization identity
\[
 c_Kq(t)\int_\mR\mathfrak g_u^H(z)\dd\widetilde W_z^{t,h}
 \overset{\eqref{eq:local-two-sided-fbm}}=\rho(t)\mathbf B_u^{t,h}
\]
handle the stochastic part. The drift satisfies
\[
 \sup_{t\in J,u\le R}\frac{|F(t+hu)-F(t)|}{h^H}
 \le[F]_{\alpha_F}R^{\alpha_F}h^{\alpha_F-H}\xrightarrow[h\to0]{}0.
\]
Since $\inf_t\rho(t)>0$, division by $\rho(t)$ proves
\textnormal{(FT2)}, uniformly also in $c$.
Finally, we verify the global increment bound and continuity.
For $0\le t\le T-h$, decompose
\begin{align*}
 \Theta_{t+h}-\Theta_t
 &=F(t+h)-F(t)
 +\int_0^t\{K(t+h-s)-K(t-s)\}q(s)\dd W_s\\
 &\qquad+\int_t^{t+h}K(t+h-s)q(s)\dd W_s.
\end{align*}
The two stochastic integrals are orthogonal. On $0<u<h$, bound the
kernel difference by $|K(u+h)|+|K(u)|$. On $u\ge h$, use the mean-value
theorem. It\^o's isometry and $\alpha_F>H$ then yield, for small $h\le1$,
\begin{align*}
 \E\left[|\Theta_{t+h}-\Theta_t|^2\right]
 &\le Ch^{2\alpha_F}
 +C\int_0^t|K(u+h)-K(u)|^2\dd u
 +C\int_0^h|K(u)|^2\dd u\\
 &\quad\overset{\eqref{eq:kernel-bounds}}\le
 Ch^{2\alpha_F}+C\int_0^{2h}u^{2H-1}\dd u+Ch^2\int_h^Tu^{2H-3}\dd u+C\int_0^hu^{2H-1}\dd u\\
 &\quad\le Ch^{2\alpha_F}+C\frac{(2h)^{2H}+h^{2H}}{2H}
 +C\frac{h^{2H}-h^2T^{2H-2}}{2-2H}\\
 &\quad\le C(h^{2\alpha_F}+h^{2H})
 \le Ch^{2H}.
\end{align*}
This proves \textnormal{(FT3)}. Gaussian moment equivalence gives
$\E\left[|\Theta_{t+h}-\Theta_t|^p\right]\le C_ph^{pH}$ for every $p\ge2$.
Choose $p>1/H$ and apply Kolmogorov's continuity theorem
\cite[Chapter~4]{kallenberg2021}. Finally,
\[
 \rho(t)^{2/(1+H)}
 =|c_Kq(t)|^{2/(1+H)}\Gamma_H^{1/(1+H)},
\]
so \cref{thm:exact-local} gives \eqref{eq:volterra-exact-equivalent}.
\end{proof}
\section{Bounds and positive quadrature}
\label{app:bounds-quadrature}
\subsection{Single-window bounds}
\label{app:single-window-bounds}
\begin{proof}[Proof of \cref{prop:analytic-bounds}]
\emph{Lower bound.}
For target $\sigma B^H$ on $[0,T]$, \cref{lem:increment} gives
\[
 D_h:=(T-h)m_1\sigma h^H\le hR+2\sqrt{TM}.
\]
With $c:=\lambda b$ and $z:=2\sqrt{TM}$, this bounds the objective by
\begin{align*}
 J_\lambda(\phi;\sigma B^H)
 \ge
   \frac{az^2}{8T}+\frac ch(D_h-z)_+
 \ge\min_{0\le z\le D_h}
       \left\{\frac{az^2}{8T}+\frac ch(D_h-z)\right\}.
\end{align*}
For $z>D_h$ the first expression increases. On $[0,D_h]$, its
unconstrained minimizer is $4Tc/(ah)$, by differentiation. Hence
\[
 J_\lambda(\phi;\sigma B^H)\ge
 \begin{cases}
 aD_h^2/(8T),&D_h\le4Tc/(ah),\\
 cD_h/h-2Tc^2/(ah^2),&D_h>4Tc/(ah).
 \end{cases}
\]
Let $\ell:=(c/(a\sigma))^{1/(1+H)}$ and $h:=q\ell$, $q>0$.
The second case applies whenever
\[
 \frac{ahD_h}{4Tc}
 =\frac{m_1}{4}\left(1-\frac{q\ell}{T}\right)q^{H+1}>1.
\]
In this case, the identity $c=a\sigma\ell^{1+H}$ gives
\begin{align*}
 \frac{J_\lambda(\phi;\sigma B^H)}{Ta\sigma^2\ell^{2H}}
 &\ge\frac{cD_h/h-2Tc^2/(ah^2)}{Ta\sigma^2\ell^{2H}}\\
 &=
 m_1\left(1-\frac{q\ell}{T}\right)q^{H-1}-2q^{-2}\\
 &\xrightarrow[\lambda\to0]{}m_1q^{H-1}-2q^{-2}.
\end{align*}
The derivative of the limiting expression is
\[
 \frac{\dd}{\dd q}\left(m_1q^{H-1}-2q^{-2}\right)
 =q^{-3}\{4-m_1(1-H)q^{H+1}\}.
\]
It changes sign from positive to negative at
\[
 q_*^{H+1}=\frac4{m_1(1-H)},\qquad
 m_1q_*^{H+1}=\frac4{1-H}>4.
\]
Thus $q_*$ is the global maximizer and satisfies the second-case
condition for small $\lambda$. Using
$m_1(1-H)q_*^{H+1}=4$ and taking infima over controls therefore gives
\begin{align*}
 \liminf_{\lambda\to0}
 \frac{V_\lambda(T;a,b,\sigma)}{Ta\sigma^2\ell^{2H}}
 &\ge m_1q_*^{H-1}-2q_*^{-2}\\
 &=\frac{m_1(1+H)}2q_*^{H-1}=C_{\rm low}(H).
\end{align*}
By \cref{thm:cell}, the left-hand side is $\kappa_H$.\\
\emph{Upper bound.}
On the full-past space, use the absolutely continuous rule
\[
 \appdef{\psi^{(h)}}{\mR}{L^2(\Omega)}{t}
 {\frac1h\int_{t-h}^t\mathbf B_s^H\dd s}.
\]
Its derivative and tracking error are
\begin{equation*}
 \dot\psi_t^{(h)}
 =\frac{\mathbf B_t^H-\mathbf B_{t-h}^H}{h},
\end{equation*}
and
\begin{equation*}
 \mathbf B_t^H-\psi_t^{(h)}
 =\frac1h\int_0^h(\mathbf B_t^H-\mathbf B_{t-u}^H)\dd u.
\end{equation*}
Stationary increments and Fubini give
\begin{align*}
 \int_0^h\int_0^h|u-v|^{2H}\dd u\dd v
 &=2\int_0^h\int_0^u(u-v)^{2H}\dd v\dd u
 =\frac{2h^{2H+2}}{(2H+1)(2H+2)}.
\end{align*}
Then,
\begin{align*}
 \E\left[|\mathbf B_t^H-\psi_t^{(h)}|^2\right]
 &=\frac1{h^2}\int_0^h\int_0^hR_H(u,v)\dd u\dd v\\
 &=\frac{h^{2H}}{2H+1}
 -\frac{h^{2H}}{(2H+1)(2H+2)}
 =\frac{h^{2H}}{2(H+1)},
\end{align*}
and
\begin{align*}
 \E\left[|\dot\psi_t^{(h)}|\right]&=m_1h^{H-1}.
\end{align*}
The initial atom and terminal reset each have expectation
$m_1h^H/\sqrt{2(H+1)}$, so their contribution per unit time vanishes.
The long-run cost is
\[
 \appdef{C}{(0,\infty)}{(0,\infty)}{h}
 {\frac{h^{2H}}{4(H+1)}+m_1h^{H-1}},
\]
with derivative
\[
 C'(h)=h^{H-2}\left\{\frac{H}{2(H+1)}h^{H+1}-m_1(1-H)\right\}.
\]
The derivative changes sign once, at
\[
 h_*^{H+1}=\frac{2m_1(1-H)(1+H)}H,
 \qquad
 m_1h_*^{H-1}=\frac{Hh_*^{2H}}{2(1-H)(1+H)}.
\]
Substitution gives
\begin{align*}
 \min_{h>0}C(h)
 =\frac{h_*^{2H}}{4(1-H)}
 =\frac1{4(1-H)}
 \left(\frac{2m_1(1-H)(1+H)}H\right)^{2H/(1+H)}
 =C_{\rm ma}(H).
\end{align*}
Apply \cref{thm:past-robustness} to conclude.
\end{proof}
\subsection{Multi-window bounds}
\label{app:multiwindow-bounds}
Fix an integer $m\ge1$ and define the \emph{affine weight space}
\[
 \cW_m:=\left\{w\in\mR^m:\sum_{i=1}^m w_i=1\right\}.
\]
For
$q=(q_1,\ldots,q_m)\in(0,\infty)^m$ and $w\in\cW_m$, define on the
two-sided fBm space
\begin{equation}
 \appdef{\psi^{w,q}}{\mR}{L^2(\Omega)}{t}
 {\sum_{i=1}^m\frac{w_i}{q_i}\int_{t-q_i}^t\mathbf B_s^H\dd s}.
 \label{eq:multiwindow-policy}
\end{equation}
Define the symmetric matrices
$A^H(q),D^H(q)\in\mR^{m\times m}$ by
\begin{equation}
 \appdef{A^H}{(0,\infty)^m}{\mR^{m\times m}}{q}
 {\left[\frac1{q_iq_j}\int_0^{q_i}\int_0^{q_j}
                  R_H(u,v)\dd u\dd v\right]_{i,j=1}^m},
 \label{eq:A-matrix}
\end{equation}
and
\begin{equation}
 \appdef{D^H}{(0,\infty)^m}{\mR^{m\times m}}{q}
 {\left[\frac{R_H(q_i,q_j)}{q_iq_j}\right]_{i,j=1}^m}.
 \label{eq:D-matrix}
\end{equation}
Both matrices are positive semidefinite because they are covariance
matrices of Gaussian linear functionals. Their quadratic forms are
\[
 \appdef{\mathsf A_H}{\cW_m\times(0,\infty)^m}{[0,\infty)}{(w,q)}
 {w^\top A^H(q)w},
\]
and
\[
 \appdef{\mathsf D_H}{\cW_m\times(0,\infty)^m}{[0,\infty)}{(w,q)}
 {w^\top D^H(q)w}.
\]
For $a,b>0$, integrating the covariance gives
\begin{align}
 \int_0^a\int_0^bR_H(u,v)\dd u\dd v
 &=\frac12\int_0^a\int_0^b
       (u^{2H}+v^{2H}-|u-v|^{2H})\dd v\dd u\notag\\
 &=
   \frac{ba^{2H+1}+ab^{2H+1}}{2(2H+1)}
 -\frac{a^{2H+2}+b^{2H+2}-|a-b|^{2H+2}}
                     {2(2H+1)(2H+2)}.
 \label{eq:integrated-covariance}
\end{align}
\begin{proposition}[multi-window upper bound]
\label{prop:multiwindow}
Define
\[
 \appdef{\mathsf C_H}{\cW_m\times(0,\infty)^m}{[0,\infty)}{(w,q)}
 {\tfrac12\mathsf A_H(w,q)+m_1\sqrt{\mathsf D_H(w,q)}}.
\]
For every $q\in(0,\infty)^m$ and every $w\in\cW_m$,
\begin{equation}
 \kappa_H\le\mathsf C_H(w,q).
 \label{eq:multiwindow-bound}
\end{equation}
\end{proposition}
\begin{proof}
Since $\sum_{i=1}^m w_i=1$,
\begin{equation*}
 \mathbf B_t^H-\psi_t^{w,q}
 =\sum_{i=1}^m \frac{w_i}{q_i}\int_0^{q_i}
       (\mathbf B_t^H-\mathbf B_{t-u}^H)\dd u,
\end{equation*}
and
\begin{equation*}
 \dot\psi_t^{w,q}
 =\sum_{i=1}^m \frac{w_i}{q_i}(\mathbf B_t^H-\mathbf B_{t-q_i}^H).
\end{equation*}
These variables are centered Gaussian, and stationarity of the
increments gives their variances:
\begin{align*}
 \E\left[|\mathbf B_t^H-\psi_t^{w,q}|^2\right]
 &=
 \sum_{i=1}^m\sum_{j=1}^m\frac{w_iw_j}{q_iq_j}
    \int_0^{q_i}\int_0^{q_j}R_H(u,v)\dd u\dd v\\
 &=w^\top A^H(q)w=\mathsf A_H(w,q),
\end{align*}
and
\begin{align*}
 \E\left[|\dot\psi_t^{w,q}|^2\right]
 &=\sum_{i=1}^m\sum_{j=1}^m\frac{w_iw_j}{q_iq_j}R_H(q_i,q_j)
 =w^\top D^H(q)w=\mathsf D_H(w,q).
\end{align*}
The initial atom and terminal reset each have expectation
$m_1\sqrt{\mathsf A_H(w,q)}$. Therefore the full-past reset cost is
\[
 J_L^{\mathbf B^H}(\psi^{w,q})
 =L\left(\frac12\mathsf A_H(w,q)+m_1\sqrt{\mathsf D_H(w,q)}\right)
   +2m_1\sqrt{\mathsf A_H(w,q)}.
\]
Divide by $L$ and use $u_H^{\past}(L)\le J_L^{\mathbf B^H}(\psi^{w,q})$
and \cref{thm:past-robustness}.
\end{proof}
\begin{corollary}[scale-free multi-window bound]
\label{cor:multiwindow-scale-free}
Fix $m\ge1$, $w\in\cW_m$, and $q\in(0,\infty)^m$.  Assume
$\mathsf A_H(w,q)>0$ and $\mathsf D_H(w,q)>0$.  For $c>0$,
\begin{equation}
 \mathsf A_H(w,cq)=c^{2H}\mathsf A_H(w,q),
 \qquad
 \mathsf D_H(w,cq)=c^{2H-2}\mathsf D_H(w,q).
 \label{eq:multiwindow-scaling}
\end{equation}
The unique minimizer of $c\mapsto \mathsf C_H(w,cq)$ is
\begin{equation}
 c_H^*(w,q)
 :=\left(
 \frac{m_1(1-H)\sqrt{\mathsf D_H(w,q)}}{H \mathsf A_H(w,q)}
 \right)^{1/(1+H)}.
 \label{eq:multiwindow-optimal-scale}
\end{equation}
Hence
\begin{equation}
 \kappa_H\le \overline{\mathsf C}_H(w,q),
 \label{eq:scale-free-multiwindow-bound}
\end{equation}
where
\begin{equation}
 \overline{\mathsf C}_H(w,q)
 :=\frac{1+H}{2(1-H)}
 \mathsf A_H(w,q)^{(1-H)/(1+H)}
 \left(
 \frac{m_1(1-H)}{H}\sqrt{\mathsf D_H(w,q)}
 \right)^{2H/(1+H)}.
 \label{eq:scale-free-multiwindow-functional}
\end{equation}
The quantity $\overline{\mathsf C}_H(w,q)$ does not change when $q$ is replaced by
$cq$.  For $m=1$, $w_1=q_1=1$, it is equal to $C_{\rm ma}(H)$.
\end{corollary}
\begin{proof}
From $R_H(cu,cv)=c^{2H}R_H(u,v)$ and the substitutions
$u=c\widetilde u$, $v=c\widetilde v$,
\begin{equation*}
 \begin{aligned}
 A_{ij}^H(cq)
 =\frac{c^{2H+2}}{c^2q_iq_j}
   \int_0^{q_i}\int_0^{q_j}R_H(\widetilde u,\widetilde v)
          \dd\widetilde u\dd\widetilde v
 =c^{2H}A_{ij}^H(q),
 \end{aligned}
\end{equation*}
and
\begin{equation*}
 D_{ij}^H(cq)=c^{2H-2}D_{ij}^H(q).
\end{equation*}
Write $A:=\mathsf A_H(w,q)>0$ and
$D:=\mathsf D_H(w,q)>0$. The scaling identities then give
\begin{align*}
 \mathsf C_H(w,cq)&=\tfrac12Ac^{2H}+m_1\sqrt D\,c^{H-1}.
\end{align*}
Then,
\begin{align*}
 \frac{\dd}{\dd c}\mathsf C_H(w,cq)
 &=c^{H-2}\{HAc^{H+1}-m_1(1-H)\sqrt D\}.
\end{align*}
The bracket is strictly increasing from negative to positive. Its zero
is \eqref{eq:multiwindow-optimal-scale}.
At $c=c_H^*(w,q)$,
\[
 m_1\sqrt D\,c^{H-1}=\frac{H}{1-H}Ac^{2H},
\]
so
\begin{align*}
 \mathsf C_H(w,cq)
 =\frac{1+H}{2(1-H)}Ac^{2H}
 =\frac{1+H}{2(1-H)}
 A^{(1-H)/(1+H)}
 \left(\frac{m_1(1-H)}H\sqrt D\right)^{2H/(1+H)}.
\end{align*}
Apply \cref{prop:multiwindow}. Under $q\mapsto cq$, the total power of
$c$ in this expression is
\[
 \frac{2H(1-H)}{1+H}+\frac{2H(H-1)}{1+H}=0,
\]
which proves scale invariance. For $m=1$ and $w_1=q_1=1$,
$A=1/[2(H+1)]$ and $D=1$. Substitution gives $C_{\rm ma}(H)$.
\end{proof}
\begin{proposition}[optimal weights for fixed windows]
\label{prop:optimal-window-weights}
Fix $q$ and assume $D=D^H(q)$ is positive definite. Put $A=A^H(q)$,
$\mathbf 1=(1,\ldots,1)^\top$, and, for $r>0$,
\[
 \appdef{M}{(0,\infty)}{\mR^{m\times m}}{r}
 {A+\frac{m_1}{r}D}.
\]
Then the best bound within these windows is
\begin{equation}
 \min_{w\in\cW_m}\mathsf C_H(w,q)
 =\min_{r>0}\left\{
     \frac{1}{2\mathbf 1^\top M(r)^{-1}\mathbf 1}
       +\frac{m_1r}{2}\right\}.
 \label{eq:window-weight-reduction}
\end{equation}
The function minimized on the right-hand side is convex in $r$. If $r_*$
is a minimizer, the optimal weight vector is
\begin{equation}
 w_*:=\frac{M(r_*)^{-1}\mathbf 1}
                {\mathbf 1^\top M(r_*)^{-1}\mathbf 1},
 \qquad r_*^2=w_*^\top Dw_*.
 \label{eq:optimal-window-weights}
\end{equation}
\end{proposition}
\begin{proof}
For $z>0$ and $r>0$, completing the square gives
\[
 \frac{z}{2r}+\frac r2-\sqrt z
 =\frac{(r-\sqrt z)^2}{2r}\ge0,
\]
with equality exactly when $r^2=z$. Define
\[
 \appdef{F}{\cW_m\times(0,\infty)}{[0,\infty)}{(w,r)}
 {\tfrac12w^\top M(r)w+\tfrac12m_1r}.
\]
This identity gives $\mathsf C_H(w,q)=\min_{r>0}F(w,r)$.
For fixed $r$, $M(r)$ is positive definite. A Lagrange multiplier
$\zeta\in\mR$ for $\mathbf1^\top w=1$ gives
\[
 M(r)w=\zeta\mathbf1,\qquad
 w(r)=\frac{M(r)^{-1}\mathbf1}{\mathbf1^\top M(r)^{-1}\mathbf1}.
\]
Thus
\[
 \frac12w(r)^\top M(r)w(r)
 =\frac\zeta2\mathbf1^\top w(r)
 =\frac1{2\mathbf1^\top M(r)^{-1}\mathbf1},
\]
which proves \eqref{eq:window-weight-reduction} after minimizing over $r$. 
For directions $z\in\mR^m$ and $s\in\mR$, the second directional
variation of the quadratic perspective is
\[
 \mathrm d^2\left(\frac{w^\top Dw}{2r}\right)[(z,s),(z,s)]
 =\frac1r\left(z-\frac srw\right)^\top
 D\left(z-\frac srw\right)\ge0.
\]
Hence $F$ is jointly convex, and partial minimization on
$\mathbf1^\top w=1$ preserves convexity. Since $D$ is positive definite and
$\|w\|^2\ge1/m$ on $\cW_m$,
\begin{align*}
 F(w,r)
 &\ge\frac{m_1}{2r}w^\top Dw+\frac{m_1r}{2}\\
 &\ge\frac{m_1\lambda_{\min}(D)}{2mr}+\frac{m_1r}{2}
 \xrightarrow[r\to0\ \text{or}\ r\to\infty]{}\infty,
\end{align*}
where $\lambda_{\min}(D)>0$ is the smallest eigenvalue.
A minimizer $r_*$ exists. Equality in the scalar identity gives
$r_*^2=w_*^\top Dw_*$, with $w_*=w(r_*)$.
\end{proof}
\subsection{Positive quadrature}
\label{app:positive-quadrature}
\begin{proof}[Proof of \cref{cor:factor-diagonal}]
\emph{OU increment kernels.}
For $x>0$ and $t\ge0$, define
\[
 \appdef{g_{t,x}}{\mR}{\mR}{s}
 {e^{-x(t-s)}\ind_{\{s<t\}}-e^{xs}\ind_{\{s<0\}}}.
\]
For $x,y>0$, splitting at $s=0$ gives
\begin{align*}
 \langle g_{t,x},g_{t,y}\rangle
 &=
 \int_0^t e^{-(x+y)r}\dd r
 +(e^{-xt}-1)(e^{-yt}-1)\int_0^\infty e^{-(x+y)r}\dd r\\
 &=\frac{1-e^{-(x+y)t}+(e^{-xt}-1)(e^{-yt}-1)}{x+y}\\
 &=\frac{2-e^{-xt}-e^{-yt}}{x+y}.
\end{align*}
Setting $y=x$ and using $1-e^{-v}\le v\wedge1$ for $v\ge0$ yields
\begin{equation}
 \|g_{t,x}\|_{L^2(\mR)}^2
 =\frac{1-e^{-xt}}x
 \le t\wedge x^{-1}.
 \label{eq:OU-increment-kernel-norm}
\end{equation}
The map $x\mapsto g_{t,x}$ is differentiable in $L^2(\mR)$.
Differentiating the preceding inner product in both variables gives
\begin{align*}
 \|\partial_xg_{t,x}\|_2^2
 &=\left.\partial_x\partial_y
       \frac{2-e^{-xt}-e^{-yt}}{x+y}\right|_{y=x}\\
 &=-\frac{te^{-xt}}{2x^2}+\frac{1-e^{-xt}}{2x^3}.
\end{align*}
Consequently,
\begin{equation}
 \|\partial_xg_{t,x}\|_{L^2(\mR)}^2
 =\frac{1-(1+xt)e^{-xt}}{2x^3},
\end{equation}
and
\begin{equation}
 x\|\partial_xg_{t,x}\|_2\le2^{-1/2}\|g_{t,x}\|_2,
 \label{eq:OU-increment-kernel-derivative}
\end{equation}
because $0\le1-(1+v)e^{-v}\le1-e^{-v}$ for $v\ge0$.\\
\emph{Bochner representation and truncation.}
Let
$\mu_H(\dd x):=\sqrt{2H}\,x^{\beta-1}\dd x/\Gamma(\beta)$.
For $t>0$, splitting at $x=t^{-1}$ gives
\begin{align*}
 \int_0^\infty\|g_{t,x}\|_2\mu_H(\dd x)
 &\le C_H\left\{\sqrt t\int_0^{1/t}x^{\beta-1}\dd x
                  +\int_{1/t}^\infty x^{\beta-3/2}\dd x\right\}\\
 &=C_H\left(\frac1\beta+\frac1H\right)t^H<\infty.
\end{align*}
Hence the Bochner integral in $L^2(\mR)$ exists, and stochastic Fubini
identifies the fBm increment kernel as
\[
 \frac1{A_H}\int_0^\infty g_{t,x}\mu_H(\dd x).
\]
For $t\le L$, the low- and high-frequency errors satisfy
\begin{align*}
 \left\|\int_0^\varepsilon g_{t,x}\mu_H(\dd x)\right\|_2
 &\le C_H\sqrt L\int_0^\varepsilon x^{\beta-1}\dd x
 =\frac{C_H}{\beta}\sqrt L\,\varepsilon^\beta,
\end{align*}
and
\begin{align*}
 \left\|\int_R^\infty g_{t,x}\mu_H(\dd x)\right\|_2
 &\le C_H\int_R^\infty x^{\beta-3/2}\dd x
 =\frac{C_H}{H}R^{-H}.
\end{align*}
\emph{Positive quadrature error.}
For $x\in[x_j,x_{j+1}]$, apply the fundamental theorem of calculus
and the derivative estimate. Since $x/x_j\le1+\eta\le2$ and
$\log(x/x_j)\le\eta$, this gives
\begin{align*}
 \|g_{t,x}-g_{t,x_j}\|_2
 &\le
   \int_{x_j}^x\|\partial_yg_{t,y}\|_2\dd y\\
 &\overset{\eqref{eq:OU-increment-kernel-derivative}}\le
   C(\sqrt t\wedge x^{-1/2})\int_{x_j}^x\frac{\dd y}{y}\\
 &\le
   C\eta(\sqrt t\wedge x^{-1/2}).
\end{align*}
The weights are the positive masses
\[
 \mu_H([x_j,x_{j+1}])
 =\frac{\sqrt{2H}}{\beta\Gamma(\beta)}
       (x_{j+1}^\beta-x_j^\beta)=w_j.
\]
Thus
\begin{align*}
 \left\|\int_\varepsilon^R g_{t,x}\mu_H(\dd x)
                    -\sum_{j=0}^{N-1}w_jg_{t,x_j}\right\|_2
 &\le\sum_{j=0}^{N-1}\int_{x_j}^{x_{j+1}}
               \|g_{t,x}-g_{t,x_j}\|_2\mu_H(\dd x)\\
 &\le C_H\eta t^H\le C_H\eta L^H.
\end{align*}
Adding the three errors, dividing by $A_H$, and applying the Wiener
isometry proves \eqref{eq:constructive-kernel-bound}.\\
\emph{Resolution and value rate.}
For the parameters in \eqref{eq:constructive-resolution},
\[
 \sqrt L\,\varepsilon_L^\beta=R_L^{-H}=\eta_LL^H=L^{-\gamma}.
\]
A geometric partition, with a possibly shorter last interval, has
\begin{align*}
 N_L&\le1+\frac{\log(R_L/\varepsilon_L)}{\log(1+\eta_L)}\\
 &\le
 1+2L^{H+\gamma}
 \left(\frac\gamma H+\frac{\gamma+1/2}{\beta}\right)\log L\\
 &=O_{H,\gamma}(L^{H+\gamma}\log L),
\end{align*}
using $\log(1+\eta)\ge\eta/2$ for $0<\eta\le1$.
The path error bound gives
\[
 A_{N_L,L}/L=O(L^{-2\gamma}),\qquad
 B_{N_L,L}=O(L^{-\gamma}).
\]
Substitution in \eqref{eq:factor-quantitative-value-bound} yields
\begin{align*}
 \left|\frac{u_{N_L}(L)}L-\kappa_H\right|
 \overset{\eqref{eq:factor-quantitative-value-bound}}\le
 C_{H,\gamma}(L^{-1}+L^{-\gamma-1}+L^{-2\gamma}+L^{-\gamma})\le C_{H,\gamma}(L^{-\gamma}+L^{-1}).
\end{align*}
\end{proof}
\section{A quadratic-hedging corollary}
\label{sec:hedging}
Let $S$ be a continuous square-integrable $\mathbf F$-martingale with
\begin{equation}
 \dd\langle S\rangle_t=a_t\dd t,
 \label{eq:martingale-qv}
\end{equation}
where $a,b$ satisfy \cref{def:tracking-value}.
For $\xi\in L^2(\mathcal F_T)$, put $\xi_0=\E[\xi\mid\mathcal F_0]$.
The \emph{Galtchouk--Kunita--Watanabe decomposition} is
\begin{equation}
 \xi=\xi_0+\int_0^T\Theta_t\dd S_t+L_T,
 \qquad L_0=0,\qquad\langle L,S\rangle=0,
 \label{eq:GKW}
\end{equation}
where $L$ is square-integrable and $\Theta$ is a predictable integrand
with $\E[\int_0^T\Theta_t^2\dd\langle S\rangle_t]<\infty$, as in
\cite{protter2005}. Assume this integrand has a continuous adapted
version, which we fix.
For initial cash $c\in L^2(\mathcal F_0)$ and
$\phi\in\cA_T(\Theta)$, define
\begin{equation}
 \appdef{\mathcal Q_\lambda^c}{\cA_T(\Theta)}{[0,\infty)}{\phi}
 {\frac12\E\left[\left(\xi-c-\int_0^T\phi_{t-}\dd S_t\right)^2\right]
       +\lambda\E\left[\int_{[0,T]}b_t\dd|\phi|_t\right]}.
 \label{eq:quadratic-hedging-objective}
\end{equation}
\begin{proposition}[penalized mean-square identity]
\label{thm:gkw-identity}
Under these assumptions,
\begin{equation}
 \inf_{\phi\in\cA_T}\mathcal Q_\lambda^c(\phi)
 -\frac12\E\bigl[(\xi_0-c)^2+L_T^2\bigr]=V_\lambda(\Theta).
 \label{eq:gkw-value-identity}
\end{equation}
If $a,b$ are continuous and $\Theta$ satisfies
\cref{ass:strong-tangent} in the same observation filtration, then
\begin{equation}
 \begin{aligned}
 \inf_{\phi\in\cA_T}\mathcal Q_\lambda^c(\phi)
 -\frac12\E\bigl[(\xi_0-c)^2+L_T^2\bigr]\underset{\lambda\to0}{\sim}
   \kappa_H\lambda^{2H/(1+H)}\int_0^Tg_H(t)\dd t.
 \end{aligned}
 \label{eq:gkw-sharp-asymptotic}
\end{equation}
In particular, the second statement applies to the Volterra class of
\cref{prop:volterra-strong-tangent} when its driving Brownian filtration
is the observation filtration.
\end{proposition}
\begin{proof}
Using $(x+y)^2\le2x^2+2y^2$, the admissibility and integrability
assumptions give
\begin{align*}
 \E\left[\int_0^T\phi_t^2\dd\langle S\rangle_t\right]
 \le
 2\E\left[\int_0^T(\phi_t-\Theta_t)^2a_t\dd t\right]
 +2\E\left[\int_0^T\Theta_t^2\dd\langle S\rangle_t\right]
 <\infty.
\end{align*}
A c\`adl\`ag finite-variation path has countably many jumps. Since
$\dd\langle S\rangle_t=a_t\dd t$, replacing $\phi_t$ by its predictable
left-limit process $\phi_{t-}$ leaves this integral unchanged.
Define the square-integrable martingale
\[
 N_t:=\int_0^t(\Theta_r-\phi_{r-})\dd S_r,\qquad N_0=0.
\]
The hedging error decomposes as
\[
 \xi-c-\int_0^T\phi_{t-}\dd S_t=(\xi_0-c)+L_T+N_T.
\]
These three terms are orthogonal in $L^2(\Omega)$: the martingale
increments have zero conditional mean given $\mathcal F_0$, while
$\langle L,S\rangle=0$ gives the last equality below.
First,
\begin{equation*}
 \E[(\xi_0-c)N_T]
 =\E[(\xi_0-c)\E[N_T\mid\mathcal F_0]]=0.
\end{equation*}
Similarly,
\begin{equation*}
 \E[(\xi_0-c)L_T]
 =\E[(\xi_0-c)\E[L_T\mid\mathcal F_0]]=0.
\end{equation*}
Finally,
\begin{equation*}
 \E[L_TN_T]
 =\E\left[\langle L,N\rangle_T\right]
 =\E\left[\int_0^T(\Theta_t-\phi_{t-})\dd\langle L,S\rangle_t\right]
 =0.
\end{equation*}
It\^o's isometry gives
\begin{align*}
 \E[N_T^2]
 =
 \E\left[\int_0^T(\Theta_t-\phi_{t-})^2\dd\langle S\rangle_t\right]
 =
 \E\left[\int_0^Ta_t(\Theta_t-\phi_t)^2\dd t\right].
\end{align*}
Consequently,
\begin{align*}
 \mathcal Q_\lambda^c(\phi)
 &=\frac12\E[(\xi_0-c)^2+L_T^2]+\frac12\E[N_T^2]
       +\lambda\E\left[\int_{[0,T]}b_t\dd|\phi|_t\right]\\
 &=\frac12\E[(\xi_0-c)^2+L_T^2]+J_\lambda(\phi;\Theta).
\end{align*}
Take infima to obtain \eqref{eq:gkw-value-identity}, then apply
\cref{thm:exact-local} for the asymptotic formula.
\end{proof}
\begin{corollary}[initial hedge position]
\label{cor:fixed-initial-hedge}
If $\E[|\Theta_0|]<\infty$, replacing the pre-initial hedge position
$\Theta_0$ by a fixed $x_0\in\mR$ changes the value by at most
$\lambda b_+\E[|x_0-\Theta_0|]$.
Thus \eqref{eq:gkw-sharp-asymptotic} is unchanged.
\end{corollary}
\begin{proof}
Apply \cref{prop:initial-inventory} to \eqref{eq:gkw-value-identity}.
\end{proof}
\section*{AI Assistance Declaration}
During the preparation of this manuscript, the author used AI tools for language editing and manuscript organization. The author takes responsibility for the mathematical statements, proofs, computations, references, and final content.
\begingroup
\fontsize{9}{10}\selectfont

\endgroup
\end{document}